\documentclass[a4paper, 12pt]{article}

\usepackage{theorem}
\usepackage{amsmath,amssymb,amscd,latexsym,eufrak,color,setspace,authblk}
\usepackage[tikz]{mdframed}
\usepackage[all]{xy}
\usepackage{multirow}

\newtheorem{thm}{Theorem}[section] 
\newtheorem{lem}[thm]{Lemma}

\newtheorem{cor}[thm]{Corollary}

\newcommand{\QED}{{\unskip\nobreak\hfil\penalty50\quad\null\nobreak\hfil
{Q.E.D.}\parfillskip0pt\finalhyphendemerits0\par\medskip}}
\newcommand{\Proof}{\noindent{\bf Proof.}\quad}

\newcommand{\bm}[1]{\mbox{\boldmath$#1$}}

\newcommand{\ds}{\displaystyle}

\newcommand{\Exp}{{\rm Exp}}
\newcommand{\id}{{\rm id}}
\newcommand{\Mat}{{\rm Mat}}
\newcommand{\mf}{\mathfrak}
\newcommand{\Nilp}{{\rm Nilp}}

\newcommand{\ord}{{\rm ord}}

\newcommand{\Rank}{{\rm Rank}\,}
\newcommand{\Rep}{{\rm Rep}}
\newcommand{\Seq}{{\rm Seq}}

\newcommand{\Tr}{{\rm Tr}}
\newcommand{\transpose}{\text{{}$^t$}}
\newcommand{\Ker}{{\rm Ker}}
\newcommand{\Span}{{\rm Span}}
\newcommand{\wt}{\widetilde}

\newcommand{\cN}{{\cal N}}
\newcommand{\cP}{{\cal P}}
\newcommand{\cQ}{{\cal Q}}

\newcommand{\cU}{{\cal U}}
\newcommand{\V}{\mathbb{V}}
\newcommand{\E}{{\cal E}}
\newcommand{\F}{F}
\newcommand{\G}{\mathbb{G}}
\newcommand{\X}{\mathbb{X}}
\newcommand{\Y}{\mathbb{Y}}
\newcommand{\Z}{\mathbb{Z}}
\newcommand{\W}{\mathbb{W}}

\newcommand{\diag}{{\rm diag}}

\definecolor{skyblue}{cmyk}{ 0.2, 0, .1, 0} 
\definecolor{OliveGreen}{cmyk}{ 0.64, 0, 0.95, 0.40} 
\definecolor{GreenYellow}{cmyk}{ 0.15, 0.00, 0.69, 0.00} 
\definecolor{Apricot}{rgb}{0.98, 0.81, 0.69}
\definecolor{Orange}{cmyk}{0, 0.61, 0.87, 0}
\definecolor{Tan}{cmyk}{0.14, 0.42, 0.56, 0.00}
\definecolor{Bisque}{rgb}{1.0, 0.89, 0.77}
\definecolor{Almond}{rgb}{0.94, 0.87, 0.8}

\allowdisplaybreaks[4]

\begin{document}

\title{\bf Exponential matrices}
\author{Ryuji Tanimoto} 
\date{}
\maketitle

\begin{abstract} 
In this article, we introduce a notion of an exponential matrix, which is a polynomial matrix with exponential properties, 
and a notion of an equivalence relation of two exponential matrices, and then we initiate to study 
classifying exponential matrices in positive characteristic, up to equivalence. 
We classify exponential matrices of Heisenberg groups in positive characteristic, up to equivalence. 
We also classify exponential matrices of size four-by-four in positive characteristic, up to equivalence. 
From these classifications, we obtain a classification of modular representations of elementary 
abelian $p$-groups into Heisenberg groups, up to equivalence, and a classification of four-dimensional modular 
representations of elementary abelian $p$-groups, up to equivalence.  

\renewcommand{\thefootnote}{\fnsymbol{footnote}}
\end{abstract}


\section{Introduction}

In this article, we introduce a notion of an exponential matrix, which is a polynomial matrix with exponential properties, 
and a notion of an equivalence relation of two exponential matrices, and then we initiate to study 
classifying exponential matrices in positive characteristic, up to equivalence. 

In characteristic zero, classifying exponential matrices up to equivalence can be rephrased as 
classifying nilpotent matrices up to equivalence, and the classification of nilpotent matrices is a beautiful result of Camille Jordan. 

On the other hand, in positive characteristic, classifying exponential matrices is complicated, even if we consider up to equivalence. 
In fact, classifying exponential matrices up to equivalence 
is closely related to classifying representations of the additive group scheme up to equivalence. 
We only have partial results concerning representations of the additive group scheme in positive characteristic 
(see \cite{Fauntleroy}, \cite{Tanimoto 2008}, \cite{Tanimoto 2013}, \cite{Tanimoto 2018}).

In this article, we keep in mind the following guidelines (1) and (2) toward classifying exponential matrices in positive characteristic, 
up to equivalence. 
Let $k$ be a field of positive characteristic $p > 0$, and let $k[T]$ be a polynomial ring in one variable over $k$. 
\begin{enumerate}
\item[\rm (1)] We do not figure out exponential matrices of a general linear group $GL(n, k[T])$, at once 
(since the figuring out is difficult). We take a step-by-step approach. 
We firstly select appropriate subgroups and subsets of $GL(n, k[T])$, 
and secondly describe exponential matrices belonging to the subgroups and the subsets, up to equivalence.
  
\item[\rm (2)] For two described exponential matrices of same size, we consider whether or not the two exponential matrices 
are equivalent. 
\end{enumerate} 

Along the above guidelines, we write this article as follows:

In Section 1, we prepare basics of exponential matrices. 
Any exponential matrix has an exponential expression (see \cite{Miyanishi 1978}, 
and also see Lemma 1.5 in this article). 
Using this expression, we can prove that any exponential matrix is triangulable (see Lemma 1.8). 
We then define subgroups $\mf{U}_n$ and $\mf{U}_{[d_1, \ldots, d_t]}$ of $GL(n, k[T])$ (see Subsubsection 1.3.2). 
These subgroups are appropriate along the above guideline (1).  
For considering these subgroups, we introduce 
the other transpose ${^\tau}A$ of a square matrix $A$ (see Subsubsection 1.3.1), 
and notations $\leadsto$ and $( \quad)^E$ (see Subsection 1.4). 
Using the other transpose, we can reduce the amount of classifying exponential matrices by half, 
and using $\leadsto$ and $( \quad)^E$, we express a link to a set of described exponential matrices 
from a set of exponential matrices belonging to any one of appropriate subgroups and subsets of $GL(n, k[T])$. 
Let $\E(n, k[T])$ be the set of all exponential matrices of $GL(n, k[T])$. 
The starting point of classifying exponential matrices of $GL(n, k[T])$ comes from the following  expression (see Lemma 1.12): 
\[
 \E(n, k[T]) \leadsto \mf{U}_n^E  =  \bigcup_{\text{$[ d_1, \ldots, d_t ]$ is an ordered partition of $n$}}  \mf{U}_{[d_1, \ldots, d_t]}^E. 
\]
So, classifying exponential matrices of $GL(n, k[T])$ reduces to classifying exponential matrices of $\mf{U}_n^E$. 
Based on the above expression of $\mf{U}_n^E$,  
we especially study, in Subsection 1.5, equivalent forms of exponential matrices belonging to subsets 
$\mf{U}_{[n]}^E$, $\mf{U}_{[n, 1]}^E$ and $\mf{U}_{[1, n]}^E$, 
and prove 
\[
\mf{U}_{[n]}^E \leadsto \mf{J}_{[n]}^E, \qquad 
\mf{U}_{[n, 1]}^E \leadsto (\mf{J}_{[n, 1]}^0)^E \cup (\mf{J}_{[n, 1]}^1)^E, \qquad 
\mf{U}_{[1, n]}^E \leadsto (\mf{J}_{[n, 1]}^0)^E \cup (\mf{J}_{[1, n]}^1)^E 
\]  
provided that $1 \leq n \leq p$ (see Corollaries 1.15, 1.22 and 1.26). 
Along the above guideline (2), we consider equivalence relations of two exponential matrices of 
$\mf{J}_{[n]}^E$, $(\mf{J}_{[n, 1]}^0)^E$, $(\mf{J}_{[n, 1]}^1)^E$ and $(\mf{J}_{[1, n]}^1)^E$, respectively (see Subsection 1.8). 
We introduce a notion of mutually $GL(n, k)$-disjoint sets, and consider 
non-equivalence of two described exponential matrices of same size (see Subsection 1.7).

Section 1 implies that, for any $1 \leq n \leq 3$, we can classify exponential matrices of $GL(n, k[T])$, up to equivalence. 
However, we cannot classify exponential matrices of $GL(4, k[T])$, up to equivalence. 
To get rid of the complexity, recall again the above guideline (1). 
We noticed that we should select a Heisenberg subgroup of $GL(4, k[T])$. 

In Section 2, we consider polynomial matrices of a Heisenberg group $\mf{H}(m + 2, k[T])$ of $GL(m + 2, k[T])$. 
The point of classifying polynomial matrices of $\mf{H}(m + 2, k[T])$ appears in the following expression (see Theorem 2.1): 
\[
\mf{H}(m + 2, k[T]) \leadsto \daleth_{m +2} \cup  \bigcup_{(\ell, r_1, r_2) \in \Omega_m} {^\ell \mf{H}}^{r_1}_{r_2}. 
\] 
We give a classification of polynomial matrices of $\daleth_{m +2}$ (see Theorem 2.5), 
and consider equivalence relations of two polynomial matrices of ${^\ell}\mf{H}^{r_1}_{r_2}$ and $\daleth_{m + 2}$ 
(see Subsection 2.4).

In Section 3, we consider exponential matrices of a Heisenberg group. 
We can describe any exponential matrix of $({^\ell}\mf{H}^{r_1}_{r_2})^E$ (see Theorem 3.3). 
Then we can consider equivalence relations of two exponential matrices of $({^\ell}\mf{H}^{r_1}_{r_2})^E$ (see Theorems 3.17 and 3.23). 
We classify exponential matrices of $(\daleth_{m + 2})^E$ up to equivalence (see Theorem 3.13), 
and consider equivalence relations of two exponential matrices of $(\daleth_{m + 2})^E$ (see Theorem 3.24)

In Section 4, we give a classification of exponential matrices of size four-by-four, up to equivalence. 
Then we know that the classification varies according to $p= 2$, $p = 3$ and $p \geq 5$. 

In Section 5, we consider modular representations of elementary abelian $p$-groups. 
We obtain a classification of modular representations of elementary 
abelian $p$-groups into Heisenberg groups, up to equivalence (see Theorem 5.7),  
and a classification of four-dimensional modular representations of elementary abelian $p$-groups, up to equivalence 
(see Subsection 5.4).  

It is an open problem to classify modular representations of elementary abelian $p$-groups, up to equivalence. 
Even so, we at least know the following: 
There are exactly $p$ inequivalent indecomposable modular representations of $\Z/p\Z$ (see for example 
\cite[Page 105]{Campbell-Wehlau}). 
Ba\v{s}ev \cite{Basev} classified indecomposable modular representations of $\Z/2\Z \times \Z/2\Z$ 
over an algebraically closed field of characteristic two. 
The representation type of any elementary abelian $p$-group $E$ except 
for $\Z/p\Z$ and $\Z/2\Z \times \Z/2\Z$ is wild (see \cite{Bondarenko-Drozd, Ringel}). 
Campbell, Shank and Wehlau \cite{Campbell-Shank-Wehlau} describe parametrizations for modular representations 
of elementary abelian $p$-groups with representation spaces in dimension two and in dimension three, respectively.  
For background material on the theory of modular representations of elementary abelian $p$-groups, see \cite{Benson}. 
\\



{\it Preliminary notations and preliminary defnitions.}  
For a commutative ring $R$ with unity, 
we denote by $\Mat_{m, n}(R)$ the module of all $m \times n$ matrices whose entries belong to $R$. 
We denote by $O_{\Mat_{m, n}(R)}$ the zero matrix of $\Mat_{m, n}(R)$. 
We write $\Mat_{n, n}(R)$ as $\Mat(n, R)$. The module $\Mat(n, R)$ naturally becomes a ring with unity. 
We denote by $I_{\Mat(n, R)}$ the identity matrix of $\Mat(n, R)$. 
The zero matrix $O_{\Mat_{m, n}(R)}$ is frequently referred as $O_{m, n}$ or $O$, 
and the zero matrix $O_{\Mat(n, R)}$ is frequently referred as $O_n$ or $O$, 
and the identity matrix $I_{\Mat(n, R)}$ is frequently referred as $I_n$ or $I$. 
For matrices $A_i \in \Mat(n_i, R)$ $(1 \leq i \leq r)$, we denote by 
$\bigoplus_{i = 1}^r A_i$ 
denote the direct sum of the matrices $A_i$ $(1 \leq i \leq r)$. 
In particular if the $A_i$'s are the same matrix $A$, we also use the notation $A^{\oplus r}$ in place of 
$\bigoplus_{i = 1}^r A$.
For any square matrix $A$ of $\Mat(n, R)$, we denote by $\det(A)$ the determinant of $A$ and 
by $\Tr(A)$ the trace of $A$. 
For $a_1, \ldots, a_n \in R$, we denote by $\diag(a_1, \ldots, a_n)$ the diagonal matrix of $\Mat(n, R)$ whose $(i, i)$-th entries  
are $a_i$ for all $1 \leq i \leq n$.

We say that a matrix $A \in \Mat(n, R)$ is {\it regular} if there exists a matrix $X \in \Mat(n, R)$ such that $AX = XA = I_n$. 
We denote by $GL(n, R)$ the group of all regular matrices of $\Mat(n, R)$. 
We denote by $SL(n, R)$ the subgroup of all regular matrices of $GL(n, R)$ whose determinants are $1$. 
We denote by $B(n, R)$ the subgroup of all regular upper triangular matrices of $GL(n, R)$.  
We say that a matrix $N$ of $\Mat(n, R)$ is {\it nilpotent} if $N^\ell = O$ for some $\ell \geq 1$.  
For any $A \in \Mat(n, R)$, we denote by $\det(A)$ the determinant of $A$. 
A nilpotent matrix $N$ of $\Mat(n, R)$ is said to be {\it $\ell$-nilpotent} if $N^\ell = O$. 
Clearly, $1$-nilpotent matrix of $\Nilp(n, R)$ is a zero matrix of $\Mat(n, R)$. 


Let $R[T]$ be a polynomial ring in one variable over $R$. 
We say that an element $A(T)$ of $\Mat(n, R[T])$ is a {\it polynomial matrix} of size $n \times n$ over $R$. 
For all polynomial matrices $A(T), B(T) \in \Mat(n, R[T])$, we say that $A(T)$ and $B(T)$ are {\it equivalent} 
if there exists a regular matrix $P \in GL(n, R)$ such that $P^{-1} A(T) P = B(T)$.

Let $k$ be a field and let $V$ be the column space $V$ in dimension $n$ over $k$. 
For all $1 \leq i \leq n$,  we denote by $\bm{e}_i $ the column vector of $V$ whose $i$-th entry is $1$ and the other 
entries are zeroes. 
For any subset $W$ of $V$, we denote by $\Span_k W$ the subspace spanned by $W$.

We denote by $\Z$ the ring of all integers. 
We denote by $\Z/n\Z$ the cyclic group of order $n$. 
For a prime number $p$ and an integer $r \geq 1$, 
a finite abelian $p$-group $(\Z/p\Z)^r$ is said to be an {\it elementary abelian $p$-group of rank $r$}. 
For a field $k$ of characteristic $p \geq 0$ and for a finite group $G$, we say that 
a representation $\rho : G \to GL(n, k)$ of $G$ is {\it modular} if $p > 0$ and $p$ divides the order of $G$.

\section{Basics of exponential matrices}

Let $k$ be a field of characteristic $p \geq 0$. For any polynomial matrix $A(T) \in \Mat(n, k[T])$, 
we say that $A(T)$ is an {\it exponential matrix} if $A(T)$ 
satisfies the following conditions (1) and (2): 
\begin{enumerate}
\item[\rm (1)] $A(0) = I_{\Mat(n, k)}$.   
\item[\rm (2)] $A(T) A(T') = A(T + T')$, where $T, T'$ are indeterminates over $k$. 
\end{enumerate} 

The above conditions (1) and (2) imply that $A(T) \in GL(n, k[T])$ and  $A(T)^{-1} = A(-T)$. 
We mention here that $A(T) \in SL(n, k[T])$ (see Lemmas 1.8 and 1.9). 

Let $\E(n, k[T])$ be the set of all exponential matrices of size $n \times n$. 
Clearly, for all $P \in GL(n, k)$ and $A(T) \in \E(n, k[T])$, we have $P^{-1} A(T) P \in \E(n, k[T])$. 
For $A(T), B(T) \in \E(n, k[T])$, we say that $A(T)$ and $B(T)$ are {\it equivalent} 
if there exists a regular matrix $P \in GL(n, k)$ such that $P^{-1} A(T) P = B(T)$. 
If exponential matrices $A(T)$ and $B(T)$ are equivalent, we write $A(T) \sim B(T)$.

\subsection{A necessary and sufficient condition for a polynomial matrix to be an exponential matrix}

We note the following necessary and sufficient condition for a given polynomial matrix to be an exponential matrix. 
The proof of the following lemma is directly obtained from the definition of an exponential matrix. 

\begin{lem} 
Let $A(T) = \sum_{i \geq 0} T^i N_i \in \Mat(n, k[T])$ be a polynomial matrix, where 
$N_i \in \Mat(n, k)$ for all $i \geq 0$. Then $A(T)$ is an exponential matrix if and only if 
$A(T)$ satisfies the following conditions {\rm (1)} and {\rm (2)}: 
\begin{enumerate}
\item[\rm (1)] $N_0 = I_n$. 
\item[\rm (2)] $\ds N_i N_j = \binom{i + j}{i} N_{i + j}$ for all $i, j \geq 0$. 
\end{enumerate} 
\end{lem}

\subsection{An expression of an exponential matrix}

In Subsection 1.2, we give an expression of an exponential matrix (see Lemmas 1.2 if $p = 0$, and see Lemma 1.5 if $p  >0$). 
The following Lemmas 1.2, 1.3 and 1.5 are obtained in more general situation 
(see for instance \cite{Miyanishi 1978}). Here, we write the lemmas and their proofs 
in terms of Matrix Theory.

If $p = 0$, we have the following expression of an exponential matrix. 

\begin{lem}
Let $k$ be a field of characteristic zero and let $A(T) \in \Mat(n, k[T])$ be an exponential matrix. 
Then there exists a unique nilpotent matrix $N \in \Mat(n, k)$ such that 
\[
 A(T) = \sum_{i \geq 0} \frac{T^i}{i!} N^i . 
\]
\end{lem}

\Proof Write $A(T) = \sum_{i \geq 0} T^i N_i$, where $N_i \in \Mat(n, k)$ for all $i \geq 0$. 
We know from Lemma 1.1 that $N_i = N_1^i/i!$ for all $i \geq 1$, which implies $A(T) = \sum_{i \geq 0} T^i  (N_1^i / i !)$. 
Since $A(T)$ is a polynomial matrix, the matrix $N_1$ is a nilpotent matrix. 
\QED 

The above Lemma 1.2 asserts that 
a classification of exponential matrices belonging to $\Mat(n, k[T])$ up to equivalence  can be given 
by the Jordan canonical forms of nilpotent matrices belonging to $\Mat(n, k)$. 
So, the classification of exponential matrices up to equivalence is settled if the characteristic of $k$ is zero.

From now on until the end of Section 5, we assume unless otherwise specified that the characteristic $p$ of $k$ is positive.

\begin{lem} 
Let $A(T) = \sum_{i \geq 0} T^i N_i \in \Mat(n, k[T])$ be an exponential matrix, where 
$N_i \in \Mat(n, k)$ for all $i \geq 0$. Then the following assertions {\rm (1), (2), (3)} hold true: 
\begin{enumerate}
\item[\rm (1)] $N_i N_j = N_j N_i$ for all $i, j \geq 0$. 
\item[\rm (2)] $N_i^p = O$ for all $i \geq 1$. 
\item[\rm (3)] Let $i $ be a non-negative integer and let $i = i_0 + i_1p + \cdots + i_r p^r$ be the $p$-adic expansion of $i$. 
Then we have
\[
N_{i_0 + i_1p + \cdots + i_r p^r}
 = N_{i_0} N_{i_1 p} \cdots N_{i_r p^r} 
 = \frac{N_1^{i_0}}{i_0!} \cdot \frac{N_p^{i_1}}{i_1!} \cdots \frac{N_{p^r}^{i_r}}{i_r!} . 
\]
\end{enumerate} 
\end{lem}

\Proof (1) The proof follows from Lemma 1.1.

(2) Since $A(T)^p = A(p T) = A(0) = I_n$, we have $\sum_{i \geq 0} T^{ip} N_i^p = I_n$, which 
implies $N_i^p = O$ for all $i \geq 1$.

(3) Using Lemma 1.1 and the following Lemma 1.4, we have
\begin{eqnarray*}
\lefteqn{N_{i_0} N_{i_1 p} \cdots N_{i_r p^r} } \\ [2mm] 
 &=& \binom{i_0 + i_1p}{i_0} N_{i_0 + i_1p} N_{i_2 p^2} \cdots N_{i_r p^r}  =  N_{i_0 + i_1p} N_{i_2 p^2} \cdots N_{i_r p^r} \\ [2mm] 
 &=& \binom{i_0 + i_1p + i_2p^2}{i_0 + i_1p} N_{i_0 + i_1p + i_2p^2} N_{i_3 p^3} \cdots N_{i_r p^r} 
  =  N_{i_0 + i_1p + i_2p^2} N_{i_3 p^3} \cdots N_{i_r p^r}  \\ [2mm] 
 &=& \cdots \\ [2mm] 
 &=&  \binom{i_0 + i_1p + \cdots + i_{r-1} p^{r-1} + i_r p^r}{i_0 + i_1p + \cdots
 + i_{r-1} p^{r-1}} N_{i_0 + i_1p + \cdots +  i_{r-1} p^{r-1} + i_r p^r} 
 = N_{i_0 + i_1p + \cdots + i_r p^r} 
\end{eqnarray*}
for all $0 \leq i_0, i_1, \ldots, i_r \leq p - 1$. 
By Lemmas 1.1 and 1.4 again, we know that $N_{p^\ell} N_{(j - 1) p^\ell} = j N_{j p^\ell}$ for all $1 \leq j \leq p-1$ and $\ell \geq 0$, 
which implies $ N_{j p^\ell} = N_{p^\ell}^j / j!$ for all $1 \leq j \leq p-1$ and $\ell \geq 0$. 
Thus we have 
\[
N_{i_0} N_{i_1 p} \cdots N_{i_r p^r} 
 = \frac{N_1^{i_0}}{i_0!} \cdot \frac{N_p^{i_1}}{i_1!} \cdots \frac{N_{p^r}^{i_r}}{i_r!} . 
\]
\QED

The following lemma is well known and is called as Lucas' theorem (see \cite{Lucas}). 

\begin{lem}
Let $p$ be a prime number. Let $i, j$ be integers. 
Let $i = i_0 + i_1p + \cdots + i_r p^r$ and $j = j_0 + j_1 p + \cdots + j_r p^r$ be the $p$-adic expansions of $i$ and $j$, 
respectively. Then we have 
\[
 \binom{i}{j} \equiv \binom{i_0}{j_0} \binom{i_1}{j_1} \cdots \binom{i_r}{j_r} \pmod p , 
\]
where we use the convention that $\ds\binom{a}{b} = 0$ if $a < b$. 
\end{lem}

Given a polynomial matrix $A(T) \in \Mat(n, k[T])$, we can define the truncated exponential $\Exp(A(T))$ of $A(T)$ as 
\[
 \Exp(A(T)) := \sum_{i = 0}^{p-1} \frac{1}{i!} \, A(T)^i . 
\]
We say that a polynomial $f(T) \in k[T]$ is a {\it $p$-polynomial} if $f(T)$ can be written in the form 
\[
 f(T) = \sum_{i \geq 0} a_i T^{p^i}, 
\]
where $a_i \in k$ for all $i \geq 0$. 
Let $\mf{P}$ be the set of all $p$-polynomials. This $\mf{P}$ becomes a ring with multiplication  
by the composition of functions, i.e., $f(T) \circ g(T) := f(g(T))$ for all $f(T), g(T) \in \mf{P}$. 
Note that $T$ is the unity of $\mf{P}$, i.e., $1_{\mf{P}} = T$.

Clearly, for any $p$-polynomial $f(T) \in \mf{P}$ and any matrix $N \in \Mat(n, k)$ satisfying $N^p = O$, 
the polynomial matrix $\Exp(f(T) N)$ becomes an exponential matrix.

\begin{lem} 
Let $A(T) = \sum_{i \geq 0} T^i N_i\in \Mat(n, k[T])$ be an exponential matrix,  
where $N_i \in \Mat(n, k)$ for all $i \geq 0$. 
Then we have 
\[
 A(T) = \prod_{i \geq 0}  \Exp(T^{p^i} N_{p^i}) . 
\]
\end{lem}

\Proof We have 
\begin{eqnarray*}
A(T)
 &=& \sum_{i_0 = 0}^{p-1} \sum_{i_1 = 0}^{p-1} \cdots \quad T^{i_0 + i_1p + \cdots} N_{i_0 + i_1p + \cdots } 
 = \sum_{i_0 = 0}^{p-1} \sum_{i_1 = 0}^{p-1} \cdots \quad T^{i_0} T^{i_1 p} \cdots \frac{N_1^{i_0}}{i_0!} \cdot \frac{N_p^{i_1}}{i_1!} \cdots \\
 &=& \sum_{i_0 = 0}^{p-1} T^{i_0} \frac{N_1^{i_0}}{i_0!} \cdot  \sum_{i_1 = 0}^{p-1} T^{i_1 p} \frac{N_p^{i_1}}{i_1!} \cdots 
 = \Exp(TN_1) \cdot \Exp(T^p N_p) \cdots . 
\end{eqnarray*}
\QED

\subsection{Basic properties of an exponential matrix}

\subsubsection{Transposing exponential matrices}

Let $A = (a_{i, j})$ be a matrix of $\Mat(n, R)$, where $R$ is a commutative ring with unity. 
We can draw the other diagonal line $l$ of the matrix $A$. 
We denote by ${^\tau}A$ the transposed matrix of $A$ with respect to the other diagonal line $l$, i.e., 
the $(i, j)$-th entry of ${^\tau A}$ is defined by $a_{n - j + 1, n - i + 1}$ for all $1 \leq i, j \leq n$.
We say that ${^\tau A}$ is {\it the other transpose} of $A$. 

In the following lemma, we note basic properties of the other transpose of matrices.

\begin{lem}
Let $A, B \in \Mat(n, R)$ and let $c \in R$. Then the following assertions hold true:  
\begin{enumerate} 
\item[\rm (1)] ${^\tau}({^\tau}A )  = A$. 
\item[\rm (2)] ${^\tau (A + B)} = {^\tau A} + {^\tau B}$. 
\item[\rm (3)] ${^\tau ( cA )} = c \, {^\tau A}$. 
\item[\rm (4)] ${^\tau  (A B)} = {^\tau B} \, {^\tau A}$. 
\item[\rm (5)] For all $P \in GL(n, R)$, we have ${^\tau ( P^{-1} ) } = ({^\tau P})^{-1}$.  
\item[\rm (6)] For all $P \in B(n, k)$, we have $P^{-1} \in B(n, k)$. 
\end{enumerate} 
\end{lem}

For exponential matrices, we have the following basic lemma: 

\begin{lem} 
Let $A(T) \in \Mat(n, k[T])$ be an exponential matrix. Then the following assertions {\rm (1)} and {\rm (2)} hold true: 
\begin{enumerate}
\item[\rm (1)] The other transpose ${^\tau}A(T)$ of $A(T)$ is also an exponential matrix. 
\item[\rm (2)] The transpose ${^t}A(T)$ of $A(T)$ is also an exponential matrix. 
\end{enumerate} 
\end{lem}

\subsubsection{Triangulability of exponential matrices}

The following lemma implies that any exponential matrix is triangulable by a regular matrix whose entries belong to $k$. 

\begin{lem}
Let $A(T) \in \Mat(n, k[T])$ be an exponential matrix. 
Then there exists a regular matrix $P \in GL(n, k)$ such that $P^{-1} A(T) P$ is an upper triangular matrix 
of $\Mat(n, k[T])$.  
\end{lem}

\Proof We can express $A(T)$ as $A(T) = \sum_{i \geq 0} T^i N_i $, where $N_i \in \Mat(n, k)$ for all $i \geq 0$. 
There exists an integer $r \geq 2$ such that $N_j = O$ for all $j \geq r + 1$. 
Since $N_i N_j = N_j N_i$ for all $1 \leq i, j \leq r$ and $N_i^p = O$ for all $1 \leq i \leq r$,   
there exists a regular matrix $P \in GL(n, k)$ such that $P^{-1} N_i P$ is an upper triangular matrix for all $1 \leq i \leq r$. 
So, $P^{-1} A(T) P$ is an upper triangular matrix. 
\QED

Based on the above Lemma 1.8, we are interested in the entries of an upper triangular exponential matrix.

\begin{lem} 
Let $A(T) = (a_{i, j}(T)) \in \Mat(n, k[T])$ be an upper triangular exponential matrix. 
Then the following assertions {\rm (1)} and {\rm (2)} hold true: 
\begin{enumerate} 
\item[\rm (1)] $a_{i, i}(T) = 1$ for all $1 \leq i \leq n$. 
\item[\rm (2)] $a_{i, i + 1}(T)$ is a $p$-polynomial for all $1 \leq i \leq n - 1$. 
\end{enumerate} 
\end{lem}

\Proof Assertion (1) follows from $A(T)^p = I_n$. Assertion (2) follows from assertion (1). 
\QED

So, we introduce the notation $\mf{U}_n$  
which is defined as the set of all upper triangular polynomial matrices $A(T) = (a_{i, j}(T))$ of $\Mat(n, k[T])$ satisfying 
$a_{i, i}(T) = 1$ for all $1 \leq i \leq n$. Clearly, $\mf{U}_n$ is a subgroup of $GL(n, k[T])$.

Let $n$ be a positive integer. 
An ordered sequence $[d_1, \ldots, d_t]$ of positive integers $d_i$ $(1 \leq i \leq t)$ is said to be an {\it ordered partition} of $n$ 
if $[d_1, \ldots, d_t]$ satisfies $\sum_{i = 1}^t d_i = n$. 
For any $A(T) = (a_{i, j}(T)) \in \mf{U}_n$, we say that {\it $A(T)$ has an ordered partition $[d_1, d_2, \ldots, d_t]$} if 
$A(T)$ satisfies 
\begin{eqnarray*}
 \{ i \in \{ 1, \ldots, n \} \mid a_{i, i + 1}(T) = 0 \} = \{ d_1,\;  d_1 + d_2, \; \ldots, \; d_1 + \cdots + d_t  \} ,  
\end{eqnarray*}
where we let $a_{n, n + 1}(T) = 0$. 
We denote by $\mf{U}_{[d_1, \ldots ,d_t]}$ the set of all polynomial matrices of $\mf{U}_n$ having the 
partition $[d_1, \ldots, d_t]$. 
Clearly, $\mf{U}_{[d_1, \ldots ,d_t]}$ becomes a subgroup of $\mf{U}_n$. 
We have 
\[
 \mf{U}_n = \bigcup_{ \text{ $[d_1, \ldots, d_t]$ is an ordered partition of $n$}  } \mf{U}_{[d_1, \ldots, d_t]} . 
\]

We remark that even if two upper triangular exponential matrices $A(T)$ and $B(T)$ are equivalent, 
the matrices $A(T)$, $B(T)$ do not necessarily have same partitions. 
For example, letting 
\begin{eqnarray*}
P := 
\left(
\begin{array}{cccc}
 1 & 0 & 0 & 0 \\
 0 & 0 & 1 & 0 \\
 0 & 1 & 0 & 0 \\ 
 0 & 0 & 0 & 1
\end{array}
\right), 
\quad 
A(T) := 
\left(
\begin{array}{cc | cc}
 1 & T & 0 & 0 \\
 0 & 1 & 0 & 0 \\ 
\hline 
 0 & 0 & 1 & T^p\\ 
 0 & 0 & 0 & 1
\end{array}
\right), 
\quad 
B(T) := P^{-1} A P = 
\left(
\begin{array}{c | c | c |c}
 1 & 0 & T & 0 \\ 
\hline 
 0 & 1 & 0 & T^p \\ 
\hline 
 0 & 0 & 1 & 0 \\ 
\hline 
 0 & 0 & 0 & 1
\end{array}
\right) , 
\end{eqnarray*}
we know that $A(T)$ has the partition $[2, 2]$ but $B(T)$ has the partition $[1, 1, 1, 1]$.

\subsection{Basic properties of $\leadsto$ and $( \quad )^E$ }

For all subsets $S, S'$ of $\Mat(n, k[T])$, we write 
$S \leadsto S'$ 
if for any $A(T) \in S$ there exists a regular matrix $P \in GL(n, k)$ such that $P^{-1} A(T) P \in S'$. 
We have the following basic lemma:

\begin{lem} 
For all subsets $S, S', S'', T, T'$ of $\Mat(n, k[T])$, the following assertions {\rm (1), (2), (3)} hold true: 
\begin{enumerate} 
\item[\rm (1)] $S \leadsto S$. 
\item[\rm (2)] If $S \leadsto S'$ and $S' \leadsto S''$, then $S \leadsto S''$. 
\item[\rm (3)] If $S \leadsto S'$ and $T \leadsto T'$, then $S \cup T \leadsto S' \cup T'$. 
\end{enumerate} 
\end{lem}

For any subset $S$ of $\Mat(n, k[T])$, we denote by $S^E$ the subset of $S$ constituting of all exponential matrices belonging to $S$. 
Clearly, $S^E = S \cap \E(n, k[T])$. So, we have the following basic lemma:

\begin{lem}
For all subsets $S, T$ of $\Mat(n, k[T])$, the following assertions {\rm (1), (2), (3), (4)} hold true: 
\begin{enumerate}
\item[\rm (1)] If $S \subset T$, then $S^E \subset T^E$. 
\item[\rm (2)] $(S \cup T)^E = S^E \cup T^E$. 
\item[\rm (3)] $(S \cap T)^E = S^E \cap T^E$. 
\item[\rm (4)] If $S \leadsto T$, then $S^E \leadsto T^E$. 
\end{enumerate}
\end{lem}

Triangulability of exponential matrices implies the following lemma: 

\begin{lem}
We have 
\[
 \E(n, k[T]) \leadsto \mf{U}_n^E  =  \bigcup_{\text{$[ d_1, \ldots, d_t ]$ is an ordered partition of $n$}}  \mf{U}_{[d_1, \ldots, d_t]}^E. 
\]
\end{lem}

\subsection{$\mf{U}_{[n]}^E$, $\mf{U}_{[n, 1]}^E$, $\mf{U}_{[1, n]}^E$}

In this Subsection 1.5, we consider necessary and sufficient conditions for a matrix $A(T)$ of $\mf{U}_n^E$ 
to have the partition $[n]$ (see Theorem 1.13), 
and for $A(T)$ of $\mf{U}_{n + 1}^E$ to have the partitions $[n, 1]$ and $[1, n]$, respectively 
(see Theorems 1.21 and 1.25, respectively). 
And we give equivalent forms of exponential matrices belonging to 
$\mf{U}_{[n]}^E$, $\mf{U}_{[n, 1]}^E$ and $\mf{U}_{[1, n]}^E$, respectively 
(see Lemmas 1.16, 1.23 and 1.27, respectively).

\subsubsection{$[n]$}

Let $\nu_n = (\varepsilon_{i, j}) \in \Mat(n, k)$ be the nilpotent matrix such that 
\begin{eqnarray*}
\varepsilon_{i, j} = 
\left\{
\begin{array}{ll}
1 \qquad & \text{ if \; $1 \leq i \leq n - 1$ \; and \; $j = i + 1$},  \\
0 &\text{ otherwise. } 
\end{array}
\right.
\end{eqnarray*}
Clearly, $\nu_n^n = O$, and $\nu_n$ has the partition $[n]$. 
For considering $n \times n$ exponential matrices having the partition $[n]$, 
we define a (striped) subset $\mf{J}_{[n]}$ of $GL(n, k[T])$ as 
\[
\mf{J}_{[n]} 
:= 
\{ I_n + s_1 \nu_n + \cdots + s_{n - 1} \nu_n^{n - 1} \in \Mat(n, k[T]) \mid s_1, \ldots, s_{n - 1} \in k[T] \text{ and } s_1 \ne 0 \} . 
\]
Let $\nu_n^0 := I_n$. 
Clearly, any $A(T) \in \mf{J}_{[n]}$ is upper triangular and has the partition $[n]$.

In the following theorem, we give a necessary and sufficient condition for a matrix $A(T)$ 
of $\mf{U}_n^E$ to have the partition $[n]$. 

\begin{thm}
Let $A(T) \in \mf{U}_n^E$.
Then $A(T) \in \mf{U}_{[n]}^E$ if and only if 
$A(T)$ satisfies the following conditions {\rm (1)} and {\rm (2)}: 
\begin{enumerate} 
\item[\rm (1)] The size $n$ of $A(T)$ satisfies $1 \leq n \leq p$. 
\item[\rm (2)] There exists a matrix $P \in B(n, k)$ such that $P^{-1} A(T) P \in \mf{J}_{[n]}$. 
\end{enumerate} 
Furthermore, in particular when $n > p$, we have $\mf{U}_{[n]}^E = \emptyset$. 
\end{thm}

\Proof Assume that $A(T)$ satisfies the conditions (1) and (2). 
Since $P^{-1} A(T) P$ has the partition $[n]$ and $P \in B(n, k)$, we know that $A(T)$ has the partition $[n]$ (cf. Lemma 1.19). 
Conversely, assume that $A(T)$ has the partition $[n]$. 
Let $m$ be the minimum of the set of all positive integers $\ell$ such that $(A(T) - I_n)^\ell = O$. 
Clearly, $m = n$. Since $(A(T) - I_n)^p = O$, we have $1 \leq n \leq p$. 
We shall prove assertion (2). 
We proceed by induction on $n$. The proof is clear for $A(T)$ of size one-by-one. 
Let $n \geq 1$ be an integer, and hypothesize that the implication holds true for $A(T)$ of size $n \times n$.
Consider any upper triangular exponential matrix $A(T) = (a_{i, j})$ of degree $n + 1$ having the partition $[n + 1]$. 
Let $A^\flat(T) = (a_{i, j})_{1 \leq i \leq n, \; 1 \leq j \leq n} \in \Mat(n, k[T])$ be the submatrix of $A(T)$. 
Clearly, $A^\flat(T)$ is an upper triangular exponential matrix having the partition $[n]$. 
By the induction hypothesis, there exists a matrix $Q \in B(n, k)$ such that $Q^{-1} A^\flat(T) Q \in \mf{J}_{[n]}$. 
So, we can write $Q^{-1}A^\flat(T) Q$ as  
\[
 Q^{-1} A^\flat(T) Q = I_n + s_1 \nu_n + \cdots + s_{n - 1} \nu_n^{n - 1}, 
\]
where $s_i \in k[T]$ $(1 \leq i \leq n -1)$. Letting 
\[
\wt{Q} := 
\left( 
\begin{array}{c | c}
Q & \bm{0} \\ 
\hline 
\bm{0} & 1
\end{array}
\right) \in B(n + 1, k), 
\]
we have 
\begin{eqnarray*}
\wt{Q}^{-1} A(T) \wt{Q}
= 
\left(
\begin{array}{c | c}
Q^{-1} A^\flat(T) Q & 
Q^{-1} 
\left(
\begin{array}{c}
 a_{1, n + 1} \\ 
 \vdots  \\
 a_{n, n +1}
\end{array}
\right)
 \\
\hline 
 \bm{0} & 1
\end{array}
\right) . 
\end{eqnarray*}
Since $Q^{-1} \in B(n, k)$, we know $\wt{Q}^{-1} A(T) \wt{Q}$ also has the partition $[n + 1]$. 
The polynomial matrix $\wt{Q}^{-1} A(T) \wt{Q}$ satisfies the hypotheses in the following Lemma 1.14. 
Thus, we have the desired matrix $P$. 
\QED

\begin{lem}
Let $A(T) = (a_{i, j}) \in \Mat( n + 1, k[T] )$ be an exponential matrix with the form 
\begin{eqnarray*}
A(T) = 
\left( 
\begin{array}{c|c}
 A^\flat(T) & \bm{b} \\
\hline 
 \bm{0} & 1 
\end{array}
\right) , 
\end{eqnarray*}
where $A^\flat(T) = I_n + s_1 \nu_1 + \cdots + s_{n - 1} \nu_n^{n-1}  \in \mf{J}_{[n]}$ $(s_1, \ldots, s_{n - 1} \in k[T])$ and 
$\bm{b} \in k[T]^n$. Assume that $A(T)$ has the partition $[n + 1]$. 
Then there exists a matrix $P \in B(n + 1, k)$ such that $P^{-1} A(T) P \in \mf{J}_{[n + 1]}$. 
\end{lem}

\Proof The proof is clear for $n = 1$. So, let $n \geq 2$. 
Write 
\[
\bm{b} = 
\left( 
\begin{array}{c}
 b_n \\
 \vdots \\ 
 b_2 \\
 b_1
\end{array} 
\right) \qquad ( b_1, b_2, \ldots, b_n \in k[T], \; b_1 \ne 0) . 
\]
Since $A(T)$ is exponential, we have $b_2(T') + s_1(T) b_1(T') + b_2(T) = b_2(T + T')$, which implies $b_1(T) = \lambda s_1(T)$ 
for some $\lambda \in k \backslash\{ 0 \}$. 
Let 
\[
 Q := 
\left(
\begin{array}{c | c}
 I_n & \bm{0} \\
\hline 
 \bm{0} & 1/\lambda
\end{array}
\right). 
\]
Thus, by regarding $Q^{-1} A(T) Q$ as $A(T)$, we may assume that $b_1 = s_1$ in the matrix $A(T)$. 
If $n = 2$, we complete the proof. 
So, let $n \geq 3$. 
Suppose that $b_i = s_i$ in $A(T)$ for all $1 \leq i \leq r $ where $r \leq n - 2$. Since 
\[
 b_{r + 2}(T') + s_1(T) b_{r + 1}(T')  + \cdots + s_r(T) b_2(T') + s_{r + 1}(T) b_1(T') + b_{r + 2}(T) = b_{r + 2}(T + T'), 
\]
we know that $s_1(T) b_{r + 1}(T') + s_{r + 1}(T) s_1(T')$ is symmetric, which implies that 
\[
b_{r + 1} = s_{r + 1} + \mu s_1
\]
for some $\mu \in k$. Let 
\[
R := 
\left(
\begin{array}{c | c}
 I_n & - \mu \bm{e}_{n - r + 1} \\
\hline 
 \bm{0} & 1 \\
\end{array}
\right)
\in B(n + 1, k), 
\]
where for $1 \leq j \leq n$, $\bm{e}_j \in k[T]^n$ is defined, as follows:  
The $j$-th entry of $\bm{e}_j \in k[T]^n$ is $1$ and the other entries are zeroes. 
Regarding $R^{-1}A(T)R$ as $A(T)$, we may assume that $b_{r + 1} = s_{r + 1}$ in the matrix $A(T)$. 
By repeating the above arguments in finitely many steps, we have the desired matrix $P \in B(n + 1, k)$. 
\QED

The following is a corollary of Theorem 1.13: 

\begin{cor}
We have 
\begin{eqnarray*}
\left\{ 
\begin{array}{ll}
 \mf{U}_{[n]}^E \leadsto \mf{J}_{[n]}^E \qquad & \text{ if \quad $1 \leq n \leq p$}, \\
 \mf{J}_{[n]}^E = \emptyset  & \text{ if \quad $n > p$} . 
\end{array}
\right. 
\end{eqnarray*} 
\end{cor}

Based on the above corollary, we are interested in expressing any matrix belonging to $\mf{J}_{[n]}^E$. 
So, for any $A(T) \in \mf{J}_{[n]}$ with the form 
$A(T) = I_n + f_1 \nu_n + \cdots + f_{n - 1} \nu_n^{n - 1} \in \mf{J}_{[n]}$ $(f_1, \ldots, f_{n - 1} \in k[T], \; f_1 \ne 0)$, 
we define a non-negative integer $d(A(T))$ as 
\[
 d(A(T)) := n - \min \{ i \in \{1, 2, \ldots, n - 1\} \mid f_i \ne 0 \} . 
\]
Clearly, $1 \leq d(A(T)) \leq n - 1$.

\begin{lem}
Let $n$ be an integer satisfying $1 \leq n \leq p$. 
For any $A(T) \in \Mat(n, k[T])$, the following conditions {\rm (1)} and {\rm (2)} are equivalent: 
\begin{enumerate} 
\item[\rm (1)] $A(T) \in \mf{J}_{[n]}^E$. 
\item[\rm (2)] $A(T)$ has the form 
\[
A(T) = \prod_{i = 1}^{n - 1} \Exp(f_i \nu_n^i) = \Exp\left( \sum_{i = 1}^{n - 1} f_i \nu_n^i \right) 
\qquad 
\left(
\begin{tabular}{l}
\text{$f_i$ $(1 \leq i \leq n - 1)$ are $p$-polynomials, } \\
\text{and $f_1 \ne 0$}. 
\end{tabular}
\right) .  
\]
\end{enumerate} 
\end{lem}

\Proof We have only to prove (1) $\Longrightarrow$ (2) (since the implication (2) $\Longrightarrow$ (1) is clear). 
We proceed by induction on the value $d := d(A(T))$ of $A(T)$ belonging to $\mf{J}_{[n]}^E$. 
If $d = 1$, then $A(T)$ clearly has the desired form.  
If $d \geq 2$, then the coefficient polynomial $f_d$ of the stripe $f_d \nu_n^d$ of $A(T)$ is a non-zero $p$-polynomial. 
Let $B(T) := A(T) \cdot \Exp( - f_d \nu_n^d )$. 
We know that $B(T) \in \mf{J}_{[n]}^E$ and $d(B(T)) < d$. 
By the induction hypothesis, we can write $B(T)$ as 
$B(T) =  \prod_{i = 1}^{n - 1} \Exp(g_i \nu_n^i)$ for some $p$-polynomials $g_i$ $(1 \leq i \leq n - 1)$. 
Thus $A(T)$ has the desired form. 
\QED

The following lemma gives another expression of the form of $A(T)$ in (2) of Lemma 1.16: 

\begin{lem}
Let $n$ be an integer satisfying $2 \leq n \leq p$ and let $f_1, \ldots, f_{n - 1} \in k[T]$. 
Then we have 
\[
 \Exp\left( \sum_{i = 1}^r f_i \nu_n^i \right)
 = \sum_{\ell = 0}^{n - 1} \left( \sum_{\substack{ i_1 + 2 i_2 + \cdots + r i_r = \ell \\ i_1, i_2, \ldots, i_r \geq 0}}
   \frac{f_1^{i_1} f_2^{i_2} \cdots f_r^{i_r}}{i_1! i_2!  \cdots  i_r!}  \right) \nu_n^\ell 
\qquad (1 \leq r \leq n - 1) .  
\]
\end{lem}

\Proof The proof is straightforward (we can prove the equality by induction on $r$ ($1 \leq r \leq n - 1$)). 
\QED

For any $0 \leq \ell \leq n$, let ${\rm UT}_{\geq \ell}$ be the set of all upper triangular polynomial matrices 
$A = (a_{i, j})$ of $\Mat(n, k[T])$ satisfying $a_{i, j} = 0$ for all $1 \leq j \leq i + \ell  - 1 \leq n$. 

\begin{lem}
Let $A(T), B(T) \in \mf{U}_{[n]}$. 
Let $P(T)$ be a polynomial matrix of $\Mat(n, k[T])$ such that 
\[
 A(T) P(T) -  P(T) B(T) \in {\rm UT}_{\geq 1}. 
\] 
Then $P(T)$ is upper triangular. 
\end{lem}

\Proof Consider the first column of $A(T) P(T) - P(T) B(T) = C(T)$, where $C(T) \in {\rm UT}_{\geq 1}$. 
The first column of $P(T) = (p_{i, j})$ satisfies $p_{i, 1} = 0$ for all $2 \leq i \leq n$. 
Using the induction on $n$, we know that $P(T)$ is upper triangular.  
\QED

Given a matrix $A(T) \in \mf{U}_{[n]}$, we shall consider a regular matrix $P$ of $GL(n, k)$ such that 
$P^{-1} A(T) P \in \mf{U}_{[n]}$. 

\begin{lem}
Let $A(T) \in \mf{U}_{[n]}$ and let $P \in GL(n, k)$. 
Then the following conditions {\rm (1)} and {\rm (2)} are equivalent: 
\begin{enumerate}
\item[\rm (1)] $P^{-1} A(T) P \in \mf{U}_{[n]}$. 
\item[\rm (2)] $P \in B(n, k)$. 
\end{enumerate} 
\end{lem}

\Proof The emplication  (1) $\Longrightarrow$ (2) is clear from Lemma 1.18. So we shall prove (2) $\Longrightarrow$ (1). 
Let $A(T) = (a_{i, j}(T))$, $P = (p_{i, j})$, $P^{-1} A P = (b_{i, j} )$ and $P^{-1} = (p_{i, j}' )$.
We know that $P^{-1} \in B(n, k)$ and $p_{i, i}' = 1/p_{i, i}$ for all $1 \leq i \leq n$. 
Since the $(i, i + 1)$-entry of $P^{-1} \cdot P$ is zero, we have 
\[
 p_{i, i +1}' \cdot {p_{i, i}} + p_{i, i +1}' \cdot p_{i + 1, i +1} = 0 . 
\]
So, the $(i, i + 1)$-entry $b_{i, i + 1}$ of $P^{-1} A P $ can be calculated as 
\[
 b_{i, i +1} = p_{i, i}' \cdot ( p_{i, i +1} + a_{i, i + 1} p_{i + 1, i + 1} ) +p_{i, i  + 1}' p_{i + 1, i + 1} =  a_{i, i + 1}  \cdot p_{i, i}' \cdot p_{i + 1, i + 1} .  
\]
\QED

\subsubsection{$[n, 1]$}

We denote by $\mf{J}_{[n, 1]}$ the set of all matrices $A(T)$ of $GL(n + 1, k[T])$ with the form 
\[
 A(T) = 
\left(
\begin{array}{c | c}
 U(T) & u_n \bm{e}_1\\
\hline 
 \bm{0} & 1 
\end{array} 
\right)
\qquad \quad 
\left( U(T) = I_n + \sum_{i = 1}^{n - 1} u_i \nu_n^i \in \mf{J}_{[n]}, \quad u_i \in k[T] \quad (1 \leq i \leq n) \right) .  
\]
Clearly, we have $\mf{J}_{[n, 1]} \subset \mf{U}_{[n, 1]}$.
We define a subset $\mf{J}_{[n, 1]}^0$ of $\mf{J}_{[n, 1]}$ as the set of all matrices $A(T)$ of $\mf{J}_{[n ,1]}$ satisfying $u_n = 0$, 
and a subset $\mf{J}_{[n, 1]}^1$ as the set of all matrices $A(T)$ of $\mf{J}_{[n, 1]}$ satisfying 
$u_1$ and $u_n$ are linearly independent over $k$. 

For any polynomial matrix $A(T)$ of $\mf{J}_{[n, 1]}$, if $u_1$ and $u_n$ are linearly dependent over $k$, 
we have $u_n = \lambda u_1$ for some $\lambda \in k$. Letting 
\[
P := 
\left(
\begin{array}{c | c}
 I_n & - \lambda \bm{e_2} \\
\hline 
 \bm{0} & 1 
\end{array} 
\right) , 
\]
we have $P^{-1} A(T) P \in \mf{J}_{[n, 1]}^0$. Now, the following lemma is obtained: 

\begin{lem} 
We have 
\[
\mf{J}_{[n, 1]} \leadsto \mf{J}_{[n ,1]}^0 \cup \mf{J}_{[n, 1]}^1. 
\] 
\end{lem}

In the following theorem, we give a necessary and sufficient condition  for a matrix $A(T)$ 
of $\mf{U}_{n + 1}^E$ to have the partition $[n, 1]$. 

\begin{thm}
Let $A(T) \in \mf{U}_{n + 1}^E$. 
Then $A(T) \in \mf{U}_{[n, 1]}^E$ if and only if $A(T)$ satisfies the following conditions {\rm  (1)} and {\rm (2)}: 
\begin{enumerate} 
\item[\rm (1)] $1 \leq n \leq p$. 
\item[\rm (2)] There exists a matrix $P \in B(n + 1, k)$ such that $P^{-1} A(T) P \in \mf{J}_{[n, 1]}$. 
\end{enumerate} 
Furthermore, in particular when $n > p$, we have $\mf{U}_{[n, 1]}^E = \emptyset$. 
\end{thm}

\Proof Assume $A(T) \in \mf{U}_{[n, 1]}^E$. 
We can write $A(T)$ as 
\[
A(T) = 
\left(
\begin{array}{c | c}
 A^\flat(T) & \bm{b} \\
 \hline 
\bm{0} & 1 
\end{array}
\right) \qquad (\; A^\flat(T) \in \mf{U}_n, \quad \bm{b} = \transpose{(b_{n - 1}, \ldots, b_2, b_1, 0)} \in k[T]^n\;) . 
\]
Note that $A^\flat(T) \in \mf{U}_{[n]}^E$.  
We know form Theorem 1.13 that $1 \leq n \leq p$ and there exists a matrix $P \in B(n, k)$ such that 
$P^{-1} A^\flat(T) P \in \mf{J}_{[n]}$. 
So, let 
\[
 P' := 
\left(
\begin{array}{c | c}
 P & \bm{0} \\
\hline 
 \bm{0} & 1
\end{array} 
\right)
\qquad \text{ and } \qquad A'(T) := P'^{-1} A(T) P' . 
\]
Clearly, we have
\begin{eqnarray*}
P' \in B(n + 1, k) \qquad \text{ and } \qquad 
A'(T) = 
\left(
\begin{array}{c | c}
P^{-1} A^\flat(T) P & P^{-1} \bm{b} \\
\hline 
 \bm{0} & 1  
\end{array} 
\right) . 
\end{eqnarray*}
So, let $\bm{\beta} := P^{-1} \bm{b}$. Write $\bm{\beta} = \transpose (\beta_{n - 1}, \ldots, \beta_2, \beta_1, 0)$, where 
$\beta_i \in k$ for all $1 \leq i \leq n - 1$. 
In the case where $\beta_i = 0$ for all $1 \leq i \leq n - 2$, then $P'^{-1} A(T) P' \in \mf{J}_{[n, 1]}$ is clear. 
So, we consider the case where there exists an integer $1 \leq r \leq n -2$ such that $b_r \ne 0$ but $b_1 = \cdots = b_{r - 1} = 0$. 
Write $P'^{-1} A(T) P'  = (a_{i, j}(T))_{1 \leq i, j \leq n + 1}$. 
Comparing the $(n - r - 2, n + 1)$-th entries of the both sides of the equality $A'(T) A'(T') = A'(T + T')$, 
where $T, T'$ are indeterminates over $k$, we have 
\[
 b_{r + 1}(T') + b_r(T') a_{r + 1, r + 2}(T) + b_{r + 1}(T) = b_{r + 1}(T + T') . 
\]
So, $b_r(T) = \lambda a_{r + 1, r  + 2}(T)$ for some $\lambda \in k \backslash \{ 0 \}$. 
Let 
\[
 Q := 
\left(
\begin{array}{c | c}
 I_n & \bm{0} \\
 \hline 
\bm{0} & 1/\lambda 
\end{array}
\right), \;
\quad 
R := 
\left( 
\begin{array}{c | c }
 I_n & -\bm{e}_{n - r + 1} \\
\hline 
 \bm{0} & 1
\end{array} 
\right),  
\quad A''(T) := R^{-1} Q^{-1} A'(T) Q R . 
\]
Clearly, $QR \in B(n + 1, k)$. We can write $A''(T)$ as 
\[
A''(T) = 
\left( 
\begin{array}{c | c}
 P^{-1} A^\flat(T) P & \bm{b}' \\ 
\hline 
 \bm{0} & 1
\end{array}
\right) 
\qquad 
(\; \bm{b}' = \transpose(b_{n - 1}', \ldots, b_{r+ 1}', 0, \ldots, 0) \in k[T]^n \; ) . 
\]
We can repeat the above arguments in finitely many steps until we have the desired form. 
Conversely, assume that $A(T)$ satisfies the conditions (1) and (2). 
By Lemma 1.19, we have $A(T) \in \mf{U}_{[n, 1]}^E$. 
\QED

The following is a corollary of the above Theorem 1.21 and Lemma 1.20: 

\begin{cor} 
We have 
\begin{eqnarray*}
\left\{ 
\begin{array}{ll}
 \mf{U}_{[n, 1]}^E \leadsto \mf{J}_{[n, 1]}^E \leadsto ( \mf{J}_{[n, 1]}^0 )^E \cup ( \mf{J}_{[n, 1]}^1 )^E 
 \qquad  & \text{ if \quad $1 \leq n \leq p$}, \\ 
 \mf{J}_{[n, 1]}^E = \emptyset  & \text{ if \quad $n > p$} . 
\end{array}
\right. 
\end{eqnarray*}
\end{cor}

The following lemma gives an expression of $\mf{J}_{[n, 1]}^E$, and the expression follows from Lemma 1.16. 

\begin{lem}
Let $n$ be an integer satisfying $1 \leq n \leq p$. 
For any $A(T) \in \Mat(n + 1, k[T])$, the following conditions {\rm (1)} and {\rm (2)} are equivalent: 
\begin{enumerate} 
\item[\rm (1)] $A(T) \in \mf{J}_{[n, 1]}^E$. 
\item[\rm (2)] $A(T)$ has the form  
\begin{eqnarray*}
A(T) = 
\left(
\begin{array}{c | c}
 \ds\Exp \left( \sum_{i = 1}^{n - 1} f_i \nu_n^i \right) & f_n \bm{e}_1 \\  
\hline 
 \bm{0} & 1
\end{array}
\right)
\qquad 
(\text{ $f_i$ $(1 \leq i \leq n)$ are $p$-polynomials, and $f_1 \ne 0$ }).
\end{eqnarray*}
\end{enumerate} 
\end{lem}

\subsubsection{$[1, n]$} 

We denote by $\mf{J}_{[1, n]}$ the set of all matrices $A(T)$ of $GL(n + 1, k[T])$ with the form 
\[
 A(T) = 
\left(
\begin{array}{c | c}
 1 & u_n \transpose\bm{e}_n \\
\hline 
\bm{0} & U(T)  
\end{array} 
\right)
\qquad \quad 
\left( U(T) = I_n + \sum_{i = 1}^{n - 1} u_i \nu_n^i \in \mf{J}_{[n]}, \quad u_i \in k[T] \quad (1 \leq i \leq n) \right) .  
\]
Clearly, we have $\mf{J}_{[1, n]} \subset \mf{U}_{[1, n]}$.
We define a subset $\mf{J}_{[1, n]}^0$ of $\mf{J}_{[1, n]}$ as the set of all matrices $A(T)$ of $\mf{J}_{[1, n]}$ satisfying $u_n = 0$, 
and a subset $\mf{J}_{[1, n]}^1$ as the set of all matrices $A(T)$ of $\mf{J}_{[1, n]}$ satisfying 
$u_1$ and $u_n$ are linearly independent over $k$.

\begin{lem} 
We have 
\[
\mf{J}_{[1, n]} \leadsto \mf{J}_{[1, n]}^0 \cup \mf{J}_{[1, n]}^1 \qquad \text{ and } \qquad \mf{J}_{[1, n]}^0 \leadsto \mf{J}_{[n, 1]}^0 . 
\] 
\end{lem}

\Proof We have only to prove the latter. 
Take any $A(T) \in \mf{J}_{[1, n]}^0$. 
We define a matrix $P = (p_{i, j} ) \in \Mat(n, k)$, as follows: 
\[
p_{i, j} 
:= 
\left\{
\begin{array}{ll}
 1  \qquad & \text{ if \quad $i  - j \equiv 1 \pmod n$}, \\ 
 0  & \text{ otherwise}. 
\end{array}
\right.
\]
Clearly, $P \in GL(n, k)$. We can calculate $P^{-1} A(T) P$ and know that $P^{-1} A(T) P \in \mf{J}_{[n, 1]}^0$. 
\QED

In the following theorem, we give a necessary and sufficient condition  
for a matrix $A(T)$ of $\mf{U}_{n + 1}^E$ to have the partition $[1, n]$.

\begin{thm}
Let $A(T) \in \mf{U}_{n + 1}^E$. 
Then $A(T) \in \mf{U}_{[1, n]}^E$ if and only if $A(T)$ satisfies the following conditions {\rm (1)} and {\rm (2)}: 
\begin{enumerate} 
\item[\rm (1)] $1 \leq n \leq p$. 
\item[\rm (2)] There exists a matrix $P \in B(n + 1, k)$ such that $P^{-1} A(T) P \in \mf{J}_{[1, n]}$. 
\end{enumerate} 
Furthermore, in particular when $n > p$, we have $\mf{U}_{[1, n]}^E = \emptyset$. 
\end{thm}

\Proof 
Note that $A(T) \in \mf{U}_{[1, n]}^E$ if and only if $^\tau A(T) \in \mf{U}_{[n, 1]}^E$. 
We know from Theorem 1.21 that the latter condition is equivalent to the following conditions (i) and (ii): 
\begin{enumerate} 
\item[\rm (i)] $1 \leq n \leq p$. 
\item[\rm (ii)] There exists a regular matrix $P \in B(n + 1, k)$ such that $P^{-1} \cdot {^\tau A}(T) \cdot P \in \mf{J}_{[n, 1]}$. 
\end{enumerate} 
By Lemmas 1.6 and 1.7, the above condition (ii) is equivalent to the condition that 
there exists a regular matrix $P' \in B(n + 1, k)$ such that $P'^{-1} \cdot A(T) \cdot P' \in \mf{J}_{[1, n]}$. 
\QED

The following is a corollary of the above Theorem 1.25 and Lemma 1.24:

\begin{cor} 
We have 
\begin{eqnarray*}
\left\{ 
\begin{array}{ll}
 \mf{U}_{[1, n]}^E
 \leadsto \mf{J}_{[1, n]}^E
 \leadsto (\mf{J}_{[1, n]}^0)^E \cup (\mf{J}_{[1, n]}^1)^E 
 \leadsto (\mf{J}_{[n, 1]}^0)^E \cup (\mf{J}_{[1, n]}^1)^E  
 \qquad & \text{ if \quad $1 \leq n \leq p$}, \\
 \mf{J}_{[1, n]}^E = \emptyset  & \text{ if \quad $p > n$}.  
\end{array}
\right. 
\end{eqnarray*}
\end{cor}

The following lemma gives an expression of $\mf{J}_{[1, n]}^E$, and the expression follows from Lemma 1.16. 

\begin{lem}
Let $n$ be an integer satisfying $1 \leq n \leq p$. 
For any $A(T) \in \Mat(n + 1, k[T])$, the following conditions {\rm (1)} and {\rm (2)} are equivalent: 
\begin{enumerate} 
\item[\rm (1)] $A(T) \in \mf{J}_{[1, n]}^E$. 
\item[\rm (2)] $A(T)$ has the form  
\begin{eqnarray*}
A(T) = 
\left(
\begin{array}{c | c}
 1 & f_n \transpose{\bm{e}}_n \\ 
\hline 
 \bm{0} & \ds\Exp \left( \sum_{i = 1}^{n - 1} f_i \nu_n^i \right)  
\end{array}
\right)  \qquad  ( \text{ $f_i$ $(1 \leq i \leq n)$ are $p$-polynomials, and $f_1 \ne 0$ }). 
\end{eqnarray*}
\end{enumerate} 
\end{lem}

\subsection{$\mf{A}(i_1, i_2, i_3)$, $\mf{A}_{i_1, i_3}$}

Let $i_1 \geq 1$, $i_2 \geq 0$, $i_3 \geq 1$ be integers.  
For any matrix $\alpha(T) \in \Mat_{i_1, i_3}(k[T])$, we define a matrix $\Lambda(i_1, i_2, i_3; \alpha(T))$ 
of $\Mat(i_1 + i_2 + i_3, k[T])$ as 
\[
 \Lambda(i_1, i_2, i_3; \alpha(T)) = 
\left(
\begin{array}{c | c | c}
 I_{i_1} & O & \alpha(T) \\
\hline 
 O &  I_{i_2} & O \\
\hline 
 O & O & I_{i_3}
\end{array}
\right). 
\]

We denote by $\mf{A}(i_1, i_2, i_3)$ the set of all polynomial matrices $A(T)$ of 
$\Mat(i_1 + i_2 + i_3, k[T])$ with the form $A(T) = \Lambda(i_1, i_2, i_3; \alpha(T)$), where $\alpha(T) \in \Mat_{i_1, i_3}(k[T])$. 
Clearly, $\mf{A}(i_1, i_2, i_3)$ becomes a commutative subgroup of $GL(i_1 + i_2 + i_3, k[T])$. 
We frequently use the notation $\Lambda(i_1, i_3; \alpha(T))$ in place of $\Lambda(i_1, i_2, i_3; \alpha(T)$ 
if we can understand the value of $i_1 + i_2 + i_3$ from the context, 
and also use the notation $\mf{A}(i_1, i_3)$ in place of $\mf{A}(i_1, i_2, i_3)$.

In the following lemma, we write equivalent conditions for a polynomial matrix $A(T)$ of $\mf{A}(i_1, i_2, i_3)$ 
to be exponential. Its proof is clear. 

\begin{lem}
Let $A(T) = \Lambda(i_1, i_2, i_3; \alpha(T)) \in \mf{A}(i_1, i_2, i_3)$. 
Then the following conditions {\rm (1)}, {\rm (2)}, {\rm (3)} are equivalent: 
\begin{enumerate}
\item[\rm (1)] $A(T) \in \mf{A}(i_1, i_2, i_3)^E$. 
\item[\rm (2)] $\alpha(T) + \alpha(T') = \alpha(T + T')$, where $T, T'$ are indeterminates over $k$. 
\item[\rm (3)] All entries of $\alpha(T)$ are $p$-polynomials. 
\end{enumerate}  
\end{lem}

For any matrix $M(T) \in \Mat_{m, n}(k[T])$, 
we denote by $\Rank M(T)$ the pair $(i, j)$ of 
the maximum number $i$ of $k$-linearly independent rows of $M(T)$ 
and the maximum number $j$ of $k$-linearly independent columns of $M(T)$. 

We denote by $\mf{A}(i_1, i_2, i_3)^\circ$ the set of all matrices $\Lambda(i_1, i_2, i_3; \alpha(T))$ of $\mf{A}(i_1, i_2, i_3)$ 
satisfying $\Rank \alpha(T) = (i_1, i_3)$.

For any $A(T) \in \mf{A}(i_1, i_2, i_3)$, where $A(T) \ne I_{i_1 + i_2 + i_3}$, 
there exists a $P \in GL(i_1 + i_2 + i_3, k[T])$ such that $P^{-1} A(T) P = \Lambda(i_1', i_2', i_3'; \beta(T))$, 
$\Rank \beta(T) = (i_1', i_3')$, $i_1' \leq i_1$, $i_2' \geq i_2$, $i_3' \leq i_3$. 
It follows that the following lemma is obtained: 

\begin{lem}
The following assertions {\rm (1)} and {\rm (2)} hold true: 
\begin{enumerate} 
\item[\rm (1)] $\ds \mf{A}(i_1, i_2, i_3) \leadsto 
 \{ I_{i_1 + i_2 + i_3} \}
 \cup 
 \bigcup_{\substack{1 \leq i_1' \leq i_1 \\  i_2 \leq i_2' \\ 1 \leq i_3' \leq i_3 }} \mf{A}(i_1', i_2', i_3')^\circ$. 
\item[\rm (2)] 
$\ds
 \mf{A}(i_1, i_2, i_3)^E \leadsto 
 \{ I_{i_1 + i_2 + i_3} \}
 \cup 
 \bigcup_{\substack{1 \leq i_1' \leq i_1 \\  i_2 \leq i_2' \\ 1 \leq i_3' \leq i_3 }} ( \mf{A}(i_1', i_2', i_3')^\circ )^E  
$. 
\end{enumerate} 
\end{lem}

\subsection{${V_n}^{A(T)}$ and $(V_n^*)^{A(T)}$, mutually $GL(n, k)$-disjoint sets}

\subsubsection{${V_n}^{A(T)}$ and $(V_n^*)^{A(T)}$}

Let $A(T)$ be a polynomial matrix of $\Mat(n, k[T])$. 
We denote by $V_n$ the column space of dimension $n$ over $k$, and by $V_n^*$ the row space of dimension $n$ over $k$. 
We can regard the polynomial matrix $A(T)$ as a $k$-linear map from $V$ to $k[T] \otimes_k V$ and also a $k$-linear map 
from $V^*$ to $k[T] \otimes_k V^*$, and naturally regard $V$ as a subspace of $k[T] \otimes_k V$, and 
$V^*$ as a subspace of $k[T] \otimes_k V^*$.  
We define the subspaces ${V_n}^{A(T)}$ and $(V_n^*)^{A(T)}$ as 
\[
 {V_n}^{A(T)} := \{ v \in V_n \mid  A(T) v = v \} \qquad \text{ and } \qquad (V_n^*)^{A(T)} := \{ v^* \in V_n^* \mid v^* A(T) = v^* \} .   
\]
If polynomial matrices $A(T)$ and $B(T)$ are equivalent, 
we have $\dim_k {V_n}^{A(T)} = \dim_k {V_n}^{B(T)}$ and $\dim_k (V_n^*)^{A(T)} = \dim_k (V_n^*)^{B(T)}$. 

For any polynomial matrix $A(T)$ belonging to the subsets $\mf{J}_{[n]}$, $\mf{J}_{[n, 1]}^0$, $\mf{J}_{[n, 1]}^1$, 
$\mf{J}_{[1, n]}^1$ and $\mf{A}(i_1, i_2, i_3)^\circ$, 
we can describe both ${V_n}^A(T)$ and $(V_n^*)^{A(T)}$, as follows:

\begin{lem} 
Let $n \geq 2$ and let $\lambda, \mu \geq 1$. 
Then the following assertions {\rm (1)}, {\rm (2)}, {\rm (3)}, {\rm (4)}, {\rm (5)} hold true: 
\begin{enumerate}
\item[\rm (1)] For any $A(T) \in \mf{J}_{[n]}$, we have 

\qquad\qquad
$
\left\{
\begin{array}{rcl}
 {V_n}^{A(T)} &=& \{ a_1 \, \bm{e}_1 \in V_n \mid a_1 \in k[T] \}, \\
 (V_n^*)^{A(T)} &=& \{ \alpha_n \, \transpose \bm{e}_n \in V_n^*  \mid \alpha_n \in k[T] \} .
\end{array}
\right. 
$

In particular, we have $\dim_k {V_n}^{A(T)} = 1$ and $\dim_k (V_n^*)^{A(T)} = 1$.

\item[\rm (2)] For any $A(T) \in \mf{J}_{[n, 1]}^0$, we have 

\qquad\qquad
$
\left\{
\begin{array}{rcl}
 {V_{n + 1}}^{A(T)} &=& \{ a_1 \, \bm{e}_1 + a_{n + 1}\, \bm{e}_{n + 1} \in V_{n + 1} \mid a_1, a_{n + 1} \in k[T] \}, \\
 (V_{n + 1}^*)^{A(T)} & = &
 \{ \alpha_n \, \transpose \bm{e}_n + \alpha_{n + 1} \, \transpose\bm{e}_{n + 1} \in V_{n + 1}^*  \mid \alpha_n, \alpha_{n + 1} \in k[T] \} . 
\end{array}
\right. 
$

In particular, we have $\dim_k {V_{n + 1}}^{A(T)} = 2$ and $\dim_k (V_{n + 1}^*)^{A(T)} = 2$.

\item[\rm (3)] For any $A(T) \in \mf{J}_{[n, 1]}^1$, 

\qquad\qquad
$
\left\{ 
\begin{array}{rcl}
 {V_{n + 1}}^{A(T)} &=& \{ a_1 \, \bm{e}_1 \in V_{n + 1} \mid a_1 \in k[T] \}, \\
 (V_{n + 1}^*)^{A(T)} &=&
 \{ \alpha_n \, \transpose \bm{e}_n + \alpha_{n + 1} \, \transpose\bm{e}_{n + 1} \in V_{n + 1}^*  \mid \alpha_n, \alpha_{n + 1} \in k[T] \} . 
\end{array}
\right. 
$

In particular, we have $\dim_k {V_{n + 1}}^{A(T)} = 1$ and $\dim_k (V_{n + 1}^*)^{A(T)} = 2$.

\item[\rm (4)] For any $A(T) \in \mf{J}_{[1, n]}^1$, we have 

\qquad\qquad
$
\left\{  
\begin{array}{rcl}
 {V_{n + 1}}^{A(T)} &=& \{ a_1 \, \bm{e}_1 + a_2 \, \bm{e}_2  \in V_{n + 1} \mid a_1, a_2 \in k[T] \}, \\
 (V_{n + 1}^*)^{A(T)} &=& 
 \{  \alpha_{n + 1} \, \transpose\bm{e}_{n + 1} \in V_{n + 1}^*  \mid  \alpha_{n + 1} \in k[T] \} . 
\end{array}
\right. 
$

In particular, we have $\dim_k {V_{n + 1}}^{A(T)} = 2$ and $\dim_k (V_{n + 1}^*)^{A(T)} = 1$.

\item[\rm (5)] For any $A(T)  \in \mf{A}(i_1, i_2, i_3)^\circ$, letting $n := i_1 + i_2 + i_3$, we have 

\qquad\qquad 
$
\left\{
\begin{array}{rcl}

 {V_n}^{A(T)} &=& \left\{ \ds \left. \sum_{\ell = 1}^{i_1 + i_2} a_\ell \, \bm{e}_\ell  \in V_n \; \right| \;  a_\ell \in k[T] \; \, (1 \leq \ell \leq i_1 + i_2 ) \right\}, \\ [5mm] 
 {V_n^*}^{A(T)} &=& \left\{ \ds \left. \sum_{\ell =  i_1 + 1}^{n} \alpha_\ell \, \transpose\bm{e}_\ell  \in V_n^* \; \right| \;  \alpha_\ell \in k[T] \; \, (i_3 + 1 \leq \ell \leq n) \right\} . 
\end{array}
\right.
$

In particular, we have 
$\dim_k {V_n}^{A(T)} = i_1 + i_2  = n - i_3$ and $\dim_k (V_{n}^*)^{A(T)} = i_2 + i_3 = n - i_1$. 
\end{enumerate}
\end{lem}

\subsubsection{Mutually $GL(n, k)$-disjoint sets}

Let $S_1, \ldots, S_N$ be subsets of $GL(n, k)$. 
We say that $S_1, \ldots, S_N$ are {\it mutually $GL(n, k)$-disjoint} 
if $(PS_i P^{-1}) \cap (P' S_j P'^{-1}) = \emptyset$ for all $1 \leq i, j \leq j$ with $i \ne j$ and for all $P, P' \in GL(n, k)$.

\begin{lem}
Let $n \geq 3$, and let $i_1, i_2, i_3$ be integers satisfying $i_1 \geq 1$, $i_3 \geq 1$ and $i_1 + i_2 + i_3 = n + 1$. 
Then the five subsets $\mf{J}_{[n + 1]}$, $\mf{J}_{[n ,1]}^0$, $\mf{J}_{[n, 1]}^1$, $\mf{J}_{[1, n]}^1$, and $\mf{A}(i_1, i_2, i_3)$ 
are mutually $GL(n + 1, k)$-disjoint. 
\end{lem}

\Proof 
Let $\mf{J} := \mf{J}_{[n + 1]} \cup \mf{J}_{[n ,1]}^0 \cup \mf{J}_{[n, 1]}^1 \cup \mf{J}_{[1, n]}^1$. 
By Lemma 1.30, it suffices to show that $\mf{J}$ and $\mf{A}(i_1, i_2, i_3)$ are mutually $GL(n + 1, k)$-disjoint. 
Choose any $A(T) \in  \mf{J}$ and choose any $B(T) \in \mf{A}(i_1, i_2, i_3)$. 
Let $P \in GL(n, k)$ satisfy $P^{-1} A(T) P = B(T)$. Squaring the both sides of the equality $P^{-1} ( A(T) - I_{n + 1} ) P = B(T) - I_{n + 1}$, 
we have  $P^{-1} ( A(T) - I_{n + 1} )^2 P = (B(T) - I_{n + 1})^2 = O_{n  + 1}$, which implies $(A(T) - I_{n + 1})^2 = O_{n + 1}$. 
This contradicts $A(T) \in \mf{J}$. 
\QED

By assertion (5) of Lemma 1.30, we have the following: 

\begin{lem} 
Let $i_1, i_2, i_3, j_1, j_2, j_3$ be integers satisfying $i_1, i_3, j_1, j_3 \geq 1$, $i_2, j_2 \geq 0$, and $i_1+ i_2 + i_3 = j_1 + j_2 + j_3$. 
Let $n := i_1 + i_2 + i_3$. Assume $(i_1, i_2, i_3) \ne (j_1, j_2, j_3)$. 
Then two subsets $\mf{A}(i_1, i_2, i_3)^\circ$ and $\mf{A}(j_1, j_2, j_3)^\circ$ are mutually $GL(n, k)$-disjoint. 
\end{lem}

\subsection{Equivalence relations of polynomial matrices}

\subsubsection{$\mf{J}_{[n]}$}

In this subsubsection, we consider an equivalence relation of two polynomial matrices of $\mf{J}_{[n]}$ (see Theorem 1.36),  
and give a necessary and sufficient condition for two exponential matrices of $\mf{J}_{[n]}^E$ to be equivalent (see Corollary 1.38). 
For stating these, we prepare notations $J(f_1, \ldots, f_n)$, $\cQ_{[n]}(R)$ and $\cP_{[n]}(R)$. 

Let $R$ be a commutative ring. 
For $f_1, \ldots, f_n \in R$, 
let 
\[
J(f_1, \ldots, f_n) := \sum_{i = 0}^n f_i \,\nu_n^i \in \Mat(n, R) .  
\]

We note the following basic lemma: 

\begin{lem}
Let $f_1, \ldots, f_n, g_1, \ldots, g_n \in R$, and let $h_\ell := \sum_{i = 1}^\ell f_i \cdot g_{\ell + 1 - i}$ for $1 \leq \ell \leq n$. 
 we have 
\[
J(f_1, \ldots, f_n) \cdot J(g_1, \ldots, g_n) = J(g_1, \ldots, g_n) \cdot J(f_1, \ldots, f_n) 
 =J(h_1, \ldots, h_n). 
\] 
In particular if $R$ is a field and $f_1 \ne 0$, $J(f_1, \ldots, f_n)$ is regular and its inverse matrix $J(f_1, \ldots ,f_n)^{-1}$ 
can be expressed as $J(g_1', \ldots, g_n')$,  
where $g_1' = 1/f_1$ and $g_\ell' := (f_2 \cdot g_{\ell - 1}' + f_3 \cdot g_{\ell - 2}' + \cdots + f_\ell \cdot g_1')/f_1$ 
for all $2 \leq \ell \leq n$. 
\end{lem}

Given elements $y_1, \ldots, y_n$ of $R$, we can inductively define a sequence of matrices 
$Q_\ell(y_1, \ldots, y_n) \in \Mat(\ell, R)$ $(1 \leq \ell \leq n)$ as 
\begin{eqnarray*}
\left\{
\begin{array}{lcl}
Q_1(y_1, \ldots, y_n) &:=& (y_1) ,  \\
Q_\ell(y_1, \ldots, y_n)  & := & 
 \left(
 \begin{array}{c | c}
 1 & \bm{0} \\
 \hline 
 \bm{0} & Q_{\ell - 1}(y_1, \ldots ,y_n)
 \end{array}
 \right) J(y_1, \ldots, y_\ell ) \qquad (2 \leq \ell \leq n). 
\end{array}
\right. 
\end{eqnarray*}
Clearly, each $Q_\ell(y_1, \ldots, y_n)$ is an upper triangular matrix, 
the first row of $Q_\ell(y_1, \ldots, y_n)$ is equal to $(y_1, \ldots, y_\ell)$, and 
the diagonal $(i, i)$-th entry $q_{i, i}$ of $Q_\ell(y_1, \ldots , y_n)$ coincides with $y_1^i$ for any $1 \leq i \leq \ell$. 
In particular when $R$ is a field, $Q_n(y_1, \ldots, y_n)$ is regular if and only if $y_1 \ne 0$. 

For any $n \geq 2$, we denote by $\cQ_{[n]}(R)$ the set of all regular matrices $Q$ of $GL(n -1, R)$ such that 
$Q = Q_{n - 1}(y_1, \ldots, y_{n - 1})$ for some $y_1, \ldots, y_{n -1} \in k$. 
We write $Q_{[n]}$ in place of $Q_{[n]}(k)$.

For a matrix $A \in \Mat(n, R)$ and for integers $1 \leq i, j \leq n$, 
we denote by $A_{(i, j)}$ the submatrix of $A$ formed by deleting the $i$-th row and the $j$-th column of $A$.

\begin{lem} 
Let $x_1, \ldots, x_n, y_1, \ldots, y_n, z_1, \ldots, z_n \in R$. 
Then the following assertions {\rm (1), (2), (3), (4), (5), (6), (7)} hold true: 
\begin{enumerate}
\item[\rm (1)] We have $Q_\ell(y_1, \ldots, y_n)_{(\ell, \ell)} =  Q_{\ell - 1}(y_1, \ldots, y_n)$.

\item[\rm (2)] Let $Z_\ell := Q_\ell(z_1, \ldots, z_n)$ for $1 \leq \ell \leq n$. 
For all $2 \leq \ell \leq n$, we have 
\[
J(y_1, \ldots, y_\ell) \cdot Z_\ell =  Z_\ell \cdot J(y_1, (y_2, \ldots, y_\ell) \cdot Z_{\ell - 1} ). 
\]

\item[\rm (3)] For all $Q \in \cQ_{[n]}(R)$, 
we have 
\[
\left(
\begin{array}{c | c}
 1 & \bm{0} \\
\hline 
\bm{0} & Q
\end{array}
\right)^{-1}
\cdot 
J(x_1, \ldots, x_n)
\cdot 
\left(
\begin{array}{c | c}
 1 & \bm{0} \\
\hline 
\bm{0} & Q
\end{array}
\right) 
= J\left(
(x_1, \ldots, x_n) 
\left(
\begin{array}{c | c}
 1 & \bm{0} \\
\hline 
\bm{0} & Q
\end{array}
\right)
\right) .
\]

\item[\rm (4)]  The $(i, j)$-th entry $q_{i, j}$ of 
$Q_\ell(y_1, \ldots, y_n)$ can be written as 
\[
 q_{i, j} = \sum_{\lambda_1 + \cdots + \lambda_i = j} y_{\lambda_1} \cdots y_{\lambda_i} . 
\]

\item[\rm (5)] For any $n \geq 2$, $\cQ_{[n]}(R)$ is a subgroup of $GL(n - 1, R)$. 

\item[\rm (6)] Let $R^\times$ be the group of all invertible elements of $R$ with respect to the multiplication. 
We define a group homomorphisms $\pi :  \cQ_{[n]}(R) \to R^\times$ as  
$\pi(Q_{n- 1}(y_1, \ldots, y_{n - 1})) := y_1$ and let $\mathbb{U}_{[n]}(R) := {\rm Ker}(\pi)$ be the kernel of $\pi$. 
Then the exact sequence 
\[
 e \longrightarrow \mathbb{U}_{[n]}(R) \longrightarrow \cQ_{[n]}(R) \overset{\pi}{\longrightarrow} (R^\times, \cdot) \longrightarrow e  
\] 
has a right split map $\iota :  R^\times \to\cQ_{[n]}(R)$ defined by $\iota(y_1) := Q_{n - 1}(y_1, 0, \ldots, 0)$. 
So, $\cQ_{[n]}(R)$ is the semidirect product $\mathbb{U}_{[n]}(R) \rtimes  R^\times$ of 
$\mathbb{U}_{[n]}(R)$ with $R^\times$. 

\item[\rm (7)] For any $y_{n}' \in R$, 
we have 
\[
Q_n(y_1, \ldots, y_n) + J(0, \ldots, 0, y_n') = Q_n(y_1, \ldots, y_{n - 1}, y_n + y_n') \in \cQ_{[n]}. 
\] 
\end{enumerate}
\end{lem}

\Proof (1) We proceed by induction on $\ell$. The proof on the first step $\ell = 2$ is clear. 
Let $3 \leq \ell \leq n$. By the induction hypothesis, we have  
\[
Q_\ell(y_1, \ldots, y_n)_{(\ell, \ell)} = \left(
 \begin{array}{c | c}
 1 & \bm{0} \\
 \hline 
 \bm{0} & Q_{\ell - 2}(y_1, \ldots, y_n)
 \end{array}
 \right) J(y_1, \ldots, y_{\ell - 1})
 = Q_{\ell - 1}(y_1, \ldots, y_n) . 
\]

(2) By deforming the both sides of the equality to be proven, we obtain 
\begin{eqnarray*}
\text{(the left hand side)} 
 &=& J(y_1, \ldots, y_\ell) \cdot 
\left(
 \begin{array}{c | c}
 1 & \bm{0} \\
 \hline 
 \bm{0} & Z_{\ell - 1}
 \end{array}
 \right) \cdot J(z_1, \ldots, z_\ell ) \\
&=& 
\left(
 \begin{array}{c | c}
 y_1 & (y_2, \ldots, y_\ell ) \cdot Z_{\ell - 1} \\
 \hline 
 \bm{0} & J(y_1, \ldots, y_{\ell -1}) \cdot Z_{\ell - 1}
 \end{array}
 \right) \cdot J(z_1, \ldots, z_\ell ) , \\
\text{(the right hand side)}
&=& 
 \left(
 \begin{array}{c | c}
 1 & \bm{0} \\
\hline 
 \bm{0} & Z_{\ell - 1}
 \end{array}
 \right)
 \cdot 
 J(z_1, \ldots, z_\ell) \cdot J(y_1, (y_2, \ldots, y_\ell) Z_{\ell - 1}) \\
&=& 
 \left(
 \begin{array}{c | c}
 y_1 & (y_2, \ldots, y_\ell) Z_{\ell - 1} \\
 \hline 
 \bm{0} &  Z_{\ell - 1} \cdot J(y_1, (y_2, \ldots, y_{\ell - 1}) \cdot Z_{\ell - 2}) 
 \end{array}
 \right)
 \cdot J(z_1, \ldots, z_\ell). 
\end{eqnarray*}
Using the induction on $\ell$, we have  
$J(y_1, \ldots, y_{\ell -1}) \cdot Z_{\ell - 1} = Z_{\ell - 1} \cdot J(y_1, (y_2, \ldots, y_{\ell - 1}) \cdot Z_{\ell - 2})$, 
which implies the desired equality.

(3) Since
\[
 J(x_1, \ldots ,x_n) \cdot 
\left(
\begin{array}{c | c}
 1 & \bm{0} \\
\hline 
\bm{0} & Q
\end{array}
\right)
= 
\left(
\begin{array}{c | c}
 x_1 & (x_2, \ldots, x_n) Q \\
 \hline 
 \bm{0} & J(x_1, \ldots, x_{n - 1})Q 
\end{array}
\right)  
\]
and 
\[
\left(
\begin{array}{c | c}
 1 & \bm{0} \\
\hline 
\bm{0} & Q
\end{array}
\right)
\cdot 
J(x_1, (x_2, \ldots, x_n) Q) 
= 
\left(
\begin{array}{c | c}
 x_1 & (x_2, \ldots, x_n) Q \\
\hline 
 \bm{0} & Q \cdot J(x_1, (x_2, \ldots , x_{n - 1}) Q_{(n - 1, n - 1)} )
\end{array}
\right), 
\]
using the above assertions (2) and (1), we have the desired equality.

(4) If $n = 1$, the assertion is clear. 
Assume $n \geq 2$. We proceed by induction on $\ell$. The proof for the first step $\ell = 1$ is clear. 
For $2 \leq \ell \leq n$, we can write $Q_\ell(y_1, \ldots, y_n)$ as 
\begin{eqnarray*}
Q_\ell(y_1, \ldots, y_n)  
= 
 \left(
 \begin{array}{c | c}
 1 & \bm{0} \\
 \hline 
 \bm{0} & Q_{\ell - 1}(y_1, \ldots, y_n)
 \end{array}
 \right) J(y_1, \ldots, y_\ell ) .  
\end{eqnarray*}
By the induction hypothesis, each $(i, j)$-th entry $q_{i, j}$ of $Q_\ell(y_1, \ldots, y_n)_{(\ell, \ell)} ( = Q_{\ell - 1}(y_1, \ldots, y_n)$, 
see assertion (1) of Lemma 1.34) 
has the desired form. We have only to show that the $(i, \ell)$-th entries $q_{i, \ell}$ $(1 \leq i \leq \ell)$ of 
$Q_\ell(y_1, \ldots, y_n)$ 
have the desired forms. In fact, 
\[
 q_{i, \ell}
 = \sum_{j = 1}^{\ell - 1} q_{i - 1, j} \cdot y_{\ell - j}
 = \sum_{j = 1}^{\ell - 1}
     \left(
          \sum_{\lambda_1 + \cdots + \lambda_{i - 1= j}}       y_{\lambda_1} \cdots y_{\lambda_{i - 1}}
     \right)
      \cdot y_{\ell - j}
 =  \sum_{\lambda_1 + \cdots + \lambda_i = \ell}       y_{\lambda_1} \cdots y_{\lambda_{i - 1}}  \cdot y_{\lambda_i} . 
\]

(5) We first prove that $Q_{\rm I}, Q_{\rm II} \in \cQ_{[n]}(R) \Longrightarrow Q_{\rm I} \cdot Q_{\rm II} \in \cQ_{[n]}(R)$. 
We can write $Q_{\rm I}$ and $Q_{\rm II}$ as 
$Q_{\rm I} = Q_{n - 1}(y_1, \ldots, y_{n - 1})$ and $Q_{\rm II} = Q_{n-1}(z_1, \ldots, z_{n -1})$ 
for some $y_1, \ldots, y_{n- 1}, z_1, \ldots, z_{n - 1} \in R$. Let $Y_\ell := Q_\ell(y_1, \ldots, y_{n - 1})$ and 
$Z_\ell := Q_\ell(z_1, \ldots, z_{n - 1})$ for $1 \leq \ell \leq n - 1$. 
For $2 \leq \ell \leq n - 1$, by the above assertion (2) and Lemma 1.33, we have 
\begin{eqnarray*}
Y_\ell \cdot Z_\ell 
 & =  &
 \left(
 \begin{array}{c | c}
 1 & \bm{0} \\
 \hline 
 \bm{0} & Y_{\ell - 1}
 \end{array}
 \right)
 J(y_1, \ldots, y_\ell) \cdot Z_\ell \\
 &\overset{\rm (2)}{ = } &  
 \left(
 \begin{array}{c | c}
 1 & \bm{0} \\
 \hline 
 \bm{0} & Y_{\ell - 1}
 \end{array}
 \right) 
 \cdot Z_\ell \cdot J(y_1, (y_2, \ldots, y_{\ell - 1}) Z_{\ell -1})  \\
 & = &  
 \left(
 \begin{array}{c | c}
 1 & \bm{0} \\
 \hline 
 \bm{0} & Y_{\ell - 1} \cdot Z_{\ell - 1}
 \end{array}
 \right) 
 \cdot J(z_1, \ldots,  z_\ell) 
 \cdot J(y_1, (y_2, \ldots, y_{\ell - 1}) Z_{\ell -1}) \\
 & = & 
  \left(
 \begin{array}{c | c}
 1 & \bm{0} \\
 \hline 
 \bm{0} & Y_{\ell - 1} \cdot Z_{\ell - 1}
 \end{array}
 \right) \cdot 
J(y_1 z_1, \; y_1 (z_2, \ldots, z_\ell) + (y_2, \ldots, y_\ell) \cdot Z_{\ell - 1} \cdot J(z_1, \ldots, z_{\ell - 1})) . 
\end{eqnarray*} 
Thus the sequence $\{ Y_\ell \cdot Z_\ell \}_{1 \leq \ell \leq n - 1}$ satisfies the recurrence formula, which implies 
$Q_{\rm I} \cdot Q_{\rm II} \in \cQ_{[n]}(R)$.

We next prove $Q \in \cQ_{[n]}(R) \Longrightarrow Q^{-1} \in\cQ_{[n]}(R)$. 
We can write $Q$ as $Q = Q_{n - 1}(y_1, \ldots, y_{n - 1})$ for some $y_1, \ldots, y_{n - 1} \in k$. 
Let $Y_\ell : = Q_\ell(y_1, \ldots, y_{n - 1})$ for $1 \leq \ell \leq n - 1$. 
We inductively define $y_1', y_2', \ldots, y_{n - 1}' \in k$ and $Y_i' \in \Mat(i, k)$ $(1 \leq i \leq n -1)$ as follows: 
Let $y_1' := 1/y_1$ and let $Y_1' := (y_1') \in \Mat(1, k)$. 
Given $y_i' \in k$ and $Y_i' \in \Mat(i, k)$ ($1 \leq i\leq \ell - 1$), 
we can choose $y_{\ell}' \in k$ such that   
\[
 Y_\ell \cdot 
 \left(
 \begin{array}{c | c}
 1 & \bm{0} \\
 \hline 
 \bm{0} & Y_{\ell - 1}
 \end{array}
 \right)
 \left(
 \begin{array}{c}
  y_\ell' \\
 \vdots \\
 y_2' \\
 y_1'
 \end{array}
 \right)
 = 
 \left(
 \begin{array}{c}
  0 \\
 \vdots \\
  0 \\
  1
 \end{array}
 \right) , 
\]
and define 
\[
 Y_\ell' := 
 \left(
 \begin{array}{c | c}
 1 & \bm{0} \\
 \hline 
 \bm{0} & Y_{\ell - 1}' 
 \end{array}
 \right)
 \cdot J(y_1', \ldots, y_\ell') . 
\]
Clearly, $Y'_{n - 1} \in \cQ_{[n]}(R)$. 
We know from the above assertion (1) that 
$(Y_\ell \cdot Y_\ell')_{(1, 1)} = (Y_\ell)_{(1, 1)} \cdot (Y_\ell')_{(1, 1)} = Y_{\ell - 1} \cdot Y_{\ell - 1}'$. 
By the definitions of $y_\ell'$ and $Y_\ell'$, the $\ell$-th column of $Y_\ell \cdot Y_\ell'$ is equal to ${^t}(0, \ldots, 0, 1)$. 
Using the induction on $\ell$, we can prove $Y_\ell \cdot Y_\ell' = I_\ell$ for all $1 \leq \ell \leq n - 1$, 
which implies $Q \cdot Y'_{n - 1} = I_{n - 1}$.

(6) The proof is straightforward.

(7) Note that $y_n$ does not appear in any entry of $Q_{n - 1}(y_1, \ldots, y_n)$ (cf. the above assertion (4)). 
So, $Q_{n - 1}(y_1, \ldots, y_n) = Q_{n - 1}(y_1, \ldots, y_{n - 1}, y_n + y_n')$ and 
\begin{eqnarray*}
\lefteqn{ Q_n(y_1, \ldots, y_n) + J(0, \ldots, 0, y_n') }\\
 &=& 
 \left(
 \begin{array}{c | c}
 1 & \bm{0} \\
 \hline 
 \bm{0} & Q_{n - 1}(y_1, \ldots, y_n)
 \end{array}
 \right)
 J(y_1, \ldots, y_n)
 + J(0, \ldots, 0, y_n') \\
 &=& 
 \left(
 \begin{array}{c | c}
 1 & \bm{0} \\
 \hline 
 \bm{0} & Q_{n - 1}(y_1, \ldots, y_n)
 \end{array}
 \right)
 J(y_1, \ldots, y_{n - 1}, y_n +  y_n') 
 = Q_n(y_1, \ldots, y_{n - 1}, y_n + y_n') . 
\end{eqnarray*}
\QED

For any integer $n \geq 2$, any matrix $Q \in \Mat(n - 1, R)$ and arbitrary elements $x_1, \ldots, x_n \in R$, let 
\[
 \varpi(Q \mid x_1, \ldots, x_n ) :=   
 \left(
 \begin{array}{c | c}
 1 & \bm{0} \\
 \hline 
 \bm{0} & Q
 \end{array}
 \right) 
 J(x_1, \ldots, x_n) . 
\]
We note that $(x_1, \ldots, x_n)$ is the first row of $\varpi(Q \mid x_1, \ldots, x_n)$. 
Clearly, $Q_\ell(y_1, \ldots , y_n) = \varpi(Q_{\ell - 1}(y_1, \ldots ,y_n) \mid y_1, \ldots, y_{\ell })$ for all $2 \leq \ell \leq n$.

For any $n \geq 2$, we denote by $\cP_{[n]}(R)$ the set of all matrices $P$ of $GL(n, R)$ 
with the form $P = \varpi(Q \mid x_1, \ldots, x_n)$, where $Q \in \cQ_{[n]}(R)$ and $x_1, \ldots, x_n \in R$. 
We simply write $\cP_{[n]}$ in place of $\cP_{[n]}(k)$.

\begin{lem} 
The following assertions {\rm (1)} and {\rm (2)} hold true: 
\begin{enumerate}
\item[\rm (1)] $\cP_{[n]}(R)$ is a subgroup of $GL(n, R)$. 
\item[\rm (2)] We can define a group homomorphism $\pi : \cP_{[n]}(R) \to \cQ_{[n]}(R)$ as 
$\pi( \varpi(Q \mid x_1, \ldots, x_n)) := Q$. Let ${\cal K}_{[n]}(R) := \Ker(\pi)$ be the kernel of $\pi$. 
Then the exact sequence 
\[
 e \longrightarrow {\cal K}_{[n]}(R) \longrightarrow \cP_{[n]}(R) \overset{\pi}{\longrightarrow} \cQ_{[n]}(R) \longrightarrow e 
\] 
has a right split map $\iota : \cQ_{[n]}(R) \to \cP_{[n]}(R)$ defined by $\iota(Q) := \varpi(Q \mid 1, 0, \ldots, 0)$, 
and the kernel ${\cal K}_{[n]}(R)$ is isomorphic to the group $R^\times \times R^{\oplus (n - 1)}$.  
So, $\cP_{[n]}(R)$ is isomorphic to the semidirect product $(R^\times \times R^{\oplus (n - 1)}) \rtimes \cQ_{[n]}(R)$ of
$R^\times \times R^{\oplus (n - 1)}$ with $\cQ_{[n]}(R)$. 
\end{enumerate} 
\end{lem}

\Proof (1) We first prove $P, P' \in \cP_{[n]}(R) \Longrightarrow P \cdot P' \in \cP_{[n]}(R)$. 
We can write $P, P'$ as $P = \varpi(Q \mid x_1, \ldots, x_n)$ and $P' = \varpi(Q' \mid x_1', \ldots, x_n')$. 
So, by assertion (3) of Lemma 1.34,  
\begin{eqnarray*}
P \cdot P' 
 &=& 
\left(
\begin{array}{c | c}
 1 & \bm{0} \\
\hline 
\bm{0} & Q
\end{array}
\right)
\left(
\begin{array}{c | c}
 1 & \bm{0} \\
\hline 
\bm{0} & Q'
\end{array}
\right)
\cdot 
\left(
\begin{array}{c | c}
 1 & \bm{0} \\
\hline 
\bm{0} & Q'
\end{array}
\right)^{-1} 
J(x_1, \ldots, x_n)
\left(
\begin{array}{c | c}
 1 & \bm{0} \\
\hline 
\bm{0} & Q'
\end{array}
\right)
\cdot 
J(x_1', \ldots, x_n') \\ 
&=& 
\left(
\begin{array}{c | c}
 1 & \bm{0} \\
\hline 
\bm{0} & Q \cdot Q'
\end{array}
\right) 
\cdot 
J(x_1, (x_2, \ldots, x_n) Q')
\cdot 
J(x_1', \ldots, x_n'), 
\end{eqnarray*}
which implies $P \cdot P' \in \cP_{[n]}$ (see Lemma 1.33). 

We next prove $P \in \cP_{[n]} \Longrightarrow P^{-1} \in \cP_{[n]}$. 
We can write $P$ as $P = \varpi(Q \mid x_1, \ldots, x_n)$ for some $x_1, \ldots, x_n \in k$ with $x_1 \ne 0$.  
So, 
\[
 P^{-1} = 
 \left(
 \begin{array}{c | c}
 1 & \bm{0} \\
 \hline 
 \bm{0} & Q^{-1} 
 \end{array}
 \right)
  \cdot 
 \left(
 \begin{array}{c | c}
 1 & \bm{0} \\
 \hline 
 \bm{0} & Q 
 \end{array}
 \right)
J(x_1, \ldots, x_n)^{-1} \cdot 
 \left(
 \begin{array}{c | c}
 1 & \bm{0} \\
 \hline 
 \bm{0} & Q 
 \end{array}
 \right)^{-1} , 
\]
and then by Lemma 1.33 and assertions (5) and (3) of Lemma 1.34, $P^{-1} \in \cP_{[n]}$.

(2) The calculation of $P \cdot P'$ in the proof of the above assertion (2) implies $\pi$ is a group homomorphism. 
Now, we can prove the assertion. 
\QED

\begin{thm}
Let $n \geq 2$ be an integer. 
Let $A(T)$ and $B(T)$ be upper triangular polynomial matrices of $\Mat(n, k[T])$ with the following forms:  
\begin{eqnarray*}
\left\{ 
\begin{array}{rcl l }
A(T) &=& \ds \sum_{i  = 1}^{n - 1} a_i \nu_n^i  \quad & (\text{$a_1, \ldots, a_{n - 1} \in k[T]$ are linearly independent over $k$} ), \\ [5mm] 
B(T) &=& \ds \sum_{i = 1}^{n - 1} b_i \nu_n^i  & (\text{$b_1, \ldots, b_{n - 1} \in k[T]$ are linearly independent over $k$}) . 
\end{array}
\right. 
\end{eqnarray*} 
Then the following conditions {\rm (1)}, {\rm (2)} and {\rm (3)} are equivalent: 
\begin{enumerate}
\item[\rm (1)] $A(T)$ and $B(T)$ are equivalent. 
\item[\rm (2)] We have $P^{-1} A(T) P = B(T)$ and $(b_1, \ldots, b_{n - 1}) = (a_1, \ldots, a_{n - 1}) Q$ 
for some $P \in \cP_{[n]}$ and $Q \in \cQ_{[n]}$.  
\item[\rm (3)] There exists a matrix $Q \in \cQ_{[n]}$ such that $(b_1, \ldots, b_{n - 1}) = (a_1, \ldots, a_{n - 1}) Q$. 
\end{enumerate}
\end{thm}

\Proof The implication (2) $\Longrightarrow$ (1) is clear. We first prove (3) $\Longrightarrow$ (2). 
The matrix $Q$ can be written as $Q = Q_{n - 1}(y_1, \ldots, y_{n - 1})$, where $y_1, \ldots, y_{n - 1} \in k$. 
Let $P := \varpi(Q \mid 1, 0, \ldots, 0) \in \cP_{[n]}$.  
By assertion (3) of Lemma 1.34, $P^{-1} \cdot A(T) \cdot P = B(T)$.

We next prove (1) $\Longrightarrow$ (3).  
We know from Lemma 1.19 that $P \in B(n, k)$. 
Since $A(T) P = PB(T)$, we have a unique matrix $Q \in B(n - 1, k)$ such that 
$(b_1, \ldots, b_{n - 1}) = (a_1, \ldots, a_{n - 1}) Q$. 
We proceed by induction on $n$ for showing 
\[
P = \varpi(Q \mid x_1, \ldots, x_n) \qquad \text{ and } \qquad Q \in \cQ_{[n]}, 
\] 
where we let $(x_1, \ldots, x_n)$ be the first row of $P$. 
The proof for the first step $n =2$ is straightforward. So let $n \geq 3$. 
Note that 
\[
P_{(n, n)}^{-1} \cdot A(T)_{(n, n)} \cdot P_{(n, n)} =B(T)_{(n, n)} \qquad \text{ and } \qquad  
(b_1, \ldots, b_{n - 2}) = (a_1, \ldots, a_{n - 2})Q_{(n - 1, n - 1)}. 
\] 
By the induction hypothesis, 
\[
P_{(n, n)} = \varpi(Q_{(n - 1, n - 1)} \mid x_1, \ldots , x_{n - 1}) \qquad \text{ and } \qquad Q_{(n - 1, n - 1)} \in \cQ_{[n - 1]} . 
\]
So, letting $(y_1, \ldots, y_{n - 1})$ be the first row of $Q$, we can write $Q_{(n - 1, n - 1)}$ as 
\[
Q_{(n - 1, n - 1)} = Q_{n - 2}(y_1, \ldots, y_{n - 2}) = Q_{n - 2}(y_1, \ldots , y_{n - 1}). 
\] 
Let $P = (p_{i, j})$. 
Comparing the $(1, n)$-th entries of the both sides of the equality $A(T) P = P  B(T)$, we have 
\[
\left( \begin{array}{c} p_{2, n} \\ \vdots \\ p_{n, n} \end{array}\right) =  Q \left( \begin{array}{c} x_{n - 1} \\ \vdots \\ x_1  \end{array}\right) . 
\]
It follows that 
\[
P = 
\left(
\begin{tabular}{c | c}
\multirow{2}{*}{ $P_{(n, n)}$ }   & $p_{1, n}$ \\    
      & $\vdots$ \\
 \cline {1-1}
 $\bm{0}$ & $p_{n, n}$
\end{tabular}
\right)
 = \varpi(Q \mid x_1, \ldots, x_n). 
\]
Then calculating $(A(T) P)_{(1, 1)} = (P B(T))_{(1, 1)}$, we can obtain 
\[
 J(0, a_1, \ldots, a_{n -2}) \cdot Q = Q \cdot J(0, b_1, \ldots, b_{n - 2}) .  
\]
Now we have 
\[
 Q^{-1} \cdot A(T)_{(n, n)} \cdot Q = B(T)_{(n, n)} \qquad \text{ and } \qquad  
(b_1, \ldots, b_{n - 2}) = (a_1, \ldots, a_{n - 2})Q_{(n - 1, n - 1)}. 
\]
Again by the induction hypothesis, 
\[
Q = \varpi(Q_{(n - 1, n - 1)} \mid y_1, \ldots, y_{n - 1}) 
 = \left(
 \begin{array}{c | c}
 1 & \bm{0} \\
 \hline 
\bm{0} & Q_{n -2}(y_1, \ldots, y_{n - 1}) 
 \end{array} 
\right) 
 J(y_1, \ldots, y_{n -1}) . 
\]
Hence $Q = Q_{n - 1}(y_1, \ldots, y_{n -1}) \in \cQ_{[n]}$. 
\QED

As a corollary of Theorem 1.36, we have a necessary and sufficient condition for two exponential matrices 
of $\mf{J}_{[n]}^E$ to be equivalent (see Corollary 1.38). We shall use the following Lemma 1.37 on proving Corollary 1.38.

For an integer $i \geq 0$, let $i = d_0  + d_1 p + d_2 p + \cdots + d_r p^r$ be the $p$-adic expansion of $i$. 
Let $\sigma_p(i) = \sum_{j = 0}^r d_j $ be the sum of digits $d_0, \ldots, d_r$.

\begin{lem}
The following assertions hold true: 
\begin{enumerate}
\item[\rm (1)] For any $1 \leq d \leq p-1$, we have $\sigma_p(\sum_{\lambda = 1}^d p^{e_\lambda}) = d$. 
\item[\rm (2)] For a monomial $m = T^i\in k[T]$ $(i > 0)$, $m $ is a $p$-polynomial if and only if $\sigma_p(i) = 1$. 
\item[\rm (3)] Let $\varphi_1, \ldots, \varphi_d$ $(2 \leq d \leq p - 1)$ be non-zero $p$-polynomials of $k[T]$. 
Then neither monomials appearing in the product $\varphi_1 \cdots \varphi_d$ are $p$-polynomials. 
\end{enumerate}
\end{lem}

\Proof Assertions (1) and (2) are clear. And assertion (3) follows from the above assertions (1) and (2). 
\QED

\begin{cor}
Let $n \geq 2$ be an integer. 
Let $A(T)$ and $B(T)$ be exponential matrices of $\mf{J}_{[n]}^E$ with the following forms: 
\begin{eqnarray*}
\left\{ 
\begin{array}{ll}
 A(T) = \ds \Exp\left( \sum_{i = 1}^{n - 1} f_i \nu_n^i \right) 
\quad & 
(\; \text{$f_i$ $(1 \leq i \leq n - 1)$ are $p$-polynomials satisfying $f_1 \ne 0$} \; ), \\ [7mm] 
 B(T) = \ds \Exp\left( \sum_{i = 1}^{n - 1} g_i \nu_n^i \right) 
 & 
(\; \text{$g_i$ $(1 \leq i \leq n -1)$ are $p$-polynomials satisfying $g_1 \ne 0$} \; ). 
\end{array} 
\right. 
\end{eqnarray*}
Then the following conditions {\rm (1)} and {\rm (2)} are equivalent: 
\begin{enumerate}
\item[\rm (1)] $A(T)$ and $B(T)$ are equivalent. 
\item[\rm (2)] There exists a matrix $Q \in \cQ_{[n]}$ such that $(g_1, \ldots, g_{n - 1}) = (f_1, \ldots, f_{n - 1}) Q$. 
\end{enumerate}
\end{cor}

\Proof For any $0 \leq \ell \leq n - 1$, we define polynomials $a_\ell$ and $b_\ell$ as 
\[
a_\ell : = 
  \sum_{\substack{i_1 +2 i_2 + \cdots + (n - 1) i_{n - 1} = \ell \\ i_1, i_2, \ldots, i_{n - 1} \geq 0}} 
 \frac{f_1^{i_1} f_2^{i_2} \cdots f_{n - 1}^{i_{n - 1}}}{i_1! i_2! \cdots i_{n - 1}!} 
\quad \text{ and } \quad 
b_\ell : = 
  \sum_{\substack{i_1 +2 i_2 + \cdots + (n - 1) i_{n - 1} = \ell \\ i_1, i_2, \ldots, i_{n - 1} \geq 0}} 
 \frac{g_1^{i_1} g_2^{i_2} \cdots g_{n - 1}^{i_{n - 1}}}{i_1! i_2! \cdots i_{n - 1}!} . 
\]
By Lemma 1.17, we have $A(T) = \sum_{\ell = 0}^{n - 1} a_\ell \nu_n^\ell$ and $B(T) = \sum_{\ell = 0}^{n - 1} b_\ell \nu_n^\ell$. 
By Theorem 1.13, we have $1 \leq n \leq p$. 
Let $m =  T^i$ be a monomial appearing among $a_0, \ldots, a_{n - 1}$. 
If $m$ appears in $f_1^{i_1} f_2^{i_2} \cdots f_{n - 1}^{i_{n - 1}}$, then 
\[
\sigma_p(i) = i_1 + i_2 + \cdots + i_{n - 1}. 
\] 
In fact, we can write $i = p^{e_1} + \cdots + p^{e_d}$, where $d = i_1 + i_2 + \cdots + i_{n-1}$. 
Since $d \leq p - 1$, using assertion (1) of Lemma 1.37, we have $\sigma_p(i) = d$. 
Thus, if $m$ appears in $f_1^{i_1}  f_2^{i_2} \cdots f_{n -1}^{i_{n - 1}}$, where $i_1 + 2 i_2 + \cdots + (n - 1) i_{n - 1} = \ell$, 
then $\sigma_p(i) \leq \ell$, and the equality $\sigma_p(i) = \ell$ holds true if and only if $m$ appears in $f_1^{\ell}$. 
So, there is no monomial $m = T^i$ of $f_1^\ell$ appears in $a_0 , \ldots, a_{\ell - 1}, a_\ell -f_1^\ell$. 
Thus $a_1, \ldots, a_{n - 1}$ are linearly independent over $k$, and $b_1, \ldots, b_{n - 1}$ are also linearly independent over $k$. 
We shall prove (1) $\Longrightarrow$ (2). 
We know from Theorem 1.36 that $(a_1, \ldots, a_{n - 1}) = (b_1, \ldots, b_{n -1}) Q$ for some $Q \in \cQ_{[n]}$. 
Any $p$-monomial appears in $a_\ell$ has to appear in $f_\ell$ (see assertion (3) of Lemma 1.37). 
Thus $(f_1, \ldots, f_{n - 1}) = (g_1, \ldots, g_{n - 1}) Q$.  
Conversely, we shall prove (2) $\Longrightarrow$ (1). 
We know from Theorem 1.36 that there exists a regular matrix $P$ of $GL(n, k)$ such that 
$P^{-1} \left( \sum_{i = 1}^{n - 1} f_i \nu_n^i \right) P = \sum_{i = 1}^{n - 1} g_i \nu_n^i$, which 
implies $P^{-1} A(T) P = B(T)$. 
\QED

\subsubsection{$\mf{J}_{[n, 1]}^0$}

\begin{thm}
Let $n \geq 2$ be an integer. 
Let $A(T)$ and $B(T)$ be polynomial matrices of $\Mat(n + 1, k[T])$ with the forms 
\begin{eqnarray*}
\left\{
\begin{array}{ll}
A(T) = 
\left(
\begin{array}{c | c}
\ds \sum_{i = 1}^{n - 1} a_i \nu_n^i & \bm{0} \\
\hline 
 \bm{0} & 0 
\end{array}
\right)
\quad & 
(\text{$a_1, \ldots, a_{n -1}$ are linearly independent over $k$}), \\ [8mm] 

B(T) = 
\left(
\begin{array}{c | c}
\ds \sum_{i = 1}^{n - 1} b_i \nu_n^i & \bm{0} \\
\hline 
 \bm{0} & 0 
\end{array}
\right)
\quad &
(\text{$b_1, \ldots, b_{n -1}$ are linearly independent over $k$}). 
\end{array}
\right.
\end{eqnarray*}
Then the following conditions {\rm (1)} and {\rm (2)} are equivalent: 
\begin{enumerate}
\item[\rm (1)] $A(T)$ and $B(T)$ are equivalent. 
\item[\rm (2)] There exists a matrix $Q$ of $\cQ_{[n]}$ such that $(b_1, \ldots, b_{n - 1}) = (a_1, \ldots, a_{n -1}) Q$. 
\end{enumerate}
\end{thm}

\Proof (1) $\iff$ there exists a regular matrix $P_1$ of $GL(n, k)$ 
such that $P_1^{-1} \cdot (\sum_{i = 1}^{n - 1} a_i \nu_n^i) \cdot P_1 = \sum_{i = 1}^{n - 1} b_i \nu_n^i \iff$ (2) 
(see Theorem 1.36 for the latter equivalence). 
\QED

As a corollary of Theorem 1.39, we have the following:

\begin{cor}
Let $n \geq 2$ be an integer. 
Let $A(T)$ and $B(T)$ are exponential matrices of $(\mf{J}_{[n, 1]}^0)^E$ with the forms 
\begin{eqnarray*}
\left\{
\begin{array}{ll}
A(T) = 
\left(
\begin{array}{c | c}
\ds \Exp\left( \sum_{i = 1}^{n - 1} f_i \nu_n^i \right) & \bm{0} \\
\hline 
 \bm{0} & 1 
\end{array}
\right)
\quad & 
(\text{$f_1, \ldots, f_{n -1}$ are $p$-polynomials, and $f_1 \ne 0$}), \\ [8mm] 

B(T) = 
\left(
\begin{array}{c | c}
\ds \Exp\left( \sum_{i = 1}^{n - 1} g_i \nu_n^i \right) & \bm{0} \\
\hline 
 \bm{0} & 1 
\end{array}
\right)
\quad & 
(\text{$g_1, \ldots, g_{n -1}$ are $p$-polynomials, and $g_1 \ne 0$}) .  
\end{array}
\right.
\end{eqnarray*}
Then the following conditions {\rm (1)} and {\rm (2)} are equivalent: 
\begin{enumerate}
\item[\rm (1)] $A(T)$ and $B(T)$ are equivalent. 
\item[\rm (2)] There exists a matrix $Q$ of $\cQ_{[n]}$ such that $(g_1, \ldots, g_{n - 1}) = (f_1, \ldots, f_{n -1}) Q$. 
\end{enumerate}
\end{cor}

\subsubsection{$\mf{J}_{[n, 1]}^1$}

In this subsubsection, we consider an equivalence relation of two polynomial matrices of $\mf{J}_{[n, 1]}^1$ (see Theorem 1.41),  
and give a necessary and sufficient condition for two exponential matrices of $(\mf{J}_{[n, 1]}^1)^E$ to be equivalent (see Corollary 1.42). 

Let $n \geq 2$ be an integer and let $\cQ_{[n, 1]}$ be the set of all matrices $Q$ of $\Mat(n, k)$ 
with the form 
\[
Q = 
\left(
\begin{array}{ c | c }
 Q^\flat & u \cdot \bm{e}_1 \\ 
\hline 
 v \cdot \transpose \bm{e}_{n - 1} & w
\end{array}
\right)
\qquad 
(\text{$Q^\flat \in \cQ_{[n]}$, \;  $u, v, w \in k$, \; $w \ne 0$}) . 
\]
We shall show that $\cQ_{[n, 1]}$ becomes a subgroup of $GL(n, k)$. 
By assertions (5) and (7) of Lemma 1.34, we have $Q, Q' \in \cQ_{[n, 1]} \Longrightarrow Q \cdot Q' \in \cQ_{[n, 1]}$. 
For any $Q \in \cQ_{[n, 1]}$ with the above form, 
defining $u', v'\in k$ as $u' := y_1^{-1} \cdot u$ and $v' := w^{-1} \cdot v$, where $y_1$ is the $(1, 1)$-th entry of $Q^\flat$,  
we have 
\[
\left(
\begin{array}{c | c}
 I_n & - u' \cdot \bm{e}_1 \\
 \hline 
 - v' \cdot {^t}\bm{e}_{n - 1} & 1
\end{array}
\right)
\left(
\begin{array}{c | c}
 (Q^\flat)^{-1} & \bm{0} \\
\hline 
 \bm{0} & w^{-1} 
\end{array}
\right)
 \cdot Q
= 
I_{n + 1} , 
\]
which implies $Q$ is regular and $Q^{-1} \in \cQ_{[n, 1]}$.

\begin{thm}
Let $n \geq 2$ be an integer. 
Let $A(T)$ and $B(T)$ be polynomial matrices of $\Mat(n + 1, k[T])$ with the following forms: 
\begin{eqnarray*}
\left\{
\begin{array}{ll}
A(T) = 
\left(
\begin{array}{ c | c }
 \ds \sum_{i = 1}^{n - 1} a_i \nu_n^i & a_n \bm{e}_1 \\
\hline 
 \bm{0} & 0 
\end{array}
\right)
\quad & 
\left(
\begin{tabular}{l}
\text{$a_1, \ldots, a_{n - 1}$ are linearly independent over $k$,} \\
\text{and $a_1$, $a_n$ are also linearly independent over $k$}
\end{tabular} 
\right) , 
\\ [10mm] 
B(T) = 
\left(
\begin{array}{ c | c }
 \ds \sum_{i = 1}^{n - 1} b_i \nu_n^i & b_n \bm{e}_1 \\
\hline 
 \bm{0} & 0 
\end{array}
\right)
\quad & 
\left(
\begin{tabular}{l}
\text{$b_1, \ldots, b_{n - 1}$ are linearly independent over $k$,} \\
\text{and $b_1$, $b_n$ are also linearly independent over $k$}
\end{tabular} 
\right) . 
\end{array}
\right. 
\end{eqnarray*}
Then the following conditions {\rm (1)} and {\rm (2)} are equivalent: 
\begin{enumerate}
\item[\rm (1)] $A(T)$ and $B(T)$ are equivalent. 
\item[\rm (2)] There exists a matrix $Q$ of $\cQ_{[n, 1]}$ such that $(b_1, \ldots, b_n) = (a_1, \ldots, a_n) Q$. 
\end{enumerate}
\end{thm}

\Proof 
Let $\alpha_1 := \sum_{i = 1}^{n - 1} a_i \nu_n^i$, $\alpha_2 := a_n \bm{e}_1$, $\beta_1 := \sum_{i = 1}^{n - 1} b_i \nu_n^i$, 
and $\beta_2 := b_n \bm{e}_1$. 

We first prove (1) $\Longrightarrow$ (2). There exists a regular matrix $P$ of $GL(n + 1, k)$ such that $P^{-1} A(T) P = B(T)$. 
We can write $P$ as 
\[
P = \left( 
\begin{array}{c | c}
 P_1 & P_2 \\
\hline 
 P_3 & P_4 
\end{array}
\right) \quad
(\text{$P_1 \in \Mat(n, k)$, $P_2 \in \Mat_{n, 1}(k)$, $P_3 \in \Mat_{1, n}(k)$, $P_4 \in \Mat(1, k)$}).
\]
Since $A(T) P = P B(T)$, we have 
\begin{eqnarray*}
(\ast) \qquad 
 \alpha_1 P_1 + \alpha_2 P_3 = P_1 \beta_1 , \quad\; 
 \alpha_1 P_2 + \alpha_2 P_4 = P_1 \beta_2 , \quad \;
 O = P_3 \beta_1 , \quad\; 
 O = P_3 \beta_2 . 
\end{eqnarray*}
We shall consider the forms of $P_1, P_2$ and $P_3$. 
By the third equality of $(\ast)$, we have $P _3 = (0, \ldots, 0, c_3)$ for some $c_3 \in k$. 
The fourth equality of $(\ast)$ holds true (since $n \geq 3$). 
By the first equality of $(\ast)$ and by Lemma 1.19, we know that $P_1$ is an upper triangular matrix. 
By the second equality of $(\ast)$, the $i$-th entries of $\alpha_1 P_2$ are zeroes for all $2 \leq i \leq n$. 
Since $a_1, \ldots, a_{n - 2}$ are linearly independent over $k$, we can write $P_2$ as 
\[
P_2 = \left(
\begin{array}{c}
 c_{21} \\
 c_{22} \\
 \bm{0}
\end{array}
\right)
\qquad \text{for some $c_{21}, c_{22} \in k$} . 
\]
So, by the forms of $P_1, P_2, P_3$, we know that $P_1$ is regular (since $P$ is regular). 
Now, the first equality of $(\ast)$ implies   
\begin{eqnarray*}
P_1^{-1} \alpha_1 P_1
 = \beta_1 - P_1^{-1} \alpha_2 P_3
 = \beta_1 -  
\left(
\begin{array}{c | c}
\bm{0} & (a_n c_3)/c_1 \\
\hline 
O & \bm{0}
\end{array}
\right), 
\end{eqnarray*}
where $c_1$ is the $(1, 1)$-th entry of $P_1$. We know from Theorem 1.36 that 
\[
 (b_1, \ldots, b_{n - 2}, b_{n - 1} - (a_n c_3)/c_1 ) = (a_1, \ldots, a_{n - 1}) Q \qquad \text{for some $Q \in \cQ_{[n]}$}. 
\]
And the second equality of $(\ast)$ implies $b_n = (c_{22}/c_1) a_1 + (P_4/c_1) a_n$. 
Thus, we have 
\[
(b_1, \ldots, b_{n -1}, b_n)
= (a_1, \ldots, a_{n - 1}, a_n)
\left(
\begin{array}{c | c}
   Q                                         &  (c_{22}/c_1) \cdot \bm{e}_1\\
\hline 
(c_3/c_1) \cdot \transpose\bm{e}_{n - 1} & P_4/c_1
\end{array}
\right) . 
\]
We next prove (2) $\Longrightarrow$ (1).
Let 
\[
P := 
\left(
\begin{array}{c | c }
 1 & \bm{0}  \\ 
\hline 
\bm{0} & Q 
\end{array}
\right)  \in GL(n + 1, k) . 
\]
Write 
\begin{eqnarray*}
\left\{
\begin{array}{l}
Q = 
\left(
\begin{array}{ c | c }
 Q^\flat & u \cdot \bm{e}_1 \\ 
\hline 
 v \cdot \transpose \bm{e}_{n - 1} & w
\end{array}
\right)
\qquad 
(\text{$Q^\flat \in \cQ_{[n]}$, \;  $u, v, w \in k$, \; $w \ne 0$}) , \\
Q^\flat = (q_{i, j}) .
\end{array}
\right. 
\end{eqnarray*}
Let $Q^{\flat\flat} := (q_{i, j})_{1 \leq i, j \leq n - 2}$ be the submatrix of $Q^\flat$ and 
let $\bm{q}$ be the $(n - 1)$-th column of $Q^\flat$. 
Since $(b_1, \ldots, b_{n - 2}, b_{n - 1}, b_n) = (a_1, \ldots, a_{n - 2}, a_{n - 1}, a_n) Q$, we have 
\begin{eqnarray*}
\left\{
\begin{array}{l}
 (b_1, \ldots, b_{n - 2}) = (a_1, \ldots, a_{n - 2}) Q^{\flat\flat}, \\
 b_{n - 1} = (a_1, \ldots, a_{n - 1}) \bm{q} + a_n v, \\
 b_n = a_1 u + a_n w .
\end{array}
\right.
\end{eqnarray*}
The first row of $A(T) \cdot P$ can be calculated as 
\begin{eqnarray*}
\lefteqn{( \, 0,\; ( a_1, \ldots, a_{n - 1}) Q^\flat + (0, \ldots, 0, a_n v), \; a_1 u + a_n w \, ) } \\
 &=& (\, 0, \;  (a_1, \ldots, a_{n - 2}) Q^{\flat\flat}, \; (a_1, \ldots, a_{n - 1}) \bm{q} + a_n v, \; a_1 u + a_n w \,) \\
 &=& (0, b_1, \ldots, b_{n - 2}, b_{n - 1}, b_n) . 
\end{eqnarray*}
Consider the $(1, 1)$-th block of $A(T) \cdot P$ of size $n \times n$. By assertion (2) of Lemma 1.34, the block is equal to 
the $(1, 1)$-th block of $P \cdot B(T)$ of size $n \times n$. In fact,  
\begin{eqnarray*}
\alpha_1 \cdot 
\left(
 \begin{array}{c | c}
 1 & \bm{0} \\
 \hline 
 \bm{0} & Q^\flat
 \end{array}
\right)
 &=& 
\left(
\begin{array}{c | c}
 0 & b_1, \ldots\ldots, b_{n -2}, b_{n - 1} \\
\hline 
\bm{0} & J(0, a_1, \ldots, a_{n - 2}) \cdot Q^\flat 
\end{array} 
\right) \\
& = &
\left(
\begin{array}{c | c}
 0 & b_1, \ldots\ldots, b_{n -2}, b_{n - 1} \\
\hline 
\bm{0} & Q^\flat  \cdot J(0, b_1, \ldots, b_{n - 2})
\end{array} 
\right) 
 = 
\left(
 \begin{array}{c | c}
 1 & \bm{0} \\
 \hline 
 \bm{0} & Q^\flat
 \end{array}
\right) 
\cdot 
\beta_1 . 
\end{eqnarray*}
The $(n + 1)$-th columns of $A(T) \cdot P$ and $P \cdot B(T)$ are $b_n \bm{e}_1$. 
Hence $A(T) \cdot P = P \cdot B(T)$. 
\QED

As a corollary of Theorem 1.41, we have the following (cf. the proof of Corollary 1.38):

\begin{cor}
Let $n \geq 2$ be an integer. 
Let $A(T)$ and $B(T)$ be exponential matrices of $\mf{J}_{[n, 1]}^E$ with the following forms: 
\begin{eqnarray*}
\left\{
\begin{array}{ll}
A(T) = \left(
\begin{array}{c | c}
 \Exp\left( \ds \sum_{i = 1}^{n - 1} f_i \nu_n^i \right) & f_n \bm{e}_1 \\
\hline 
 \bm{0} & 1 
\end{array}
\right)
\quad &  
\left( 
\begin{tabular}{l}
\text{$f_i$ $(1 \leq i \leq n)$ are $p$-polynomials, and} \\
\text{$f_1$ and $f_n$ are linearly independent over $k$} 
\end{tabular} 
\right), \\ [10mm] 
B(T) = \left(
\begin{array}{c | c}
 \Exp\left( \ds \sum_{i = 1}^{n - 1} g_i \nu_n^i \right) & g_n \bm{e}_1 \\
\hline 
 \bm{0} & 1 
\end{array}
\right)
\quad &  
\left( 
\begin{tabular}{l}
\text{$g_i$ $(1 \leq i \leq n)$ are $p$-polynomials, and} \\
\text{$g_1$ and $g_n$ are linearly independent over $k$} 
\end{tabular} 
\right) . 
\end{array}
\right. 
\end{eqnarray*}
Then the following conditions {\rm (1)} and {\rm (2)} are equivalent: 
\begin{enumerate}
\item[\rm (1)] $A(T)$ and $B(T)$ are equivalent. 
\item[\rm (2)] There exists a matrix $Q$ of $\cQ_{[n, 1]}$ such that $(g_1, \ldots, g_n) = (f_1, \ldots, f_n) Q$. 
\end{enumerate}
\end{cor}

\subsubsection{$\mf{J}_{[1, n]}^1$}

As corollaries of Theorem 1.41 and Corollary 1.42, we have the following Corollaries 1.43 and 1.44, respectively (see Lemmas 1.6 and 
1.7).

\begin{cor}
Let $n \geq 2$ be an integer. 
Let $A(T)$ and $B(T)$ be polynomial matrices of $\Mat(n + 1, k[T])$ with the following forms: 
\begin{eqnarray*}
\left\{
\begin{array}{ll}
A(T) = 
\left(
\begin{array}{ c | c }
 0      & a_n {^t}\bm{e}_n \\
\hline 
\bm{0} & \ds \sum_{i = 1}^{n - 1} a_i \nu_n^i  
\end{array}
\right)
\quad & 
\left(
\begin{tabular}{l}
\text{$a_1, \ldots, a_{n - 1}$ are linearly independent over $k$,} \\
\text{and $a_1$, $a_n$ are also linearly independent over $k$}
\end{tabular} 
\right) , 
\\ [10mm] 
B(T) = 
\left(
\begin{array}{ c | c }
 0      & b_n \transpose \bm{e}_n \\
\hline 
\bm{0} & \ds \sum_{i = 1}^{n - 1} b_i \nu_n^i  
\end{array}
\right)
\quad & 
\left(
\begin{tabular}{l}
\text{$b_1, \ldots, b_{n - 1}$ are linearly independent over $k$,} \\
\text{and $b_1$, $b_n$ are also linearly independent over $k$}
\end{tabular} 
\right) . 
\end{array}
\right. 
\end{eqnarray*}
Then the following conditions {\rm (1)} and {\rm (2)} are equivalent: 
\begin{enumerate}
\item[\rm (1)] $A(T)$ and $B(T)$ are equivalent. 
\item[\rm (2)] There exists a matrix $Q$ of $\cQ_{[n, 1]}$ such that $(b_1, \ldots, b_n) = (a_1, \ldots, a_n) Q$. 
\end{enumerate}
\end{cor}

\begin{cor}
Let $n \geq 2$ be an integer. 
Let $A(T)$ and $B(T)$ be exponential matrices of $(\mf{J}_{[1, n]}^1)^E$ with the following forms: 
\begin{eqnarray*}
\left\{
\begin{array}{ll}
A(T) = \left(
\begin{array}{c | c}
1 & f_n \transpose \bm{e}_n \\ 
\hline 
\bm{0} & \Exp\left( \ds \sum_{i = 1}^{n - 1} f_i \nu_n^i \right) 
\end{array}
\right)
\quad &  
\left( 
\begin{tabular}{l}
\text{$f_i$ $(1 \leq i \leq n)$ are $p$-polynomials, and} \\
\text{$f_1$ and $f_n$ are linearly independent over $k$} 
\end{tabular} 
\right), \\ [10mm] 
B(T) = \left(
\begin{array}{c | c}
1 & g_n \transpose \bm{e}_n \\ 
\hline 
\bm{0} & \Exp\left( \ds \sum_{i = 1}^{n - 1} g_i \nu_n^i \right)   
\end{array}
\right)
\quad &  
\left( 
\begin{tabular}{l}
\text{$g_i$ $(1 \leq i \leq n)$ are $p$-polynomials, and} \\
\text{$g_1$ and $g_n$ are linearly independent over $k$} 
\end{tabular} 
\right) . 
\end{array}
\right. 
\end{eqnarray*}
Then the following conditions {\rm (1)} and {\rm (2)} are equivalent: 
\begin{enumerate}
\item[\rm (1)] $A(T)$ and $B(T)$ are equivalent. 
\item[\rm (2)] There exists a matrix $Q$ of $\cQ_{[n, 1]}$ such that $(g_1, \ldots, g_n) = (f_1, \ldots, f_n) Q$. 
\end{enumerate}
\end{cor}

\subsubsection{$\mf{A}(i_1, i_2, i_3)$}

\begin{lem}
Let $A(T), B(T) \in \mf{A}(i_1, i_2, i_3)^\circ$. 
Write $A(T) = \Lambda(i_1, i_2, i_3; \alpha(T))$ and $B(T) = \Lambda(i_1, i_2, i_3; \beta(T))$. 
Then the following conditions {\rm (1)} and {\rm (2)} are equivalent: 
\begin{enumerate}
\item[\rm (1)] $A(T)$ and $B(T)$ are equivalent. 
\item[\rm (2)] We have $\beta(T) = Q \, \alpha(T) \, R$ for some $Q \in GL(i_1, k)$ and $R \in GL(i_3, k)$.  
\end{enumerate}
\end{lem}

\Proof We first prove (1) $\Longrightarrow$ (2). There exists a regular matrix $P \in GL(i_1 + i_2 + i_3, k)$ such that 
$P^{-1} A(T) P = B(T)$. Using the rank conditions of $\alpha(T)$ and $\beta(T)$, we can write $P$ as 
\[
 P = 
\left(
\begin{array}{c | c | c}
 P_1 & P_2 & P_3 \\
\hline 
 O & P_4 & P_5 \\ 
\hline 
 O & O & P_6
\end{array}
\right) \quad (\; P_1 \in GL(i_1, k), \; P_4 \in GL(i_2, k), \; P_6 \in GL(i_3, k) \;) . 
\]
Letting $Q := P_1^{-1}$ and $R := P_6$, we have $\beta(T) = Q \, \alpha(T) \, R$. 

We next prove (2) $\Longrightarrow$ (1). 
Let 
\[
P := 
\left(
\begin{array}{c | c | c}
 Q^{-1} & O & O \\
\hline 
 O & I_{i_2} & O \\
\hline 
 O & O & R \\
\end{array}
\right)
\in GL(i_1 + i_2 + i_3, k) .
\]
Then $P$ satisfies $P^{-1} A(T) P = B(T)$. 
\QED

\section{Polynomial matrices belonging to Heisenberg groups}

Let $R$ be a commutative ring with unity.  
Let $m \geq 1$ be an integer. 
For $x_1, \ldots, x_m, y_1, \ldots, y_m, z \in R$, we define an upper triangular matrix $\eta(x_1, \ldots, x_m, y_1, \ldots, y_m, z) \in \Mat(m + 2, R)$ as 
\begin{eqnarray*}
\eta(x_1, \ldots, x_m, y_1, \ldots, y_m, z):= 
\left( 
\begin{array}{cccccc} 
 1 & x_1 & \cdots & \cdots & x_m  & z \\
   &   1 &         0 & \cdots  &     0 & y_1 \\
   &      & 1        & \ddots  & \vdots & \vdots \\
   &      &   & \ddots          &   0        & \vdots \\
   &    &             &          & 1   & y_m \\
   &     &              &    &        & 1
\end{array} 
\right) . 
\end{eqnarray*}
Let $\mf{H}(m + 2, R)$ be the set of all matrices of $\Mat(m + 2, R)$ 
with the form $\eta(x_1, \ldots, x_m, y_1, \ldots, y_m, z)$. 
Clearly, $\mf{H}(m + 2, k[T])$ is a subgroup of the general linear group $GL(m + 2, R)$. 
We call the subgroup $\mf{H}(m + 2, R)$ the {\it Heisenberg group} over $R$ of {\it degree} $m + 2$.

In this section, we study polynomial matrices of $\mf{H}(m + 2, k[T])$, where $k$ is a field of positive characteristic $p$. 

\subsection{On $\mf{H}(m + 2, k[T])$}

In order to consider equivalent forms of polynomial matrices of $\mf{H}(m + 2, k[T])$, 
we define subsets ${^\ell} \mf{H}^{r_1}_{r_2}$ indexed by $(\ell, r_1, r_2) \in \Omega_m$ 
and subsets $\daleth_{m + 2}^{(i)}$ indexed by $1 \leq i \leq m + 1$ as follows:

Let 
\[
\Omega_m := \{ (\ell, r_1, r_2) \in \Z^3 \mid 
1 \leq \ell\leq m, \;\,  
1 \leq r_1 \leq \ell, \;\, 
0 \leq r_2 \leq m - \ell  
\}.
\] 
For any $(\ell, r_1, r_2) \in \Omega_m$, 
we denote by ${^\ell \mf{H}}^{r_1}_{r_2}$ the set of all elements $\eta(a_1, \ldots, a_{2m + 1})$ 
of $\mf{H}(m + 2, k[T])$ satisfying the following conditions (1), (2), (3): 
\begin{enumerate}
\item[\rm (1)] $a_1, \ldots, a_\ell$ are linearly independent over $k$. 

\item[\rm (2)] $a_{m + 1}, \ldots, a_{m + r_1}, a_{m + \ell + 1}, \ldots, a_{m + \ell + r_2}$ are linearly independent over $k$. 

\item[\rm (3)] $a_i = 0$ for all integers $i$ within one of the following ranges: 
\[
\ell + 1 \leq i \leq m, \qquad m + r_1 + 1 \leq i \leq m + \ell, \qquad \quad m + \ell + r_2 + 1 \leq i \leq 2m . 
\]
\end{enumerate}

We can describe any element $A(T)$ of ${^\ell}\mf{H}^{r_1}_{r_2}$ as 
\begin{eqnarray*}
A(T) = 
\left(
\begin{array}{ c | c | c | c | c | c}
1 & \bm{a}_1 &  \bm{a}_2  & \bm{0} & \bm{0} & a_{2m + 1} \\ 
\hline 
 \bm{0} & I_{r_1} & O & O & O & {^t} \bm{\alpha}_1 \\ 
\hline 
 \bm{0} & O & I_{\ell - r_1} & O & O & \bm{0} \\ 
\hline 
 \bm{0} & O & O & I_{r_2} & O & {^t} \bm{\alpha}_2  \\ 
\hline 
 \bm{0} & O & O & O & I_{m - \ell - r_2} & \bm{0} \\ 
\hline 
 0        & \bm{0}    & \bm{0}    & \bm{0} & \bm{0} & 1
\end{array}
\right) , \qquad  
\left( 
\begin{array}{l}
\bm{a}_1 := (a_1, \ldots, a_{r_1}), \\
\bm{a}_2 := (a_{r_1 + 1}, \ldots, a_\ell), \\
\bm{\alpha}_1 := (a_{m + 1}, \ldots, a_{m + r_1}), \\
\bm{\alpha}_2 := (a_{m + \ell + 1}, \ldots, a_{m + \ell + r_2}) 
\end{array} 
\right) . 
\end{eqnarray*}

For any integer $1 \leq i \leq m + 1$, we denote by $\daleth_{m + 2}^{(i)}$ the set of all elements 
$\eta(a_1, \ldots, a_{2m + 1})$ of $\mf{H}(m + 2, k[T])$ such that the $j$-th entry $a_j$ is zero if $j$ is located outside of the ranges 
$j  = i, \ldots,  m + i - 1$ and $j = 2m + 1$, 
i.e., if we see the sequence $a_1, \ldots, a_{2m + 1}$ in the matrix $\eta(a_1, \ldots, a_{2m + 1})$, then  
\[
(a_1, \ldots, a_m \mid  a_{2m + 1} \mid a_{m + 1}, \ldots, a_{2m})
 = (0, \ldots, 0, \overset{m}{\overbrace{ a_i, \ldots, a_m \mid a_{2m + 1} \mid a_{m + 1}, \ldots, a_{m + i -1}}}, 0, \ldots, 0). 
\] 
Let $\daleth_{m + 2} := \bigcup_{i = 1}^{m + 1} \daleth_{m + 2}^{(i)}$.

\begin{thm} 
The following assertions {\rm (1)} and {\rm (2)} hold true: 
\begin{enumerate}
\item[\rm (1)] $\ds \mf{H}(m + 2, k[T]) \leadsto \daleth_{m +2} \cup  \bigcup_{(\ell, r_1, r_2) \in \Omega_m} {^\ell \mf{H}}^{r_1}_{r_2}$. 
\item[\rm (2)] Two sets $\daleth_{m + 2}$ and $\ds \bigcup_{(\ell, r_1, r_2) \in \Omega_m} {^\ell \mf{H}}^{r_1}_{r_2}$ 
are mutually $GL(m + 2, k)$-disjoint. 
\end{enumerate} 
\end{thm}

In order to prove Theorem 2.1, we prepare lemmas in Subsubsection 2.1.1. 
Subsequently, Subsubsection 2.1.2, we give a proof of Theorem 2.1.

\subsubsection{Lemmas} 

For all integers $0 \leq \ell, r \leq m$, we denote by $\mf{H}_{\ell, r}$ the set of all polynomial matrices $A(T)$ of 
$\mf{H}(m + 2, k[T])$ satisfying $\ell_{A(T)} = \ell$ and $r_{A(T)} = r$. 
So, $\mf{H}(m + 2, k[T]) = \bigcup_{0 \leq \ell, r \leq m} \mf{H}_{\ell, r}$.

\begin{lem} 
We have 
\[
\mf{H}(m + 2, k[T]) = \daleth_{m + 2} \cup \bigcup_{1 \leq \ell, r \leq m} \mf{H}_{\ell, r}. 
\] 
\end{lem} 

\Proof The proof follows from $(\bigcup_{0 \leq r \leq m} \mf{H}_{0, r}) \cup ( \bigcup_{0 \leq \ell \leq m} \mf{H}_{\ell, 0}) \subset \daleth_{m + 2}$. 
\QED

For integers $1 \leq \ell, r \leq m$, we define a subset ${^\ell \mf{H}_r}$ of $\mf{H}_{\ell, r}$ as 
the set consisting of all elements $\eta(a_1, \ldots, a_{2m + 1}) \in \mf{H}_{\ell, r}$ satisfying $a_i = 0$ for all $\ell + 1 \leq i \leq m$.

For $A(T) \in {^\ell} \mf{H}_r$, 
we define subvector spaces $V_+(A(T))$ and $V_-(A(T))$ of $k[T]$ as 
\begin{eqnarray*}
\left\{ 
\begin{array}{rcl}
 V_+ (A(T)) &:=& \Span \{ a_i \mid m + 1 \leq i \leq m + \ell \}, \\
 V_- (A(T)) &:=& \Span \{ a_i \mid m + \ell + 1 \leq i \leq 2m \}. 
\end{array}
\right. 
\end{eqnarray*} 
Clearly, $r \leq V_+(A(T)) + V_-(A(T)) \leq m$. 
For integers $r_1, r_2$ satisfying $0 \leq r_1 \leq \ell$ and $0 \leq r_2 \leq m - \ell$, 
we define a subset ${^\ell \mf{H}_{r, r_1, r_2}}$ of ${^\ell} \mf{H}_r$ as the 
set consisting of all elements $A(T) \in {^\ell \mf{H}}_r$ satisfying 
$\dim_k V_+(A(T)) = r_1$ and $\dim_k V_-(A(T)) = r_2$. 
So, ${^\ell} \mf{H}_r = \bigcup_{r \leq r_1 + r_2 \leq m} {^\ell} \mf{H}_{r, r_1, r_2}$.

\begin{lem} 
For all integers $1 \leq \ell, r \leq m$, we have 
\[
 \mf{H}_{\ell, r} \leadsto  \daleth_{m + 2}
 \cup 
 \bigcup_{\substack{r \leq r_1 + r_2 \leq m \\ 1 \leq r_1 \leq \ell}}{^\ell \mf{H}_{r, r_1, r_2}} , 
\] 
\end{lem}

\Proof We have $\mf{H}_{\ell, r} \leadsto {^\ell \mf{H}_r}
 = {^\ell \mf{H}}_{r, 0, r} \cup \bigcup_{r \leq r_1 + r_2 \leq m , \; 1 \leq r_1 \leq \ell} {^\ell \mf{H}_{r, r_1, r_2}}$.  
Since $ {^\ell \mf{H}}_{0, r} \leadsto \daleth_{m + 2}$, 
we obtain the desired formula. 
\QED

\begin{lem} 
For all integers $r_1, r_2$ satisfying $1 \leq r_1 \leq \ell$, $0 \leq r_2 \leq m - \ell$ and $r \leq r_1 + r_2 \leq m$,  we have 
\[
{^\ell \mf{H}_{r, r_1, r_2}} \leadsto \bigcup_{(\ell, r_1, r_2) \in \Omega_m}  {^\ell \mf{H}^{r_1}_{r_2}}. 
\]
\end{lem}

\Proof For any $A(T) \in {^\ell} \mf{H}_{r, r_1, r_2}$, 
it is possible for some lower triangular matrix $P \in GL(m + 2, k)$  
that $P^{-1} A(T) P \in {^\ell} \mf{H}^{s_1}_{s_2}$, where $s_1, s_2$ satisfy $1 \leq s_1 \leq r_1$ and $0 \leq s_2 \leq r_2$. 
\QED

\subsubsection{A proof of Theorem 2.1}

We shall prove assertion (1). Using Lemmas 2.2, 2.3 and 2.4, we have 
\[
 \mf{H}(m + 2, k[T])
 \leadsto \daleth_{m + 2} \cup \bigcup_{1 \leq \ell, r \leq m} \bigcup_{\substack{r \leq r_1 + r_2 \leq m \\ 1 \leq r_1 \leq \ell}} {^\ell \mf{H}_{r, r_1, r_2}} 
 \leadsto \daleth_{m + 2} \cup \bigcup_{(\ell, r_1, r_2) \in \Omega_m} {^\ell \mf{H} }^{r_1}_{r_2} . 
\]

We shall prove assertion (2). 
Supposing $A(T)$ and $B(T)$ are equivalent, we shall find a contradiction. 
There exists a regular matrix $P$ of $GL(m + 2, k[T])$ such that $P^{-1} A(T) P = B(T)$. 
We have $P^{-1} ( A(T) - I ) P = B(T) - I$. 
Since $A(T) \in \daleth_{m + 2}$, we know $( P^{-1} ( A(T) - I ) P ) \cdot ( P^{-1} ( A(T') - I ) P ) = O$, 
where $T, T'$ are indeterminates over $k$. So, $(B(T) - I) \cdot (B(T') - I) = O$. 
We can express $B(T)$ as $B(T) = \eta(b_1, \ldots, b_{2m + 1})$, where $b_1, \ldots, b_{2m + 1} \in k[T]$. 
So, we have $\sum_{i = 1}^{r_1} b_i(T) \cdot b_{m + i}(T') = 0$. 
This equality implies that $b_1(T), \ldots, b_{r_1}(T)$ are linearly dependent over $k$, 
since there exists a non-zero polynomial among $b_{m + 1}(T'), \ldots, b_{m + r_1}(T')$ (see the definition of ${^\ell \mf{H}}^{r_1}_{r_2}$). 
However, $b_1(T), \ldots, b_{r_1}(T)$ are linearly independent over $k$, since $B(T) \in {^\ell \mf{H}}^{r_1}_{r_2}$. 
We have a contradiction.

\subsection{On $\daleth_{m + 2}$}

For considering equivalent forms of polynomial matrices of $\daleth_{m +2}$, we define several notations as follows: 

For $i \geq 1$ and $j \geq 1$, we denote by ${^i}\mf{L}_j$ the set of all matrices $L(a_1, \ldots, a_i \mid a_{i + 1}, \ldots,  a_{i + j} \mid  a)$ of 
$\Mat_{j + 1, i + 1}(k[T])$ with the following form 
\begin{eqnarray*}
\begin{array}{c}
L(a_1, \ldots, a_i \mid a_{i + 1}, \ldots, a_{i + j} \mid a) = 
\left(
\begin{array}{ccc | c}
 a_1 & \cdots & a_i & a \\
\hline 
      &           &     & a_{i + 1} \\
      &     O    &     & \vdots \\
      &           &     & a_{i + j}
\end{array}
\right) , 
 \\ [13mm] 
\left(
\begin{array}{l}
\text{$a_1, \ldots, a_i$ are linearly independent over $k$, and } \\
\text{$a_{i + 1}, \ldots, a_{i + j}$ are also linearly independent over $k$}
\end{array}
\right) . 
\end{array} 
\end{eqnarray*}
We denote by $({^1}\mf{L}_{1})^{\dim = d}$ the set of all matrices $L(a_1 \mid a_2 \mid a)$ of 
${^1} \mf{L}_{1}$ satisfying $\dim_k \Span \{ a_1, a_2, a \}  = d$. 
So, the value of $d$ is two or three.

For $i \geq 1$ and $j \geq 1$ satisfying $i + j \leq m$, we denote by ${^i}\daleth_j$ 
the set of all matrices $A(T)$ of $\daleth_{m  +2}$ with the form $A(T) = \Lambda(j + 1, i + 1; L)$ 
for some $L \in {^i}\mf{L}_j$. 
We denote by $({^1}\daleth_1)^{\dim = d}$ the set of all matrices $\Lambda(2, 2; L)$, 
where $L \in ({^1}\mf{L}_{1})^{\dim = d}$.

Let ${^1} \mf{L}_1^0$ denote the set of all matrices $L \in {^1}\mf{L}_1$ 
with the form $L = L(a_1 \mid a_2 \mid 0)$, 
and let ${^1}\daleth_1^0$ denote the set of all matrices $\Lambda(2, 2; L)$, where $L \in {^1}\mf{L}_1^0$. So, ${^1}\daleth_1^0 \subset ({^1}\daleth_1)^{\dim = 2}$. 

Now, let 
\begin{eqnarray*}
\left\{
\begin{array}{l}
\daleth_{m + 2}^{\rm I} := \{ I_{m + 2} \} ,  \\ [3mm] 
\ds \daleth_{m + 2}^{\rm II} := \bigcup_{1 \leq j \leq m + 1} \mf{A}(1, j)^\circ, \qquad 
\ds \daleth_{m + 2}^{\rm III} := \bigcup_{2 \leq i \leq m + 1} \mf{A}(i, 1)^\circ, \\  [7mm] 
\ds \daleth_{m + 2}^{\rm IV} := {^1}\daleth_1^0 , \qquad 
\ds \daleth_{m + 2}^{\rm V} := ({^1}\daleth_1)^{\dim = 3}, \qquad 
\ds \daleth_{m + 2}^{\rm VI} := \bigcup_{3 \leq i + j \leq m} {^i}\daleth_j .
\end{array}
\right.
\end{eqnarray*}
are mutually $GL(m + 2, k)$-disjoint union.  
We know $\daleth_{m + 2}^{\rm II} \subset \mf{A}(1, m + 1)$ and $\daleth_{m + 2}^{\rm III} \subset \mf{A}(m + 1, 1)$, 
and 
we index $\daleth_{m + 2}^{\rm III}$ from $i  =2$ so that $\daleth_{m + 2}^{\rm II} \cap \daleth_{m + 2}^{\rm III} = \emptyset$.  

The following theorem tells us that any polynomial matrix of $\daleth_{m + 2}$ 
has an equivalent form belonging to one of the six subsets $\daleth_{m + 2}^{\rm I}, \daleth_{m + 2}^{\rm II}, \daleth_{m + 2}^{\rm III}, 
\daleth_{m + 2}^{\rm IV}, \daleth_{m + 2}^{\rm V}, \daleth_{m + 2}^{\rm VI}$, and 
this belonging is unique.

\begin{thm}
The following assertions {\rm (1)} and {\rm (2)} hold true: 
\begin{enumerate}
\item[\rm (1)] 
$\daleth_{m + 2}
 \leadsto \daleth_{m + 2}^{\rm I} \cup \daleth_{m + 2}^{\rm II} \cup \daleth_{m + 2}^{\rm III}
 \cup \daleth_{m + 2}^{\rm IV} \cup \daleth_{m + 2}^{\rm V} \cup \daleth_{m + 2}^{\rm VI}$. 
\item[\rm (2)] The six subsets $\daleth_{m + 2}^{\rm I}, \daleth_{m + 2}^{\rm II}, \daleth_{m + 2}^{\rm III}, 
\daleth_{m + 2}^{\rm IV}, \daleth_{m + 2}^{\rm V}, \daleth_{m + 2}^{\rm VI}$ are mutually $GL(m + 2, k)$-disjoint. 
\end{enumerate}
\end{thm}

We shall give a proof of the above theorem after considering $({^1}\daleth_1)^{\dim = 2}$ 
(see the following Subsubsections 2.2.1 and 2.2.2).

\subsubsection{$({^1}\daleth_1)^{\dim = 2}$}

We have the following lemma concerning equivalent forms of elements of $({^1}\daleth_1)^{\dim = 2}$. 

\begin{lem}
We have $({^1}\daleth_1)^{\dim = 2} \leadsto {^1}\daleth_1^0$. 
\end{lem} 

\Proof The proof follows from Lemma 1.45 and the following Lemma 2.7. 
\QED

\begin{lem}
For any $L = L(a_1 \mid a_2 \mid a) \in {^1}\mf{L}_1$, 
the following conditions {\rm (1)} and {\rm (2)} are equivalent: 
\begin{enumerate}
\item[\rm (1)] $L \in ({^1}\mf{L}_1)^{\dim = 2}$. 
\item[\rm (2)] $Q \cdot L \cdot R  =  L(a_1 \mid a_2 \mid 0)$ for some $Q, R \in GL(2, k)$. 
\end{enumerate} 
\end{lem}

\Proof We first prove (1) $\Longrightarrow$ (2). Write $a = \lambda a_1 + \mu a_2$ $(\lambda, \mu \in k)$. 
Letting 
\[
Q = 
\left( 
\begin{array}{cc}
 1 & - \lambda \\
 0 & 1
\end{array}
\right) \qquad \text{and} \qquad
R = 
\left( 
\begin{array}{cc}
 1 & - \mu \\
 0 & 1
\end{array}
\right), 
\]
we have the desired equality. 
We next prove (2) $\Longrightarrow$ (1). Since $L =  Q^{-1} \cdot L(a_1 \mid a_2 \mid 0) \cdot R^{-1}$, 
we know that $a$ is a linear combination of $a_1$ and $a_2$ over $k$. 
\QED

\subsubsection{A proof of Theorem 2.5}

We first prove assertion (1). 
We have $\daleth_{m +2} \leadsto 
\daleth_{m + 2}^{\rm I} \cup \daleth_{m + 2}^{\rm II} \cup \daleth_{m + 2}^{\rm III}
 \cup  ({^1}\daleth_1)^{\dim = 2} \cup \daleth_{m + 2}^{\rm V} \cup \daleth_{m + 2}^{\rm VI}$ (cf. Lemma 1.29), 
and then by Lemma 2.7, we have the desired result.

We next prove assertion (2). 
We know from Lemma 1.32 that five sets $\daleth_{m + 2}^{\rm I}$, $\daleth_{m + 2}^{\rm II}$, $\daleth_{m + 2}^{\rm III}$, 
$\daleth_{m + 2}^{\rm IV} \cup \daleth_{m + 2}^{\rm V}$, $\daleth_{m + 2}^{\rm VI}$ are mutually $GL(m + 2, k)$-disjoint. 
And we can obatain from Lemmas 1.45 and 2.6 that two sets $\daleth_{m + 2}^{\rm IV}$ and $\daleth_{m + 2}^{\rm V}$ 
are mutually $GL(m + 2, k)$-disjoint. 
\QED

\subsection{${V_{m + 2}}^{A(T)}$ and $(V_{m + 2}^*)^{A(T)}$ for $A(T) \in \mf{H}(m + 2, k[T])$}

\begin{lem}
Let $A(T)$ be a polynomial matrix of $\mf{H}(m + 2, k[T])$ 
and regard $A(T)$ as $k$-linear maps from $V_{m + 2}$ to $k[T] \otimes_k V_{m + 2}$ 
and from $V_{m + 2}^*$ to $k[T] \otimes_k V^*$. 
For simplicity, let $V := V_{m + 2}$, $\ell := \ell_{A(T)}$, $r := r_{A(T)}$, $t := t_{A(T)}$. 
Then we have 
\begin{eqnarray*}
\lefteqn{ (\dim_k V^{A(T)}, \; \dim_k (V^*)^{A(T)} ) } \\
& = &
\left\{
\begin{array}{lll }

 (m - t + 2,  & m - t + 2) \qquad & \text{ if \quad $\ell = 0$ \quad and \quad $r = 0$}, \\

 (m + 1,  & m - t + 2) \qquad & \text{ if \quad $\ell = 0$ \quad and \quad $r \geq 1$}, \\

 (m - t + 2, &  m  + 1) \qquad & \text{ if \quad $\ell \geq 1$ \quad and \quad $r = 0$}, \\

 (m - \ell + 1, & m - r + 1 ) \qquad & \text{ if \quad $\ell \geq 1$ \quad and \quad $r \geq 1$} . 
\end{array}
\right. 
\end{eqnarray*}
\end{lem}

\Proof The proof is straightforward. 
\QED

\begin{lem}
Let $A(T)$ be a polynomial matrix of $\mf{H}(m + 2, k[T])$. 
Then the following assertions {\rm (1), (2), (3)} hold true: 
\begin{enumerate} 
\item[\rm (1)] The following conditions {\rm (1.1)} and {\rm (1.2)} are equivalent: 
\begin{enumerate} 
\item[\rm (1.1)] $\ell_{A(T)} + r_{A(T)} + t_{A(T)} = 0$. 
\item[\rm (1.2)] $(\dim_k V^{A(T)}, \; \dim_k (V^*)^{A(T)}) = (m + 2,\; m + 2)$. 
\end{enumerate}

\item[\rm (2)] The following conditions {\rm (2.1)} and {\rm (2.2)} are equivalent: 
\begin{enumerate} 
\item[\rm (2.1)] $1 \leq \ell_{A(T)} + r_{A(T)} + t_{A(T)} \leq 2$. 
\item[\rm (2.2)] $(\dim_k V^{A(T)}, \; \dim_k (V^*)^{A(T)}) = (m + 1,\; m + 1)$. 
\end{enumerate}

\item[\rm (3)] The following conditions {\rm (3.1)} and {\rm (3.2)} are equivalent: 
\begin{enumerate} 
\item[\rm (3.1)] $\ell_{A(T)} + r_{A(T)} + t_{A(T)} \geq 3$. 
\item[\rm (3.2)] $\max\{ \dim_k V^{A(T)}, \; \dim_k (V^*)^{A(T)} \} \leq m  + 1$ and 
$\min\{ \dim_k V^{A(T)}, \; \dim_k (V^*)^{A(T)} \} \leq m$, 
where for all integers $\alpha, \beta$, we denote by $\max\{\alpha, \beta\}$ the maximum of $\alpha$ and $\beta$, 
and by $\min\{\alpha, \beta\}$ the minimum of $\alpha$ and $\beta$. 
\end{enumerate} 
\end{enumerate} 
\end{lem}

\Proof We have $\ell_{A(T)} + r_{A(T)} + t_{A(T)} = 0$ if and only if $(\ell_{A(T)}, r_{A(T)}, t_{A(T)}) = (0, 0, 0)$.  
By Lemma 2.8, we have the following array: 
\begin{eqnarray*}
\begin{array}{| @{\hspace{5mm}} c @{\hspace{5mm}} c @{\hspace{5mm}} c @{\hspace{5mm}} | @{\hspace{5mm}} c  @{\hspace{5mm}} | @{\hspace{5mm}} c @{\hspace{5mm}} |}
\hline 
 \ell_{A(T)} & r_{A(T)} & t_{A(T)}  & \dim_k V^{A(T)} & \dim_k (V^*)^{A(T)} \\
\hline 
 0     & 0      & 0    & m + 2   & m + 2 \\
\hline 
\end{array}
\end{eqnarray*} 
We have $1 \leq \ell_{A(T)} + r_{A(T)} + t_{A(T)} \leq 2$ if and only if $(\ell_{A(T)}, r_{A(T)}, t_{A(T)}) = (0, 0, 1), (0, 1, 1), (1, 0, 1)$. 
By Lemma 2.8, we have the following array: 
\begin{eqnarray*}
\begin{array}{| @{\hspace{5mm}} c @{\hspace{5mm}} c @{\hspace{5mm}} c @{\hspace{5mm}} | @{\hspace{5mm}} c  @{\hspace{5mm}} | @{\hspace{5mm}} c @{\hspace{5mm}} |}
\hline 
 \ell_{A(T)} & r_{A(T)} & t_{A(T)}  & \dim_k V^{A(T)} & \dim_k (V^*)^{A(T)} \\
\hline 
 0     & 0      & 1    & m + 1   & m + 1 \\
 0     & 1      & 1    & m + 1   & m + 1 \\
 1     & 0      & 1    & m + 1  & m + 1 \\
\hline 
\end{array}
\end{eqnarray*} 
We have $\ell_{A(T)} + r_{A(T)} + t_{A(T)} \geq 3$ if and only if $(\ell_{A(T)}, r_{A(T)}, t_{A(T)})$ coincides with one of the following tuples: 
\begin{eqnarray*}
\begin{array}{l}
(0, 1, 2), \qquad (1, 0, 2),  \qquad (0, r, t) \quad (r \geq 2, \; t \geq 2) , \qquad 
 (\ell, 0,  t) \quad (\ell \geq 2, \; t \geq 2) , \qquad \\
 (\ell, r, t) \quad (\ell \geq 1, \; r \geq 1, \; t \geq 1). 
\end{array} 
\end{eqnarray*}
By Lemma 2.8, we have the following array: 
\begin{eqnarray*}
\begin{array}{| @{\hspace{5mm}} c @{\hspace{5mm}} c @{\hspace{5mm}} c @{\hspace{5mm}} | @{\hspace{5mm}} c  @{\hspace{5mm}} | @{\hspace{5mm}} c @{\hspace{5mm}} |}
\hline 
 \ell_{A(T)} & r_{A(T)} & t_{A(T)}  & \dim_k V^{A(T)} & \dim_k (V^*)^{A(T)} \\
\hline 
 0     & 1      & 2    & m  +1   & m \\
 1     & 0      & 2    & m        & m + 1 \\
 0 & \geq 2 & \geq 2 & m + 1 & \leq m \\
 \geq 2 & 0 & \geq 2 & \leq m & m + 1 \\
\geq 1 & \geq 1 & \geq 1 & \leq m & \leq m \\
\hline  
\end{array}
\end{eqnarray*}
Thus the assertions (1), (2), (3) are proved. 
\QED

\subsection{Equivalence relations of polynomial matrices of Heisenberg groups}

\subsubsection{${^\ell \mf{H}}^{r_1}_{r_2}$}

Given two polynomial matrices of $\bigcup_{(\ell, r_1, r_2) \in \Omega_m} {^\ell \mf{H}}^{r_1}_{r_2}$ which are equivalent, 
we can prove that the two polynomial matrices are in the same ${^\ell \mf{H}}^{r_1}_{r_2}$. 
This is ensured by the following lemma:  

\begin{lem} 
Let $A(T) \in {^\ell \mf{H}}^{r_1}_{r_2}$ and $B(T) \in {^l \mf{H}}^{s_1}_{s_2}$. 
Assume that $A(T)$ and $B(T)$ are equivalent. 
Then we have $\ell = l$, $r_1 = s_1$, and $r_2 = s_2$. 
\end{lem}

\Proof We know from Lemma 2.8 that $\ell = l$ and $r_1 + r_2 = s_1 + s_2$. 
So, we have only to show that $r_1 = s_1$. 
We can write $A(T)$ and $B(T)$ as  $A(T) := \eta(a_1, \ldots, a_{2m + 1})$ and $B(T) = \eta(b_1, \ldots, b_{2m + 1})$ 
for some $a_1, \ldots, a_{2m + 1}, b_1, \ldots, b_{2m + 1} \in k[T]$. 
Let 
\begin{eqnarray*}
\left\{ 
\begin{array}{ccl}
\bm{a} &=& (a_1, \ldots, a_\ell),                                     \\
\bm{f}_1 &=& (a_{m + 1}, \ldots, a_{m + \ell}), \\
\bm{f}_2 &=& (a_{m + \ell + 1}, \ldots, a_{2m}), \\ 
\alpha &=& a_{2m + 1}, 
\end{array}
\right. 
\qquad \text{ and } \qquad 
\left\{
\begin{array}{ccl}
\bm{b} &=& (b_1, \ldots, b_\ell),   \\
\bm{g}_1 &=& (b_{m + 1}, \ldots, b_{m + \ell}),   \\  
\bm{g}_2 &=& (b_{m + \ell + 1}, \ldots, b_{2m}), \\  
\beta &=& b_{2m + 1} . 
\end{array}
\right. 
\end{eqnarray*}
There exists a regular matrix $P \in GL(m + 2, k)$ such that $P^{-1} A(T) P = B(T)$. 
We can write $P$ as $P = (P_{i, j})_{1 \leq i, j \leq 4}$, where 
the sizes of the rows of these blocks are $1$, $\ell$, $m - \ell$, and $1$ from the top to the bottom, 
and the sizes of the columns of these blocks are $1$, $\ell$, $m - \ell$, and $1$ from the left to the right. 
Since $(A(T) - I_{m  + 2} ) P = P (B(T) - I_{m + 2})$, we have 
\[
\left(
\begin{array}{c | c | c | c}
 0       & \bm{a} & \bm{0}       & f \\
\hline 
 \bm{0} & O    & O             & {^t} \bm{f}_1 \\
\hline 
 \bm{0} & O       & O & {^t} \bm{f}_2 \\
\hline 
 0       & \bm{0}  & \bm{0}       & 0  
\end{array}
\right) 
\left(
\begin{array}{c | c | c |c}
 P_{1, 1} & P_{1, 2} & P_{1, 3} & P_{1, 4} \\ 
\hline 
 P_{2, 1} & P_{2, 2} & P_{2, 3} & P_{2, 4} \\ 
\hline 
 P_{3, 1} & P_{3, 2} & P_{3, 3} & P_{3, 4} \\ 
\hline 
 P_{4, 1} & P_{4, 2} & P_{4, 3} & P_{4, 4}
\end{array}
\right)
 = 
\left(
\begin{array}{c | c | c |c}
 P_{1, 1} & P_{1, 2} & P_{1, 3} & P_{1, 4} \\ 
\hline 
 P_{2, 1} & P_{2, 2} & P_{2, 3} & P_{2, 4} \\ 
\hline 
 P_{3, 1} & P_{3, 2} & P_{3, 3} & P_{3, 4} \\ 
\hline 
 P_{4, 1} & P_{4, 2} & P_{4, 3} & P_{4, 4}
\end{array}
\right)
\left(
\begin{array}{c | c | c | c}
 0        & \bm{b} & \bm{0}       & g \\
\hline 
 \bm{0} & O    & O             & {^t} \bm{g}_1 \\
\hline 
 \bm{0} & O       & O  & {^t} \bm{g}_2 \\
\hline 
 0       & \bm{0}  & \bm{0}       & 0   
\end{array}
\right) . 
\]
So, 
\begin{eqnarray*}
\lefteqn{
\left(
\begin{array}{ c | c | c | c}
 \bm{a} \cdot P_{2, 1} + f \cdot P_{4, 1} & \bm{a} \cdot P_{2, 2} + f \cdot P_{4, 2} & \bm{a} \cdot P_{2, 3} + f \cdot P_{4, 3} & \bm{a} \cdot P_{2, 4} + f \cdot P_{4, 4} \\ 
\hline 
 {^t} \bm{f}_1 \cdot P_{4, 1}              & {^t} \bm{f}_1 \cdot P_{4, 2}            & {^t} \bm{f}_1 \cdot P_{4, 3}            & {^t} \bm{f}_1 \cdot P_{4, 4} \\
\hline 
{^t} \bm{f}_2 \cdot P_{4, 1}              &  {^t} \bm{f}_2 \cdot P_{4, 2}            & {^t} \bm{f}_2 \cdot P_{4, 3}            & {^t} \bm{f}_2 \cdot P_{4, 4} \\
\hline  
 0                                        & \bm{0}                                 & \bm{0}                                & 0                         
\end{array}
\right)
}\\
& = &   
\left(
\begin{array}{ @{\hspace{8mm}} c @{\hspace{8mm}}  | @{\hspace{7mm}} r @{\hspace{7mm}}  | @{\hspace{8mm}}  c @{\hspace{8mm}}  |  c}
 0 & P_{1, 1} \cdot \bm{b} & O & P_{1, 1} \cdot g + P_{1, 2} \cdot {^t} \bm{g}_1 + P_{1, 3} \cdot {^t} \bm{g}_2 \\
\hline 
 \bm{0} & P_{2, 1} \cdot \bm{b} & O & P_{2, 1} \cdot g + P_{2, 2} \cdot {^t} \bm{g}_1 + P_{2, 3} \cdot {^t} \bm{g}_2 \\
\hline 
 \bm{0} & P_{3, 1} \cdot \bm{b} & O & P_{3, 1} \cdot g + P_{3, 2} \cdot {^t} \bm{g}_1 + P_{3, 3} \cdot {^t} \bm{g}_2 \\
\hline 
 0 & P_{4, 1} \cdot \bm{b} & O & P_{4, 1} \cdot g + P_{4, 2} \cdot {^t} \bm{g}_1 + P_{4, 3} \cdot {^t} \bm{g}_2 \\
\end{array}
\right) . 
\end{eqnarray*}
The $(2, 3)$-th entries imply $P_{4, 3} = O$. 
The $(1, 3)$-th entries imply $P_{2, 3} = O$. 
The $(4, 2)$-th entries imply $P_{4, 1} = O$. 
The $(1, 1)$-th entries imply $P_{2, 1} = O$. 
The $(2, 2)$-th entries imply $P_{4, 2} = O$. 
The $(3, 2)$-th entries imply $P_{3, 1} = O$. 
Since $P$ is regular, $P_{4, 4} \ne ( 0 )$. 
The $(2, 4)$-th entries imply ${^t} \bm{f}_1 \cdot P_{4, 4} = P_{2, 2} \cdot {^t}\bm{g}_1$. 
So, we have 
\[
r_1 = \dim \Span_k \{ a_{m + i} \mid 1 \leq i \leq r_1 \} \leq \dim \Span_k \{ b_{m + i} \mid 1 \leq i \leq s_1 \} = s_1. 
\] 
Similarly, using $(B(T) - I_{m + 2}) P^{-1} = P^{-1} (A(T) - I_{m + 2})$, we have $s_1 \leq r_1$. Hence, we obtain $r_1 = s_1$. 
\QED

Now, we are interested in equivalence relations of two polynomial matrices of ${^\ell}\mf{H}^{r_1}_{r_2}$. 
For this study, we introduce a notation ${^\ell}\V^{r_1}_{r_2}$ which parametrizes 
${^\ell}\mf{H}^{r_1}_{r_2}$. 

Let $R'$ be a not necessarily commutative $k$-algebra. 
A sequence $(f_1, \ldots, f_i)$ of $R'$ is said to be a {\it frame} of $R'$ 
if $f_1, \ldots, f_i$ are linearly independent over $k$. 
For any integer $i \geq 1$, we denote by $\F(i, R')$ the set of all frames $(f_1, \ldots, f_i)$ of $R'$, i.e., 
\[
\F(i, R') := \{ (f_1, \ldots, f_i) \in R'^i \mid \text{$f_1, \ldots, f_i$ are linearly independent over $k$} \} .
\]
For any $(\ell, r_1, r_2) \in \Omega_m$, we define sets ${^\ell}\V^{r_1}_{r_2}$ and ${^\ell}\W^{r_1}_{r_2}$ as 
\begin{eqnarray*}
\left\{ 
\begin{array}{rcl}
 {^\ell}\V^{r_1}_{r_2}  & := & \F(\ell, k[T]) \times \F(r_1 + r_2, k[T]) \times k[T] , \\
 {^\ell}\W^{r_1}_{r_2}  & := &  k[T]^{r_1} \times k[T]^{\ell - r_1} \times k[T]^{r_1} \times k[T]^{r_2} \times k[T] . 
\end{array}
\right. 
\end{eqnarray*}
Clearly, ${^\ell}\V^{r_1}_{r_2} \subset {^\ell}\W^{r_1}_{r_2}$. 
We can define a map $h : {^\ell}\W^{r_1}_{r_2} \to \mf{H}(m + 2, k[T])$ as 
\begin{eqnarray*}
h(\bm{a}_1, \bm{a}_2, \bm{\alpha}_1, \bm{\alpha}_2, \varphi) := 
\left(
\begin{array}{ c | c | c | c | c | c}
1 & \bm{a}_1 &  \bm{a}_2  & \bm{0} & \bm{0} &\varphi  \\ 
\hline 
 \bm{0} & I_{r_1} & O & O & O & {^t} \bm{\alpha}_1 \\ 
\hline 
 \bm{0} & O & I_{\ell - r_1} & O & O & \bm{0} \\ 
\hline 
 \bm{0} & O & O & I_{r_2} & O & {^t} \bm{\alpha}_2  \\ 
\hline 
 \bm{0} & O & O & O & I_{m - \ell - r_2} & \bm{0} \\ 
\hline 
 0        & \bm{0}    & \bm{0}    & \bm{0} & \bm{0} & 1
\end{array}
\right) . 
\end{eqnarray*}
For any element $(\bm{a}_1, \bm{a}_2, \bm{\alpha}_1, \bm{\alpha}_2, \varphi ) \in {^\ell}\W^{r_1}_{r_2}$, 
we observe that 
$h(\bm{a}_1, \bm{a}_2, \bm{\alpha}_1, \bm{\alpha}_2, \varphi ) \in {^\ell}\mf{H}^{r_1}_{r_2}$ if and only if 
$(\bm{a}_1, \bm{a}_2, \bm{\alpha}_1, \bm{\alpha}_2, \varphi ) \in {^\ell}\V^{r_1}_{r_2}$ 
(see the definition of ${^\ell}\mf{H}^{r_1}_{r_2}$).  
Thus, $h$ defines a bijection from ${^\ell}\V^{r_1}_{r_2}$ to ${^\ell}\mf{H}^{r_1}_{r_2}$. 
We have the following commutative diagram, where upward arrows are inclusions: 
\begin{eqnarray*}
\xymatrix@!C=36pt{                               & \mf{H}(m + 2, k[T]) \\
  {^\ell}\W^{r_1}_{r_2} \ar[r]^(.4)h \ar[ru]^{h}               &  h( {^\ell}\W^{r_1}_{r_2}) \ar[u] \\
 {^\ell}\V^{r_1}_{r_2} \ar[r]^{\cong} \ar[u]     &  {^\ell}\mf{H}^{r_1}_{r_2} \ar[u] \\ 
}
\end{eqnarray*}

We consider an equivalence relation of two polynomial matrices $A(T)$, $B(T)$ of $h({^\ell}\W^{r_1}_{r_2})$ 
in terms of two corresponding elements $(\bm{a}_1, \bm{a}_2, \bm{\alpha}_1, \bm{\alpha}_2, \varphi), 
(\bm{b}_1, \bm{b}_2, \bm{\beta}_1, \bm{\beta}_2, \psi)$ of ${^\ell}\W^{r_1}_{r_2}$ (see 
the following Lemma 2.11), 
and then consider equivalence relations of two polynomial matrices of ${^\ell}\mf{H}^{r_1}_{r_2}$ 
in terms of two corresponding elements of ${^\ell}\V^{r_1}_{r_2}$ (see the following Lemma 2.12). 
For the considerations, we define a subgroup ${^\ell}\cP^{r_1}_{r_2}$ of $GL(m + 2, k)$ and 
a subgroup ${^\ell}\cQ^{r_1}_{r_2}$ of $GL(\ell + r_1 + r_2 + 1, k)$  as follows: 

Let $(\ell, r_1, r_2) \in \Omega_m$. 
We denote by ${^\ell \cP}^{r_1}_{r_2}$ the set of all regular matrices $P$ 
of $GL(m + 2, k)$ with the form
\[
 P = 
\left( 
\begin{array}{ c | c | c | c | c | c }
 P_{1, 1} & P_{1, 2} & P_{1, 3} & P_{1, 4} & P_{1, 5} & P_{1, 6} \\
\hline 
\bm{0} & P_{2, 2} & P_{2, 3} & O & O & P_{2, 6} \\
\hline 
\bm{0} & O & P_{3, 3} & O & O & P_{3, 6} \\ 
\hline 
\bm{0} & P_{4, 2} & P_{4, 3} & P_{4, 4} & P_{4, 5} & P_{4, 6} \\ 
\hline 
\bm{0} & O & P_{5, 3} & O & P_{5, 5} & P_{5, 6} \\ 
\hline 
\bm{0} & \bm{0} & \bm{0} & \bm{0} & \bm{0} & P_{6, 6}
\end{array}
\right), 
\]
where the sizes of square matrices $P_{1, 1}$, $P_{2, 2}$, $P_{3, 3}$, $P_{4, 4}$, $P_{5, 5}$, $P_{6, 6}$ are 
$1$, $r_1$, $\ell - r_1$, $r_2$, $m - \ell - r_2$, $1$, respectively. 
Clearly, ${^\ell \cP}^{r_1}_{r_2}$ is a subgroup of $GL(m + 2, k)$ and $\det(P) = \prod_{i = 1}^6\det(P_{i, i})$.

We denote by ${^\ell} \cQ^{r_1}_{r_2}$ be the set of all matrices $Q$ of $GL(\ell + r_1 + r_2 + 1, k)$ with the form 
\[
 Q = 
\left( 
\begin{array}{c | c | c | c | c}
 Q_{1, 1} & Q_{1, 2} & \bm{0} & \bm{0} & Q_{1, 5} \\ 
\hline 
 \bm{0} & Q_{2, 2} & O        & O        & Q_{2, 5} \\ 
\hline 
 \bm{0} & O        & Q_{3, 3} & Q_{3, 4} & Q_{3, 5} \\ 
\hline 
 \bm{0} & O        & O        & Q_{4, 4}  & Q_{4, 5} \\ 
\hline 
 0        & \bm{0}  & \bm{0} & \bm{0}   & Q_{5, 5} 
\end{array}
\right) , 
\]
where the sizes of square matrices $Q_{1, 1}$, $Q_{2, 2}$, $Q_{3, 3}$, $Q_{4, 4}$, and $Q_{5, 5}$ 
are $r_1$, $\ell - r_1$, $r_1$, $r_2$, and $1$, respectively, and 
\[
Q_{1, 1} \cdot {^t} Q_{3, 3} = Q_{5, 5}^{\oplus r_1}.
\]
Clearly, ${^\ell} \cQ^{r_1}_{r_2}$ becomes a subgroup of $GL(\ell + r_1 + r_2 + 1, k)$.

\paragraph{Equivalence relations of polynomial matrices of $h({^\ell} \W^{r_1}_{r_2})$ and ${^\ell}\mf{H}^{r_1}_{r_2}$}
\quad \medskip 

The following lemma gives a necessary and sufficient condition for two given polynomial matrices 
of $h({^\ell} \W^{r_1}_{r_2})$ to be equivalent. 

\begin{lem}
Let $A(T) = h(\bm{a}_1, \bm{a}_2, \bm{\alpha}_1, \bm{\alpha}_2, \varphi) \in h( {^\ell}\W^{r_1}_{r_2} )$ 
and let $B(T) = h(\bm{b}_1, \bm{b}_2, \bm{\beta}_1, \bm{\beta}_2, \psi) \in h( {^\ell}\W^{r_1}_{r_2} )$.  
Then the following conditions {\rm (1)} and {\rm (2)} are equivalent: 
\begin{enumerate} 
\item[\rm (1)] There exists a regular matrix $P$ of ${^\ell}\cP^{r_1}_{r_2}$ such that $P^{-1} A(T) P = B(T)$. 
\item[\rm (2)] We have 
$(\bm{b}_1, \bm{b}_2, \bm{\beta}_1, \bm{\beta}_2, \psi)
= (\bm{a}_1, \bm{a}_2,  \bm{\alpha}_1, \bm{\alpha}_2, \varphi)\, Q$ for some $Q \in {^\ell \cQ}^{r_1}_{r_2}$. 
\end{enumerate} 
\end{lem}

\Proof We first prove (1) $\Longrightarrow$ (2). Since $(A(T) - I_{m + 2}) P = P(B(T) - I_{m + 2})$, 
we have 
\begin{eqnarray*}
\lefteqn{
\left(
\begin{array}{c | c | c | c | c | c}
 0 & \bm{a}_1 \cdot P_{2, 2} & \bm{a}_1 \cdot P_{2, 3} + \bm{a}_2 \cdot P_{3, 3} & \bm{0}  & \bm{0} & \bm{a}_1 \cdot P_{2, 6} + \bm{a}_2 \cdot P_{3, 6} + \varphi \cdot P_{6, 6} \\ 
\hline 
\bm{0} & O & O & O &O & {^t}\bm{\alpha}_1 \cdot P_{6, 6} \\ 
\hline 
\bm{0} & O & O & O &O & \bm{0} \\ 
\hline 
\bm{0} & O & O & O &O & {^t}\bm{\alpha}_2 \cdot P_{6, 6} \\ 
\hline 
\bm{0} & O & O & O &O & \bm{0} \\ 
\hline 
\bm{0} & O & O & O &O & 0 \\ 
\end{array}
\right) } \\ 
&=& 
\left(
\begin{array}{c | c | c | c | c | c}
 0 & P_{1, 1} \cdot \bm{b}_1 & P_{1, 1} \cdot \bm{b}_2 & \bm{0}  & \bm{0} & P_{1, 1} \cdot \psi + P_{1, 2} \cdot {^t}\bm{\beta}_1 + P_{1, 4} \cdot {^t}\bm{\beta}_2  \\ 
\hline 
\bm{0} & O & O & O &O & P_{2, 2} \cdot {^t}\bm{\beta}_1 \\ 
\hline 
\bm{0} & O & O & O &O & \bm{0} \\ 
\hline 
\bm{0} & O & O & O &O & P_{4, 2} \cdot {^t}\bm{\beta}_1 + P_{4, 4} \cdot {^t}\bm{\beta}_2 \\ 
\hline 
\bm{0} & O & O & O &O & \bm{0} \\ 
\hline 
\bm{0} & O & O & O &O & 0 \\ 
\end{array}
\right) . 
\end{eqnarray*}
Thus we have 
\begin{eqnarray*}
\left\{
\begin{array}{ccl}
 \bm{b}_1 &=& P_{1, 1}^{-1} \cdot \bm{a}_1 \cdot P_{2, 2} , \\ [2mm] 
 \bm{b}_2 & =& P_{1, 1}^{-1} \cdot \bm{a}_1 \cdot P_{2, 3} + P_{1, 1}^{-1} \cdot \bm{a}_2 \cdot P_{3, 3} , \\ [2mm] 
 {^t} \bm{\beta}_1 &=& P_{2, 2}^{-1} \cdot {^t} \bm{\alpha}_1 \cdot P_{6, 6}, \\ [2mm] 
 {^t} \bm{\beta}_2 &=& - P_{4, 4}^{-1} \cdot P_{4, 2} \cdot P_{2, 2}^{-1} \cdot {^t} \bm{\alpha}_1 \cdot P_{6, 6}
                               + P_{4, 4}^{-1} \cdot {^t} \bm{\alpha}_2 \cdot P_{6, 6}, \\  [2mm] 
 \psi &=& P_{1, 1}^{-1} \cdot \bm{a}_1 \cdot P_{2, 6} + P_{1, 1}^{-1} \cdot \bm{a}_2 \cdot P_{3, 6} \\ 
 & & + P_{1, 1}^{-1} \cdot ( - P_{1, 2} \cdot P_{2, 2}^{-1}  +  P_{1, 4} \cdot P_{4, 4}^{-1} \cdot P_{4, 2} \cdot P_{2, 2}^{-1}) \cdot {^t} \bm{\alpha}_1 \cdot P_{6, 6} \\
         &  & - P_{1, 1}^{-1} \cdot P_{1, 4} \cdot P_{4, 4}^{-1} \cdot {^t} \bm{\alpha}_2 \cdot P_{6, 6} 
                + P_{1, 1}^{-1} \cdot \varphi \cdot P_{6, 6} . 
\end{array}
\right. 
\end{eqnarray*}
From these equalities, we have $Q \in {^\ell} \cQ^{r_1}_{r_2}$ such that 
$(\bm{b}_1, \bm{b}_2, \bm{\beta}_1, \bm{\beta}_2, \psi) = (\bm{a}_1, \bm{a}_2, \bm{\alpha}_1, \bm{\alpha}_2, \varphi) \, Q$.

We next prove (2) $\Longrightarrow$ (1). 
We can write the matrix $Q$ as $Q = ( Q_{i, j} )_{1 \leq i, j \leq 5}$, 
where the sizes of square matrices $Q_{1, 1}$, $Q_{2, 2}$, $Q_{3, 3}$, $Q_{4,4}$, and $Q_{5, 5}$ are $r_1$, $\ell - r_1$, 
$r_1$, $r_2$, and $1$, respectively. Let 
\[
P := 
\left(
\begin{array}{c | c | c | c | c |c}
 1 & {^t} ( - Q_{3, 3}^{-1} \cdot Q_{3, 5} + Q_{3, 3}^{-1} \cdot Q_{3, 4} \cdot Q_{4, 4}^{-1} \cdot Q_{4, 5}) 
 & \bm{0} & - {^t} (Q_{4, 4}^{-1} \cdot Q_{4, 5}) & \bm{0} & 0 \\
\hline 
 \bm{0} & Q_{1, 1} & Q_{1, 2} & O & O & Q_{1, 5} \\
\hline 
 \bm{0} & O & Q_{2, 2} & O & O & Q_{2, 5} \\
\hline 
 \bm{0} & - \, {^t}(Q_{3, 4} \cdot Q_{4, 4}^{-1}) \cdot Q_{1, 1} & O & Q_{5, 5} \cdot {^t}(Q_{4, 4}^{-1}) & O & \bm{0} \\ 
\hline 
\bm{0} & O & O & O & I_{m - \ell - r_2} & \bm{0} \\
\hline 
 0 & \bm{0} & \bm{0} & \bm{0} & \bm{0} & Q_{5, 5} 
\end{array}
\right) . 
\]
Then we have $(A(T) - I_{m + 2} ) \cdot P = P \cdot (B(T) - I_{m + 2})$, which implies (1). 
\QED

Given two polynomial matrices 
$A(T) = h(\bm{a}_1, \bm{a}_2, \bm{\alpha}_1, \bm{\alpha}_2, \varphi)$ and 
$B(T) = h(\bm{b}_1, \bm{b}_2, \bm{\beta}_1, \bm{\beta}_2, \psi)$ of ${^\ell \mf{H}}^{r_1}_{r_2}$ 
which are equivalent, 
we shall see in the following lemma that any matrix $P$ of  $GL(m + 2, k)$ satisfying $P^{-1} A(T) P = B(T)$ 
has to belong to ${^\ell} \cP^{r_1}_{r_2}$.

\begin{lem}
Let $A(T) = h(\bm{a}_1, \bm{a}_2, \bm{\alpha}_1, \bm{\alpha}_2, \varphi) \in {^\ell \mf{H}}^{r_1}_{r_2}$ and 
$B(T) = h(\bm{b}_1, \bm{b}_2, \bm{\beta}_1, \bm{\beta}_2, \psi) \in {^\ell \mf{H}}^{r_1}_{r_2}$, 
where $(\bm{a}_1, \bm{a}_2, \bm{\alpha}_1, \bm{\alpha}_2, \varphi) \in {^\ell}\V^{r_1}_{r_2}$ and 
$(\bm{b}_1, \bm{b}_2, \bm{\beta}_1, \bm{\beta}_2, \psi) \in {^\ell}\V^{r_1}_{r_2}$.  
Then the following conditions {\rm (1)}, {\rm (2)} and {\rm (3)} are equivalent: 
\begin{enumerate} 
\item[\rm (1)] $A(T)$ and $B(T)$ are equivalent. 
\item[\rm (2)] There exists a regular matrix $P$ of ${^\ell}\cP^{r_1}_{r_2}$ such that $P^{-1} A(T) P = B(T)$. 
\item[\rm (3)] We have 
$(\bm{b}_1, \bm{b}_2, \bm{\beta}_1, \bm{\beta}_2, \psi)
= (\bm{a}_1, \bm{a}_2,  \bm{\alpha}_1, \bm{\alpha}_2, \varphi)\, Q$ for some $Q \in {^\ell \cQ}^{r_1}_{r_2}$. 
\end{enumerate} 
\end{lem}

\Proof  
By Lemma 2.11, we have only to show the implication (1) $\Longrightarrow$ (2). 
So, assume that $A(T)$ ane $B(T)$ are equivalent, and let $P$ be a regular matrix of $GL(m +2, k)$ such that $P^{-1} A(T) P = B(T)$. 
Let $P = (P_{i, j})_{1 \leq i, j \leq 6}$ be a block matrix, where we assign  
the numbers of rows (and columns) of these blocks $1$, $r_1$, $\ell - r_1$, $r_2$, $m - \ell - r_2$, $1$ from the top to the bottom 
(from the left to the right). 
We know from the proof of Lemma 2.10 that $P$ has the form 
\[
 P = 
\left( 
\begin{array}{ c | c | c | c | c | c }
 P_{1, 1} & P_{1, 2} & P_{1, 3} & P_{1, 4} & P_{1, 5} & P_{1, 6} \\
\hline 
\bm{0} & P_{2, 2} & P_{2, 3} & O & O & P_{2, 6} \\
\hline 
\bm{0} & P_{3, 2} & P_{3, 3} & O & O & P_{3, 6} \\ 
\hline 
\bm{0} & P_{4, 2} & P_{4, 3} & P_{4, 4} & P_{4, 5} & P_{4, 6} \\ 
\hline 
\bm{0} & P_{5, 2} & P_{5, 3} & P_{5, 4} & P_{5, 5} & P_{5, 6} \\ 
\hline 
\bm{0} & \bm{0} & \bm{0} & \bm{0} & \bm{0} & P_{6, 6}
\end{array}
\right) .
\]
Since $(A(T) - I_{m +2} ) P = P (B(T) - I_{m + 2})$, we have 
\begin{eqnarray*}
\lefteqn{ 
\left(
\begin{array}{ c | c | c | c | c | c}
0 & \bm{a}_1 \cdot P_{2, 2} + \bm{a}_2 \cdot P_{3, 2} & \bm{a}_1 \cdot P_{2, 3} + \bm{a}_2 \cdot P_{3, 3} & \bm{0} & \bm{0} & \bm{a}_1 P_{2, 6}  + \bm{a}_2 \cdot P_{3, 6} + \varphi \cdot P_{6, 6} \\ 
\hline 
 \bm{0} & O & O & O & O & {^t} \bm{\alpha}_1 \cdot P_{6, 6} \\ 
\hline 
 \bm{0} & O & O & O & O & \bm{0} \\ 
\hline 
 \bm{0} & O & O & O & O &{^t}  \bm{\alpha}_2 \cdot P_{6, 6} \\ 
\hline 
 \bm{0} & O & O & O & O & \bm{0} \\ 
\hline 
 0        & \bm{0}    & \bm{0}    & \bm{0} & \bm{0} & 0
\end{array}
\right)
} \qquad \quad \\
& = & 
\left(
\begin{array}{c | c | c | c | c | c}
 0 & P_{1, 1} \cdot \bm{b}_1  & P_{1, 1} \cdot \bm{b}_2 & \bm{0} & \bm{0} & P_{1, 1} \cdot \psi + P_{1, 2} \cdot {^t} \bm{\beta}_1 + P_{1, 4} \cdot {^t} \bm{\beta}_2 \\ 
\hline 
 \bm{0} & O & O & O & O & P_{2, 2} \cdot {^t} \bm{\beta}_1 \\ 
\hline 
 \bm{0} & O & O & O & O & P_{3, 2} \cdot {^t}\bm{\beta}_1 \\ 
\hline 
 \bm{0} & O & O & O & O & P_{4, 2}  \cdot {^t} \bm{\beta}_1 + P_{4, 4} \cdot {^t} \bm{\beta}_2 \\
\hline  
 \bm{0} & O & O & O & O & P_{5, 2} \cdot {^t} \bm{\beta}_1 + P_{5, 4} \cdot {^t} \bm{\beta}_2 \\
\hline 
 0        & \bm{0}    & \bm{0}    & \bm{0} & \bm{0} & 0
\end{array}
\right) . 
\end{eqnarray*}
The $(3, 6)$-th entries imply $P_{3, 2} = O$. The $(5, 6)$-th entries imply $P_{5, 2} = O$ and $P_{5, 4} = O$. 
Thus $P$ has the desired form. 
\QED

\paragraph{$\mf{H}(m + 2, k[T]) \curvearrowleft {^\ell}\cP^{r_1}_{r_2}$ 
and ${^\ell}\W^{r_1}_{r_2} \curvearrowleft {^\ell}\cQ^{r_1}_{r_2}$ } 
\quad \medskip

The group ${^\ell}\cP^{r_1}_{r_2}$ acts from right on $\mf{H}(m + 2, k[T])$ by conjugation, i.e., 
$\mf{H}(m + 2, k[T]) \times {^\ell}\cP^{r_1}_{r_2} \to \mf{H}(m + 2, k[T])$, $(A(T), P) \mapsto P^{-1} \cdot A(T) \cdot P$. 
For all $P, P' \in {^\ell}\cP^{r_1}_{r_2}$ and $\bm{w} \in {^\ell}\W^{r_1}_{r_2}$, 
we have $P' \cdot h(\bm{w}) \cdot P \in h({^\ell}\W^{r_1}_{r_2})$. 
So, the subset $h({^\ell}\W^{r_1}_{r_2})$ of $\mf{H}(m + 2, k[T])$ is ${^\ell}\cP^{r_1}_{r_2}$-invariant. 
For any $A(T) \in {^\ell}\mf{H}^{r_1}_{r_2}$ and any $P \in {^\ell}\cP^{r_1}_{r_2}$, 
we have $P^{-1} \cdot A(T) \cdot P = P^{-1} \cdot (A(T) - I_{m + 2}) \cdot P + I_{m + 2} \in {^\ell}\mf{H}^{r_1}_{r_2}$ (see the proof (1) $\Longrightarrow$ (2) of Lemma 2.11 for the calculation $(A(T) - I_{m + 2}) P$). 
So, the subset ${^\ell}\mf{H}^{r_1}_{r_2}$ of $\mf{H}(m + 2, k[T])$ is also ${^\ell}\cP^{r_1}_{r_2}$-invariant.

The group ${^\ell} \cQ^{r_1}_{r_2}$ acts on ${^\ell}\W^{r_1}_{r_2}$ by right translation, i.e., 
${^\ell}\W^{r_1}_{r_2} \times {^\ell} \cQ^{r_1}_{r_2} \to {^\ell}\W^{r_1}_{r_2}$, $(\bm{w}, Q) \mapsto \bm{w} Q$. 
The subset ${^\ell}\V^{r_1}_{r_2}$ of ${^\ell}\W^{r_1}_{r_2}$ is ${^\ell}\cQ^{r_1}_{r_2}$-invariant (by the definitions of ${^\ell}\V^{r_1}_{r_2}$ and ${^\ell}\cQ^{r_1}_{r_2}$).

The bijection $h :  {^\ell}\W^{r_1}_{r_2} \to  h( {^\ell}\W^{r_1}_{r_2}) $ induces a bijection 
$h :  {^\ell}\W^{r_1}_{r_2}/{^\ell} \cQ^{r_1}_{r_2} \to  h( {^\ell}\W^{r_1}_{r_2})/{^\ell}\cP^{r_1}_{r_2} $ (see Lemma 2.11), 
and the bijection $h : {^\ell}\V^{r_1}_{r_2}  \to  {^\ell}\mf{H}^{r_1}_{r_2} $ induces a bijection 
$h : {^\ell}\V^{r_1}_{r_2}/{^\ell} \cQ^{r_1}_{r_2} \to  {^\ell}\mf{H}^{r_1}_{r_2}/{^\ell}\cP^{r_1}_{r_2} $ (see Lemma 2.12).

The following lemma in particular shows a parametrization of ${^\ell}\mf{H}^{r_1}_{r_2}/\sim$. 

\begin{lem} 
The following diagram is commutative, where the upward arrows are natural injections and 
the notation $\twoheadrightarrow$ means surjective: 
\begin{eqnarray*}
\xymatrix@C=24pt@R=35pt{                               & \mf{H}(m + 2, k[T])/{^\ell}\cP^{r_1}_{r_2}  \ar[r]^{\twoheadrightarrow} 
 &  \mf{H}(m + 2, k[T])/\sim  \\
  {^\ell}\W^{r_1}_{r_2}/{^\ell} \cQ^{r_1}_{r_2} \ar[r]^(.45)h               &  h( {^\ell}\W^{r_1}_{r_2})/{^\ell}\cP^{r_1}_{r_2} \ar[u] \ar[r]^{\twoheadrightarrow}  & 
h( {^\ell}\W^{r_1}_{r_2})/\sim \ar[u] \\
 {^\ell}\V^{r_1}_{r_2}/{^\ell} \cQ^{r_1}_{r_2} \ar[r]^(.45)h \ar[u]     &  {^\ell}\mf{H}^{r_1}_{r_2}/{^\ell}\cP^{r_1}_{r_2} \ar[u] \ar[r]^{=} & 
{^\ell}\mf{H}^{r_1}_{r_2}/\sim \ar[u] 
\\ 
}
\end{eqnarray*}
\end{lem}

\subsubsection{$\daleth_{m + 2}$}

In this susubsection, we consider an equivalence relation of two polynomial matrices belonging to just 
one of the five sets $\daleth_{m + 2}^{\rm II}$, $\daleth_{m + 2}^{\rm III}$, $\daleth_{m + 2}^{\rm IV}$, 
$\daleth_{m + 2}^{\rm V}$, $\daleth_{m + 2}^{\rm VI}$ (see assertion (2) of Theorem 2.5).

\paragraph{$\daleth_{m + 2}^{\rm II} ( = \bigcup_{1 \leq j \leq m + 1} \mf{A}(1, j)^\circ )$}
\quad \medskip

\begin{thm}
Let $A(T) = \Lambda(1, j; \alpha(T)) \in \daleth_{m + 2}^{\rm II}$ 
and $B(T) = \Lambda(1, j'; \beta(T))\in \daleth_{m + 2}^{\rm II}$, where $1 \leq j, j' \leq m + 1$. 
Then the following conditions {\rm (1)} and {\rm (2)} are equivalent: 
\begin{enumerate}
\item[\rm (1)] $A(T)$ and $B(T)$ are equivalent. 
\item[\rm (2)] $j = j'$, and $\beta(T) = \alpha(T) Q$ for some $Q \in GL(j, k)$. 
\end{enumerate} 
\end{thm}

\Proof 
We shall prove (1) $\Longrightarrow$ (2). 
By Lemma 1.32, $j = j'$. By Lemma 1.45, $\beta(T) = \alpha(T) Q$ for some $Q \in GL(j, k)$. 
The implication (2) $\Longrightarrow$ (1) follows from Lemma 1.45. 
\QED

\paragraph{$\daleth_{m + 2}^{\rm III} ( = \bigcup_{2 \leq i \leq m + 1} \mf{A}(i, 1)^\circ )$} \quad \medskip

\begin{thm}
Let $A(T) = \Lambda(i, 1; \alpha(T)) \in \daleth_{m + 2}^{\rm II}$ and 
$B(T) = \Lambda(i', 1; \beta(T)) \in \daleth_{m + 2}^{\rm II}$, where $2 \leq i, i' \leq m + 1$. 
Then the following conditions {\rm (1)} and {\rm (2)} are equivalent: 
\begin{enumerate}
\item[\rm (1)] $A(T)$ and $B(T)$ are equivalent. 
\item[\rm (2)] $i = i'$, and $\beta(T) = Q \cdot \alpha(T)$ for some $Q \in GL(i, k)$. 
\end{enumerate} 
\end{thm}

\Proof The proof follows from Theorem 2.14 and Lemmas 1.6 and 1.7. 
\QED

\paragraph{$\daleth_{m + 2}^{\rm IV} ( = {^1}\daleth_1^0 )$} \quad \medskip

Let $R$ be a commutative ring with unity and let $n \geq 1$ be an integer. 
We denote by $GLDD(n, R)$ the set of all matrices $D = (d_{i, j})$ of $GL(n, R)$ with one of the following forms 
(1) and (2): 
\begin{enumerate} 
\item[\rm (1)] $D$ is a diagonal matrix. 
\item[\rm (2)] $d_{i, j} = 0$ except for all $(i, j)$'s satisfying $i + j = n + 1$.  
\end{enumerate} 
Clearly, $GLDD(n, R)$ becomes a subgroup of $GL(n, R)$.

\begin{lem}
Let $\alpha(T) = L(a_1 \mid a_2 \mid 0) \in {^1}\mf{L}_1^0$ 
and $\beta(T) = L(b_1 \mid b_2 \mid 0) \in {^1}\mf{L}_1^0$. 
Then the following conditions {\rm (1)} and {\rm (2)} are equivalent: 
\begin{enumerate}
\item[\rm (1)] $Q \cdot \alpha(T) \cdot R = \beta(T)$ for some $Q, R \in GL(2, k)$. 
\item[\rm (2)] $(b_1, b_2) = (a_1, a_2) \, g $ for some $g \in GLDD(2, k)$. 
\end{enumerate}
\end{lem} 

\Proof The proof is straightforward. 
\QED

\begin{thm} 
Let $A(T) = \Lambda(2, 2; \alpha(T)) \in {^1}\daleth_1^0$ 
and let $B(T) = \Lambda(2, 2; \beta(T)) \in {^1}\daleth_1^0$, where 
$\alpha(T) = L(a_1 \mid a_2 \mid 0) \in {^1}\mf{L}_1^0$ 
and $\beta(T) = L(b_1 \mid b_2 \mid 0) \in {^1}\mf{L}_1^0$.   
Then the following conditions {\rm (1)} and {\rm (2)} are equivalent: 
\begin{enumerate}
\item[\rm (1)] $A(T)$ and $B(T)$ are equivalent. 
\item[\rm (2)] $(b_1, b_2) = (a_1, a_2) \, g $ for some $g \in GLDD(2, k)$. 
\end{enumerate} 
\end{thm}

\Proof The proof follows from Lemmas 1.45 and 2.16. 
\QED

\paragraph{$\daleth_{m + 2}^{\rm V}$ ( = $({^1}\daleth_1)^{\dim = 3}$)} \quad \medskip

Let $R$ be a commutative ring with unity. 
For $i_1, i_2, i_3 \geq 1$, we denote by $V(i_1, i_2, i_3 ; R)$ the set of all matrices $g \in GL(i_1 + i_2 + i_3, R)$ with the (V-)form 
\[
g = 
\left(
\begin{array}{ccc}
 g_{1, 1} & O        & g_{1, 3} \\
 O        & g_{2, 2}  & g_{2, 3} \\
 O        & O         & g_{3, 3}
\end{array}
\right)  
\qquad 
(\,  g_{1, 1} \in GL(i_1, R), \quad  g_{2, 2} \in GL(i_2, R), \quad  g_{3, 3} \in GL(i_3, R) \, ) . 
\]

\begin{lem}
Let $\alpha(T) = L(a_1 \mid a_2 \mid a) \in ({^1}\mf{L}_1)^{\dim = 3}$ 
and $\beta(T) = L(b_1 \mid b_2 \mid b) \in ({^1}\mf{L}_1)^{\dim = 3}$. 
Then the following conditions {\rm (1)} and {\rm (2)} are equivalent: 
\begin{enumerate}
\item[\rm (1)]$Q \cdot \alpha(T) \cdot R = \beta(T)$ for some $Q, R \in GL(2, k)$. 
\item[\rm (2)] $(b_1, b_2, b) = (a_1, a_2, a) \, g$ for some $g \in V(1, 1, 1; k)$. 
\end{enumerate}
\end{lem}

\Proof 
Write 
\[
 Q = \left( \begin{array}{cc} q_1 & q_2 \\ q_3 & q _4 \end{array} \right) \quad (\, q_1, q_2, q_3, q_4 \in k \,) , \qquad 
 R = \left( \begin{array}{cc} r_1 & r_2 \\ r_3 & r _4 \end{array} \right) \quad (\, r_1, r_2, r_3, r_4 \in k \,) .  
\]
The implication (1) $\Longrightarrow$ (2) is clear by observing $q_3 = r_3 = 0$. The converse implication (2) $\Longrightarrow$ (1) 
is straightforward. 
\QED

\begin{thm} 
Let $A(T) = \Lambda(2, 2; \alpha(T)) \in ({^1}\daleth_1)^{\dim = 3}$ 
and let $B(T) = \Lambda(2, 2; \beta(T)) \in ({^1}\daleth_1)^{\dim = 3}$, where 
$\alpha(T) = L(a_1 \mid a_2 \mid a) \in ({^1}\mf{L}_1)^{\dim = 3}$ 
and $\beta(T) = L(b_1 \mid b_2 \mid b) \in ({^1}\mf{L}_1)^{\dim = 3}$.   
Then the following conditions {\rm (1)} and {\rm (2)} are equivalent: 
\begin{enumerate}
\item[\rm (1)] $A(T)$ and $B(T)$ are equivalent. 
\item[\rm (2)] $(b_1, b_2, b) = (a_1, a_2, a) \, g$ for some $g \in V(1, 1, 1; k)$. 
\end{enumerate} 
\end{thm}

\Proof The proof follows from Lemmas 1.45 and 2.18. 
\QED

\paragraph{$\daleth_{m + 2}^{\rm VI} ( = \bigcup_{3 \leq i + j \leq m} {^i}\daleth_j )$} \quad \medskip

We prepare the following Lemma 2.20 and 2.21 for considering an equivalence relation of two 
polynomial matrices of $\daleth_{m + 2}^{\rm VI}$ (see Theorem 2.22). 

\begin{lem}
Let $i, j$ be positive integers. 
Assume that 
\[
Q \cdot L(a_1,\ldots, a_i \mid a_{i + 1}, \ldots, a_{i + j} \mid a) \cdot R
 = L(b_1, \ldots, b_i \mid b_{i + 1}, \ldots, b_{i + j} \mid b), 
\]
where $Q \in GL(j + 1, k)$ and $R \in GL(i + 1, k)$. 
Write $Q, R$ as 
\begin{eqnarray*}
\left\{ 
\begin{array}{ll}
Q = 
\left( 
\begin{array}{c | c}
 Q_1 & Q_2 \\
\hline 
 Q_3 & Q_4 
\end{array}
\right)
 \qquad
& (\, Q_1 \in \Mat(1, k), \quad Q_4 \in \Mat(j, k)  \, ) , 
\\ [6mm] 
R = 
\left( 
\begin{array}{c | c}
 R_1 & R_2 \\
\hline 
 R_3 & R_4 
\end{array}
\right)
& (\, R_1 \in \Mat(i, k), \quad R_4 \in \Mat(1, k)  \, ) . 
\end{array}
\right. 
\end{eqnarray*}
Then the following assertions hold true:  
\begin{enumerate}
\item[\rm (1)] $Q_3 = O$ if and only if $R_3 = O$. 
\item[\rm (2)] If $i + j \geq 3$, then $Q_3 = O$. 
\end{enumerate}
\end{lem}

\Proof (1) Since 
\[
Q \cdot L(a_1,\ldots, a_i \mid a_{i + 1}, \ldots, a_{i + j} \mid a)
 = L(b_1, \ldots, b_i \mid b_{i + 1}, \ldots, b_{i + j} \mid b) \cdot R^{-1}, 
\]
by writing $R^{-1}$ as 
\[
R^{-1}= 
\left( 
\begin{array}{c | c}
 R_1' & R_2' \\
\hline 
 R_3' & R_4' 
\end{array}
\right)
\qquad  (\, R_1' \in \Mat(i, k), \quad R_4' \in \Mat(1, k)  \, ) , 
\]
we obtain the following four equalities: 
\begin{eqnarray*}
(\ast) \qquad 
\left\{
\begin{array}{l}
 Q_1 \cdot ( a_1, \ldots, a_i ) = ( b_1, \ldots, b_i ) \cdot R_1' + b \cdot R_3', \\
 Q_1 \cdot (a) + Q_2 \cdot {^t} (a_{i + 1}, \ldots, a_{i + j}) = (b_1, \ldots, b_i) \cdot R_2' + b \cdot R_4', \\
 Q_3 \cdot (a_1, \ldots, a_i) = {^t}(b_{i + 1}, \ldots, b_{i + j}) \cdot R_3', \\
 Q_3 \cdot (a) + Q_4 \cdot {^t}(a_{i + 1}, \ldots, a_{i + j}) = {^t}(b_{i + 1}, \ldots, b_{i + j}) \cdot R_4' . 
\end{array}
\right.
\end{eqnarray*}
Note that $R_3 = O$ if and only if $R_3' = O$. The proof is clear from the third equality of $(\ast)$. 

(2) We prove the assertion by contraposition. So, assume $Q_3 \ne O$. 
Write $Q_3 = {^t}(q_1, \ldots, q_j)$, where $q_1, \ldots, q_j \in k$. 
We have $q_\lambda \ne 0$ for some $1 \leq \lambda \leq j$. 
Since $q_\lambda \cdot (a_1, \ldots, a_i) = b_{i + \lambda} \cdot R_3'$ (see the third equality of $(\ast)$) and 
$a_1, \ldots, a_i$ are linearly independent over $k$, we have $i = 1$. 
Now, the third equality $Q_3 \cdot (a_1) = {^t} (b_2, \ldots, b_{1 + j}) \cdot R_3'$ implies $j = 1$ (since 
$b_2, \ldots, b_{1 + j}$ are linearly independent over $k$). 
Thus, we have $i + j = 2$. 
\QED

\begin{lem} 
Let $i, j$ be positive integers satisfying $i + j \geq 3$. 
Let $\alpha(T) = L(a_1, \ldots, a_i \mid a_{i + 1}, \ldots, a_{i + j} \mid a) \in {^i}\mf{L}_j$ 
and $\beta(T) = L(b_1, \ldots, b_i \mid b_{i + 1}, \ldots, b_{i + j} \mid b) \in {^i}\mf{L}_j$.  
Then the following conditions {\rm (1)} and {\rm (2)} are equivalent: 
\begin{enumerate}
\item[\rm (1)] $Q \cdot \alpha(T) \cdot R = \beta(T)$ for some $Q \in GL(j + 1, k)$ and $R \in GL(i + 1, k)$. 
\item[\rm (2)] $(b_1, \ldots, b_i, b_{i + 1}, \ldots, b_{i + j}, b) = (a_1, \ldots, a_i, a_{i + 1}, \ldots, a_{i + j}, a) \, g$ for some $g \in V(i, j, 1; k)$. 
\end{enumerate}
\end{lem}

\Proof We first prove (1) $\Longrightarrow$ (2). 
Since $i + j \geq 3$, we know from Lemma 2.20 that $Q$ and $R$ has the following form: 
\begin{eqnarray*}
\left\{ 
\begin{array}{ll}
Q = 
\left( 
\begin{array}{c | c}
 Q_1 & Q_2 \\
\hline 
 O & Q_4 
\end{array}
\right)
 \qquad
& (\, Q_1 \in \Mat(1, k), \quad Q_4 \in \Mat(j, k)  \, ) , 
\\ [6mm] 
R = 
\left( 
\begin{array}{c | c}
 R_1 & R_2 \\
\hline 
 O & R_4 
\end{array}
\right)
& (\, R_1 \in \Mat(i, k), \quad R_4 \in \Mat(1, k)  \, ) . 
\end{array}
\right. 
\end{eqnarray*}
Since $Q \cdot \alpha(T) \cdot R = \beta(T)$, we have the desired form in the condition (2). 

We next prove (2) $\Longrightarrow$ (1). 
We can write $g$ as 
\[
g = 
\left(
\begin{array}{ccc}
 g_{1, 1} & O        & g_{1, 3} \\
 O        & g_{2, 2}  & g_{2, 3} \\
 O        & O         & g_{3, 3}
\end{array}
\right)  
\qquad 
(\,  g_{1, 1} \in GL(i, k), \quad  g_{2, 2} \in GL(j, k), \quad  g_{3, 3} \in k \backslash \{0 \} \, ) . 
\]
Let 
\begin{eqnarray*}
Q := 
\left(
\begin{array}{c | c}
 1 & {^t}g_{2, 3} \\
\hline 
 \bm{0} & g_{3, 3}^{-1} \cdot {^t}g_{2, 2} 
\end{array}
\right)
\in GL(j + 1, k) 
\quad 
\text{ and }
\quad 
R := 
\left(
\begin{array}{c | c}
 g_{1, 1} & g_{1, 3} \\
\hline 
 \bm{0} & g_{3, 3}
\end{array}
\right) 
\in GL(i + 1, k).
\end{eqnarray*}
Then we have $Q \cdot \alpha(T) \cdot R = \beta(T)$. 
\QED

\begin{thm}
Let $i, j$ be positive integers satisfying $3 \leq i + j \leq m$. 
Let $A(T) = \Lambda(j + 1, i + 1; \alpha(T) ) \in \daleth_{m + 2}^{\rm VI}$ 
and let $B(T) = \Lambda( j' + 1, i' + 1; \beta(T) ) \in \daleth_{m + 2}^{\rm VI}$, where 
$\alpha(T) = L(a_1, \ldots, a_i \mid a_{i + 1}, \ldots, a_{i + j} \mid a) \in {^i}\mf{L}_j$ 
and $\beta(T) = L(b_1, \ldots, b_{i'} \mid b_{i' + 1}, \ldots, b_{i' + j'} \mid b) \in {^{i'}}\mf{L}_{j'}$.   
Then the following conditions {\rm (1)} and {\rm (2)} are equivalent: 
\begin{enumerate}
\item[\rm (1)] $A(T)$ and $B(T)$ are equivalent. 
\item[\rm (2)] $(i, j) = (i', j')$, and $(b_1, \ldots, b_i, b_{i + 1}, \ldots, b_{i + j}, b) = (a_1, \ldots, a_i, a_{i + 1}, \ldots, a_{i + j}, a) \, g$ for some $g \in V(i, j, 1)$. 
\end{enumerate} 
\end{thm}

\Proof The proof follows from Lemmas 1.32, 1.45 and 2.21. 
\QED

\section{Exponential matrices belonging to Heisenberg groups}

We have the following theorem (see Theorem 2.1): 

\begin{thm} 
The following assertions {\rm (1)} and {\rm (2)} hold true: 
\begin{enumerate}
\item[\rm (1)] $\ds \mf{H}(m + 2, k[T])^E
 \leadsto (\daleth_{m +2})^E \cup  \bigcup_{(\ell, r_1, r_2) \in \Omega_m} ( {^\ell \mf{H}}^{r_1}_{r_2} )^E$. 
\item[\rm (2)] Two sets $( \daleth_{m + 2} )^E$ and $\ds \bigcup_{(\ell, r_1, r_2) \in \Omega_m} ({^\ell \mf{H}}^{r_1}_{r_2})^E$ 
are mutually $GL(m + 2, k)$-disjoint. 
\end{enumerate} 
\end{thm}

So, we are interested in exponential matrices of belonging to any one of sets $(\daleth_{m +2})^E$ and 
 $( {^\ell \mf{H}}^{r_1}_{r_2} )^E$, where $(\ell, r_1, r_2) \in \Omega_m$. 
We consider $( {^\ell \mf{H}}^{r_1}_{r_2} )^E$ in Subsection 3.2, and consider $(\daleth_{m +2})^E$ in Subsection 3.3. 
Subsequently, we consider an equivalence relation of two exponential matrices belonging to 
any one of sets $( {^\ell \mf{H}}^{r_1}_{r_2} )^E$ and $(\daleth_{m +2})^E$ in Subsection 3.4.

\subsection{A necessary and sufficient condition for a polynomial matrix $\eta(x_1, \ldots, x_m, y_1, \ldots, y_m, z)$ to be exponential}

The following lemma gives a necessary and sufficient condition for a polynomial matrix 
of a Heisenberg group to be an exponential matrix.

\begin{lem} 
Let $x_1, \ldots, x_m, y_1, \ldots, y_m, z \in k[T]$. 
Then the following conditions are equivalent: 
\begin{enumerate} 
\item[\rm (1)] $\eta(x_1, \ldots, x_m, y_1, \ldots, y_m, z)$ is an exponential matrix. 
\item[\rm (2)] $x_1, \ldots, x_m, y_1, \ldots, y_m$ are $p$-polynomials, 
and the equality 
\[
z(T + T') - z(T) - z(T') = \sum_{i = 1}^m x_i(T) \cdot y_i(T') 
\]
holds true, where $T, T'$ are indeterminates over $k$. 
\end{enumerate} 
\end{lem}

\Proof The proof follows from the definition of an exponential matrix. 
\QED

In the following Subsection 3.2, 
we shall study $p$-polynomials $x_1, \ldots, x_m, y_1, \ldots, y_m$ and a polynomial $z$ satisfying the above condition (2).  
Since $z(T + T') - z(T) - z(T')$ is a symmetric polynomial of $k[T, T']$, where the transposition $\sigma = (1 \; \, 2)$ of the  
symmetric group $\mf{S}_2$ of order $2$ acts on $k[T, T']$ as $\sigma(T) := T'$ and $\sigma(T') := T$. 
Therefore, $\sum_{i = 1}^m x_i(T) \cdot y_i(T')$ is also a symmetric polynomial. 
Based on this symmetry, we can describe any exponential matrix of ${^\ell} \mf{H}^{r_1}_{r_2}$ 
(see the following Theorem 3.3).

\subsection{$({^\ell \mf{H}}^{r_1}_{r_2})^E$}

\begin{thm} 
Let $A(T) = \eta(a_1, \ldots, a_m, a_{m + 1}, \ldots, a_{2m},  a_{2m + 1}) \in {^\ell \mf{H}}^{r_1}_{r_2}$.  
Then $A(T) \in ({^\ell \mf{H}}^{r_1}_{r_2})^E$ 
if and only if the following conditions {\rm (1), (2), (3)} hold true: 
\begin{enumerate}

\item[\rm (1)] The polynomials $a_i$ $(1 \leq i \leq 2m)$ are $p$-polynomials.

\item[\rm (2)] There exists a unique regular symmetric matrix $S = (s_{i, j})_{1 \leq i, j \leq r_1}\in GL(r_1, k)$ such that 
\[
(a_{m + 1}, \ldots, a_{m + r_1}) =  (a_1, \ldots, a_{r_1}) \, S, 
\]
and furthermore if $p = 2$, the diagonal entires $s_{i, i}$ $(1 \leq i \leq r_1)$ of $S$ are zeroes.

\item[\rm (3)]
There exists a unique $p$-polynomial $\alpha(T) \in k[T]$ such that 
\begin{eqnarray*}
a_{2m + 1}(T) = \left\{ 
\begin{array}{ll}
\ds \alpha(T) + \sum_{1 \leq i < j \leq r_1} s_{i, j} \cdot a_i(T) \cdot a_j(T)  & \quad \text{ if \; $p = 2$},  \\ [6mm] 
\ds \alpha(T) +  \frac{1}{2}  \cdot \bm{a}(T) \cdot S \cdot {^t} \bm{a}(T) \qquad & \quad  \text{ if  \; $p \geq 3$},
\end{array} 
\right. 
\end{eqnarray*}
where $\bm{a}(T) := (a_1(T), \ldots, a_{r_1}(T)) \in k[T]^{r_1}$.  
\end{enumerate} 
\end{thm}

For proving the above theorem, we prepare Subsubsections 3.2.1, 3.2.2, 3.2.3.  
Subsequently, in Subsubsection 3.2.4, we give a proof of Theorem 3.3.

\subsubsection{On a symmetric polynomial $\ds \sum_{i = 1}^n a_i (T) \cdot b_i (T')$ of $k[T, T']$}

The following theorem is crucial on determining the form of an exponential matrix belonging to Heisenberg groups:

\begin{thm} 
Let $n \geq 1$ be an integer and let $a_1, \ldots, a_n, b_1, \ldots, b_n \in k[T]$. 
Assume that $a_1, \ldots, a_n$ are linearly independent over $k$. 
Then the following conditions {\rm (1)} and {\rm (2)} are equivalent: 
\begin{enumerate} 
\item[\rm (1)] $\ds\sum_{i  = 1}^n a_i(T) \cdot b_i(T') $ is a symmetric polynomial of $k[T, T']$. 
\item[\rm (2)] There exists a symmetric matrix $S \in \Mat(n, k)$ such that $(b_1, \ldots, b_n ) = (a_1, \ldots, a_n) S$. 
\end{enumerate} 
\end{thm}

We shall give a remark concerning Theorem 3.4. 
Let $a_1(T) := T$, $a_2(T) := T^2$, $a_3(T) :=0$,  $b_1(T) := T$, $b_2(T) := 0$, $b_3(T) := T^3$. 
Then $\sum_{i = 1}^3 a_i(T) \cdot b_i(T') = T \cdot T'$ is symmetric.  
However, there is no symmetric matrix $S \in \Mat(3, k)$ such that $(b_1, b_2, b_3) = (a_1, a_2, a_3) S$. 
So, we cannot delete the assumption that $a_1, \ldots, a_n$ are linearly independent over $k$, from Theorem 3.4.

In order to prove the above Theorem 3.4, we prepare 3.2.1.1, 3.2.1.2, 3.2.1.3, and then prove Theorem 3.4 in 3.2.1.4.

\paragraph{Lemmas}
\quad \medskip

For a polynomial $a(T) \in k[T]$, we define the {\it order} $\ord(a(T))$ of $a(T)$ as 
\begin{eqnarray*}
\ord(a(T)) 
:=
\left\{
\begin{array}{ll}
 \max\{d \in \mathbb{Z}_{\geq 0} \mid a(T) \in T^d \cdot k[T] \} \qquad & \text{ if \quad  $a(T) \ne 0$}, \\
 \infty & \text{ if \quad $a(T) = 0$} . 
\end{array} 
\right.
\end{eqnarray*}
We use the symbol $\infty > d$ for any integer $d$.

\begin{lem}
Let $k[T, T']$ be a polynomial ring in two variables over $k$. 
Assume that polynomials $a_1, \ldots, a_n, b_1, \ldots, b_n \in k[T]$ satisfy the following conditions {\rm (i), (ii), (iii)}: 
\begin{enumerate} 
\item[\rm (i)] $\ds \sum_{i = 1}^n a_i(T) \cdot b_i(T')$ is a symmetric polynomial of $k[T, T']$. 
\item[\rm (ii)] $a_1 \ne 0$. 
\item[\rm (iii)] $a_1$ has the minimum order among $a_1, \ldots, a_n, b_1, \ldots, b_n$, i.e., 
\[
 \ord(a_1) \leq \min\{ \ord(a_1), \, \ldots, \ord(a_n),\, \ord(b_1), \, \ldots, \, \ord(b_n) \} . 
\]
\end{enumerate} 
Then the following assertions {\rm (1), (2), (3)} hold true: 
\begin{enumerate} 
\item[\rm (1)] There exist $\lambda_i \in k$ $(2 \leq i \leq n)$ and $\mu_i \in k$ $(1 \leq i \leq n)$ 
such that 
\begin{eqnarray*}
\left\{ 
\begin{array}{ll}
\ord(a_i - \lambda_i \cdot a_1) > \ord(a_1) \quad & (2 \leq i \leq n), \\
\ord(b_i - \mu_i \cdot a_1) > \ord(a_1) & (1 \leq i \leq n). 
\end{array}
\right. 
\end{eqnarray*} 
\item[\rm (2)] Let $A_i := a_i - \lambda_i \cdot a_1$ $(2 \leq i \leq n)$ and let $B_i := b_i - \mu_i \cdot a_1$ $(1 \leq i \leq n)$. 
Then the polynomial  
\[
\ds\sum_{i = 2}^n A_i(T) \cdot B_i(T')
\]
is a symmetric polynomial of $k[T, T']$. 
\item[\rm (3)] If there exists a symmetric matrix $S \in \Mat(n - 1, k)$ such that $(B_2, \ldots, B_n) = (A_2, \ldots, A_n) S$, 
then there exists a symmetric matrix $\wt{S} \in \Mat(n, k)$ such that $(b_1, \ldots, b_n) = (a_1, \ldots, a_n) \wt{S}$. 
\end{enumerate}
\end{lem}

\Proof Assertion (1) is clear from the conditions (ii), (iii). 
We shall prove assertion (2). 
Note that 
\begin{eqnarray*}
\lefteqn{ \sum_{i = 1}^n a_i(T) \cdot b_i(T')}\\
 & = & a_1(T) \cdot b_1(T') + \sum_{i = 2}^n a_i(T) \cdot b_i(T') \\
 & = & a_1(T) \cdot (B_1(T') + \mu_1 \cdot a_1(T')) 
  + \sum_{i = 2}^n (A_i(T) + \lambda_i \cdot a_1(T)) \cdot (B_i(T') + \mu_i \cdot a_1(T')) \\ 
 & = & a_1(T) \cdot \left(  B_1(T')  +  \sum_{i = 2}^n \lambda_i \cdot B_i(T') \right)  + a_1(T') \cdot \sum_{i = 2}^n \mu_i \cdot A_i(T) 
 + \sum_{i = 2}^n A_i(T) \cdot B_i(T') \\
 & &  + \left( \mu_1 + \sum_{i = 2}^n \lambda_i \cdot \mu_i \right) a_1(T) \cdot a_1(T') .  
\end{eqnarray*}
By the assumption (i), we have 
\begin{eqnarray*}
\lefteqn{
a_1(T) \cdot \left(  B_1(T')  +  \sum_{i = 2}^n \lambda_i \cdot B_i(T') \right)  + a_1(T') \cdot \sum_{i = 2}^n \mu_i \cdot A_i(T) 
 + \sum_{i = 2}^n A_i(T) \cdot B_i(T') } \\
&=&
a_1(T') \cdot \left(  B_1(T)  +  \sum_{i = 2}^n \lambda_i \cdot B_i(T) \right)  + a_1(T) \cdot \sum_{i = 2}^n \mu_i \cdot A_i(T') 
 + \sum_{i = 2}^n A_i(T') \cdot B_i(T) . 
\end{eqnarray*}
Let $m(T)$ be the monomial appearing in $a_1(T)$ whose order is $\ord(a_1(T))$. 
Comparing the coefficients of the monomial $m(T)$ appearing in the both sides of the above equality, we have 
\[
 B_1(T')  +  \sum_{i = 2}^n \lambda_i \cdot B_i(T')  = \sum_{i = 2}^n \mu_i \cdot A_i(T') , 
\]
and thereby have 
\[
  \sum_{i = 2}^n A_i(T) \cdot B_i(T')  =\sum_{i = 2}^n A_i(T') \cdot B_i(T) . 
\]
So, assertion (2) is proved. We shall prove assertion (3). Write $S = (s_{\alpha, \beta})_{2 \leq \alpha, \beta \leq n}$, where $s_{\alpha, \beta} \in k$ $(2 \leq \alpha, \beta\leq n)$. 
Since $B_i = \sum_{j = 2}^n s_{j, i} \, A_j$ for all $2 \leq i \leq n$, we have 
\begin{eqnarray*}
 b_i = \left( \mu_i - \sum_{j = 2}^n s_{j, i} \cdot \lambda_j \right) \cdot a_1 + \sum_{j = 2}^n s_{j, i} \cdot a_j 
 \qquad (2 \leq \forall i \leq n).  
\end{eqnarray*}
Note that 
\begin{eqnarray*}
B_1(T)
 &=& \sum_{i = 2}^n \mu_i \cdot A_i(T) - \sum_{i = 2}^n \lambda_i \cdot B_i(T) \\
 &=& \sum_{i = 2}^n \mu_i \cdot (a_i - \lambda_i \cdot a_1) - \sum_{i = 2}^n \lambda_i \cdot (b_i - \mu_i \cdot a_1) \\
 &=& \sum_{i = 2}^n \mu_i \cdot a_i - \sum_{i = 2}^n \lambda_i \cdot b_i . 
\end{eqnarray*}
Since $B_1 = b_1 - \mu_1 \cdot a_1$, we have 
\begin{eqnarray*}
 b_1
 &=& \mu_1 \cdot a_1 + \sum_{i = 2}^n \mu_i \cdot a_i - \sum_{i = 2}^n \lambda_i \cdot b_i \\
 &=& \left( \mu_1 - \sum_{i = 2}^n  \lambda_i \cdot ( \mu_i - \sum_{j = 2}^n s_{j, i} \cdot \lambda_j) \right) \cdot a_1 
 + \sum_{i = 2}^n \left(\mu_i - \sum_{j = 2}^n \lambda_j \cdot s_{i, j}  \right) \cdot a_i .  
\end{eqnarray*}  
Now we find a matrix $\wt{S} = (\wt{s}_{\alpha, \beta})_{1 \leq \alpha, \beta \leq n} \in \Mat(n, k)$ 
so that $(b_1, \ldots, b_n) = (a_1, \ldots ,a_n) \wt{S}$. In fact, let 
\begin{eqnarray*}
\wt{s}_{\alpha, \beta} 
:= 
\left\{ 
\begin{array}{ll}
 \ds \mu_1 - \sum_{i = 2}^n  \lambda_i \cdot ( \mu_i - \sum_{j = 2}^n s_{j, i} \cdot \lambda_j)  \qquad & \text{if \quad $\alpha = 1$, \;  $\beta =1$},  \\
 \ds \mu_\beta - \sum_{j = 2}^n s_{j, \beta} \cdot  \lambda_j  & \text{if \quad $\alpha = 1$, \; $2 \leq \beta \leq n$}.  \\
 \ds \mu_\alpha - \sum_{j = 2}^n \lambda_j \cdot s_{\alpha, j}    & \text{if \quad $2 \leq \alpha \leq n$, \; $\beta = 1$},   \\
 s_{\alpha, \beta}    & \text{if \quad $2 \leq \alpha \leq n$, \; $2 \leq \beta \leq n$} .   \\
\end{array} 
\right. 
\end{eqnarray*}
Since $S$ is symmetric, $\wt{S}$ is also symmetric. 
\QED

\begin{lem}
Let $\sigma$ be a permutation of the set $\{1, 2, \ldots, n\}$. 
Assume that $a_1, \ldots, a_n, b_1, \ldots, b_n \in k[T]$ satisfy 
$(b_{\sigma(1)}, \ldots, b_{\sigma(n)}) = (a_{\sigma(1)}, \ldots, a_{\sigma(n)}) S$ for some symmetric matrix $S \in \Mat(n, k)$. 
Then there exists a symmetric matrix $S' \in \Mat(n, k)$ such that $(b_1, \ldots, b_n) = (a_1, \ldots, a_n) S'$. 
\end{lem}

\Proof The proof is  straightforward. 
\QED

\begin{lem}
Assume that $a_1, \ldots, a_n, b_1, \ldots, b_n \in k[T]$ satisfy the following conditions {\rm (i), (ii), (iii)}: 
\begin{enumerate} 
\item[\rm (i)] $\ds\sum_{i = 1}^n a_i(T) \cdot b_i(T')$ is a symmetric polynomial of $k[T, T']$. 
\item[\rm (ii)] $a_1, \ldots, a_n$ are linearly independent over $k$. 
\item[\rm (iii)] $b_1, \ldots, b_n$ are linearly dependent over $k$. 
\end{enumerate}
Then there exist polynomials $A_1, \ldots, A_{n - 1}, B_1, \ldots, B_{n - 1} \in k[T]$ such that the following conditions {\rm (1), (2), (3), (4)} 
are satisfied: 
\begin{enumerate}
\item[\rm (1)] $\ds \sum_{i = 1}^{n - 1} A_i(T) \cdot B_i(T')$ is a symmetric polynomial of $k[T, T']$. 
\item[\rm (2)] $A_1, \ldots, A_{n - 1}$ are linearly independent over $k$. 
\item[\rm (3)] $\Span_k\{b_1, \ldots, b_n\} = \Span_k\{B_1, \ldots, B_{n - 1}\}$. 
\item[\rm (4)] If there exists a symmetric matrix $S \in \Mat(n - 1, k)$ such that 
$(B_1, \ldots, B_{n - 1}) = (A_1, \ldots, A_{n - 1}) S$, then there exists a symmetric matrix $\wt{S} \in \Mat(n, k)$ such that 
$(b_1, \ldots, b_n) = (a_1, \ldots, a_n) \wt{S}$.  
\end{enumerate} 
\end{lem}

\Proof There exists an integer $1 \leq i \leq n$ such that $b_i \in \Span_k \{ b_1, \ldots, b_{i -1}, b_{i + 1}, \ldots, b_n \}$.  
Let $\sigma$ be the permutation of $\{1, 2, \ldots, n\}$ defined by 
\[
\sigma := 
\left\{ 
\begin{array}{ll}
 (i \; \, n) \qquad & \text{ if \quad $i < n$}, \\
 \text{the identity permutation}  \qquad & \text{ if \quad $i = n$} . 
\end{array}
\right.
\] 
Let $a_i' := a_{\sigma(i)}$ and $b_i' := b_{\sigma(i)}$ for all $1 \leq i \leq n$. Then we have the following: 
\begin{enumerate}
\item[\rm (i)$'$] $\ds\sum_{i = 1}^n a_i'(T) \cdot b_i'(T')$ is a symmetric polynomial of $k[T, T']$. 
\item[\rm (ii)$'$] $a_1', \ldots, a_n'$ are linearly independent over $k$. 
\item[\rm (iii)$'$] $b_n' = \ds \sum_{i  = 1}^{n - 1} c_i \cdot b_i' $ for some $c_1, \ldots c_{n - 1} \in k$. 
\end{enumerate} 
Let $A_i := a_i' + c_i \cdot a_n'$ for $1 \leq i \leq n - 1$, and let $B_i := b_i'$ for all $1 \leq i \leq n - 1$. 
We shall show that the $A_i$ and $B_i$ satisfy the conditions (1), (2), (3) and (4). 
By the definitions of $A_i$ and $B_i$ and the above (iii)$'$, we have
\begin{eqnarray*}
\sum_{i = 1}^{n - 1} A_i(T) \cdot B_i(T')
 = \sum_{i = 1}^n a_i'(T) \cdot b_i'(T'). 
\end{eqnarray*}
So, (1) follows from (i)$'$.  Clearly, (2) follows from (ii)$'$, and (3) follows from 
$\Span_k\{b_1, \ldots b_n\} = \Span_k \{ b_1', \ldots, b_{n - 1}'\}$.  
We shall prove (4). 
Write the symmetric matrix $S$ in the assumption of (4) as 
$S= (s_{\alpha, \beta})_{1 \leq \alpha, \beta \leq n - 1}$, where $s_{\alpha, \beta} \in k$ ($1 \leq \alpha, \beta \leq n - 1$). 
Since $B_i = \sum_{j = 1}^{n - 1} s_{j, i} \, A_j$ for all $1 \leq i \leq n - 1$, we have 
\begin{eqnarray*}
 b_i' = \sum_{j = 1}^{n - 1} s_{j, i} a_j' + \sum_{j = 1}^{n - 1} c_j s_{j, i} a_n'  \qquad (1 \leq  i \leq n -1) . 
\end{eqnarray*}
Thus we have 
\begin{eqnarray*}
 b_n' 
 = \sum_{i = 1}^{n - 1} c_i b_i'
 = \sum_{j = 1}^{n - 1} \left( \sum_{i = 1}^{n - 1}  c_i s_{j, i} \right) a_j'  + \left( \sum_{i = 1}^{n - 1} \sum_{j = 1}^{n - 1} c_i c_j s_{j, i} \right) a_n'. 
\end{eqnarray*}
We define $s_{i, j}' \in k$ $(1 \leq \alpha, \beta \leq n)$ as 
\begin{eqnarray*}
s_{\alpha, \beta}' := 
\left\{
\begin{array}{ll}
 s_{\alpha, \beta}  & \text{ if \quad $1 \leq \alpha \leq n - 1$, \; $1 \leq \beta \leq n - 1$}, \\
 \ds \sum_{i = 1}^{n - 1}  c_i s_{\beta, i}  & \text{ if \quad $1 \leq \alpha \leq n - 1$, \; $\beta = n$},  \\
 \ds \sum_{j = 1}^{n - 1} c_j s_{j, \beta} & \text{ if \quad $\alpha = n$, \; $1 \leq \beta \leq n - 1$}, \\
 \ds \sum_{i = 1}^{n -1} \sum_{j = 1}^{n - 1} c_i c_j s_{j, i}  \quad  & \text{ if \quad $\alpha = n$, \; $\beta = n$} . 
\end{array}
\right.
\end{eqnarray*}
Let $S' := (s_{i, j}')_{1 \leq i, j \leq n} \in \Mat(n, k)$. 
We then have $(b_1', \ldots, b_{n - 1}', b_n') = (a_!', \ldots, a_{n - 1}', a_n') S'$. 
Since $S = (s_{i, j})$ is symmetric, $S'$ is also symmetric. 
By Lemma 3.6, we can obtain the desired symmetric matrix $\wt{S}$ appearing in (4). 
\QED

\paragraph{The case where $b_1, \ldots, b_n$ are linearly independent over $k$}

\begin{lem}
Assume that $a_1, \ldots, a_n, b_1, \ldots, b_n \in k[T]$ satisfy the following conditions {\rm (i), (ii), (iii)}: 
\begin{enumerate}
\item[\rm (i)] $\ds\sum_{i  = 1}^n a_i(T) \cdot b_i(T') $ is a symmetric polynomial of $k[T, T']$. 
\item[\rm (ii)] $a_1, \ldots, a_n$ are linearly independent over $k$. 
\item[\rm (iii)] $b_1, \ldots, b_n$ are linearly independent over $k$. 
\end{enumerate}
Then there exists a regular symmetric matrix $S \in \Mat(n, k)$ such that $(b_1, \ldots, b_n) = (a_1, \ldots, a_n) S$. 
\end{lem}

\Proof We proceed by induction on $n$. If $n = 1$, the proof is clear.  
So, let $n \geq 2$. 
We know from Lemma 3.6 that, for obtaining the symmetric matrix $S$ in the conclusion,  
we may change $a_1, \ldots, a_n$ and $b_1, \ldots, b_n$ 
with $a_{\sigma(1)}, \ldots, a_{\sigma(n)}$ and $b_{\sigma(1)}, \ldots, b_{\sigma(n)}$, respectively.  
So, we may assume that $a_1$ or $b_1$ has the minimum order among the polynomials $a_1, \ldots, a_n, b_1, \ldots, b_n$. 
At least one of the following cases can occur: 
\begin{enumerate} 
\item[\rm (a)]  $\ord(a_1) \leq \min\{ \ord(a_1), \, \ldots, \, \ord(a_n), \, \ord(b_1), \, \ldots, \, \ord(b_n) \}$.  
\item[\rm (b)]  $\ord(b_1) \leq \min\{ \ord(a_1), \, \ldots, \, \ord(a_n), \, \ord(b_1), \, \ldots, \, \ord(b_n) \}$.  
\end{enumerate} 
We have only to consider the case (a). 
Let $ A_i$ $(2 \leq i \leq n)$ and $B_i $ $(1 \leq i \leq n)$ be as in assertion (2) of Lemma 3.5. 
We know the following: 
\begin{enumerate} 
\item[\rm (i)$'$] $\ds \sum_{i = 2}^n A_i(T) \cdot B_i(T') = \sum_{i = 2}^n A_i(T') \cdot B_i(T) $. 
\item[\rm (ii)$'$] $A_2, \ldots, A_n$ are linearly independent over $k$. 
\end{enumerate} 
Clearly, (ii)$'$ follows from the assumption (ii). 
We first consider the case where $B_2, \ldots B_n$ are linearly independent over $k$. 
By the induction hypothesis, there exists a symmetric matrix $S \in \Mat(n - 1, k)$ such that 
$(B_2, \ldots ,B_n) = (A_2, \ldots, A_n) S$. 
By Lemma 3.5 (3), we have the desired symmetric matrix.  
Now, we consider the case where $B_2, \ldots B_n$ are linearly dependent over $k$. 
If $n = 2$, we have $B_2 = 0$, which implies $(B_2) =  (A_2)(0)$, and we have the desired symmetric matrix by Lemma 3.5 (3).  
If $n \geq 3$, there exists a non-zero polynomial $B_i$ among $B_2, \ldots , B_n$, by the condition (iii). 
Using Lemma 3.7 in finitely many steps, 
we can reduce polynomials $A_2, \ldots, A_n, B_2, \ldots, B_n$ until we have polynomials 
$\alpha_1, \ldots, \alpha_{n^\flat}, \beta_1, \ldots, \beta_{n^\flat} \in k[T]$ $({n^\flat} < n - 1)$ satisfying the following conditions (i)$''$, (ii)$''$, (iii)$''$, (iv)$''$: 
\begin{enumerate} 
\item[\rm (i)$''$] $\ds \sum_{i = 1}^{n^\flat} \alpha_i(T) \cdot \beta_i (T')$ is a symmetric polynomial of $k[T, T']$. 
\item[\rm (ii)$''$] $\alpha_1, \ldots, \alpha_{n^\flat}$ are linearly independent over $k$. 
\item[\rm (iii)$''$] $\beta_1, \ldots, \beta_{n^\flat}$ are linearly independent over $k$. 
\item[\rm (iv)$''$] If there exists a symmetric matrix $S \in \Mat(n^\flat, k)$ such that 
$(\beta_1, \ldots, \beta_{n^\flat}) = (\alpha_1, \ldots, \alpha_{n^\flat}) S$, 
then there exists a symmetric matrix $S^\sharp \in \Mat(n - 1, k)$ such that 
$(B_2, \ldots, B_{n - 1}) = (A_2, \ldots, A_{n - 1}) S^\sharp$. 
\end{enumerate} 
Note that $\alpha_1, \ldots, \alpha_{n^\flat}, \beta_1, \ldots, \beta_{n^\flat}$ satisfy the induction hypothesis (see 
the above conditions (i)$''$, (ii)$''$, (iii)$''$). 
So, there exists a symmetric matrix $S \in \Mat(n^\flat, k)$ satisfying the assumption of (iv)$''$. 
Hence we have the symmetric matrix $S^\sharp$ appearing in the conclusion of (iv)$''$. 
By Lemma 3.5 (3), we have the desired symmetric matrix. 
\QED

\paragraph{The case where $b_1, \ldots, b_n$ are linearly dependent over $k$}

\begin{lem} 
Assume that $a_1, \ldots, a_n, b_1, \ldots, b_n \in k[T]$ satisfy the following conditions {\rm (i), (ii), (iii)}: 
\begin{enumerate}
\item[\rm (i)] $\ds\sum_{i  = 1}^n a_i(T) \cdot b_i(T') $ is a symmetric polynomial of $k[T, T']$. 
\item[\rm (ii)] $a_1, \ldots, a_n$ are linearly independent over $k$. 
\item[\rm (iii)] $b_1, \ldots, b_n$ are linearly dependent over $k$. 
\end{enumerate}
Then there exists a symmetric matrix $S \in \Mat(n, k)$ such that 
$(b_1, \ldots, b_n) = (a_1, \ldots, a_n) S$. 
\end{lem}

\Proof 
We proceed by induction on $n$. If $n = 1$, the proof is clear. 
So, let $n \geq 2$. 
We apply Lemma 3.7 to the polynomials $a_1, \ldots, a_n, b_1, \ldots, b_n$, and  then we 
have $A_1, \ldots, A_{n - 1}, B_1, \ldots, B_{n - 1}\in k[T]$ as in Lemma 3.7. 
If $B_1, \ldots, B_{n - 1}$ are linearly independent over $k$, Lemma 3.8 and Lemma 3.5 (3) yield the desired symmetric matrix.  
On the ohter hand, if $B_1, \ldots, B_{n - 1}$ are linearly dependent over $k$, 
the induction hypothesis and Lemma 3.5 (3) yield the desired symmetric matrix. 
\QED

\paragraph{A proof of Theorem 3.4}
\quad \medskip

The implication (2) $\Longrightarrow$ (1) follows from a straightforward calculation. 
We shall prove the implication (1) $\Longrightarrow$ (2). 
If $b_1, \ldots, b_n$ are linearly independent over $k$, 
the implication follows from Lemma 3.8. 
On the other hand, if $b_1, \ldots, b_n$ are linearly dependent over $k$, 
the implication follows from Lemma 3.9. 
\QED

\subsubsection{On a symmetric matrix $S$ in characteristic two}

The proof of the following Lemma 3.10 is straightforward. 
Lemma 3.10 shall be used in the proof of Lemma 3.11.

\begin{lem} 
Let $\alpha(T), \beta(T)$ are $p$-polynomials. 
Let $\gamma(T) := \alpha(T) \cdot \beta(T)$. 
Then we have 
\[
 \gamma(T + T') - \gamma(T) - \gamma(T') = \alpha(T) \cdot \beta(T') + \alpha(T') \cdot \beta(T) . 
\]
\end{lem}

\begin{lem}  
Assume $p = 2$. 
Let $a_1, \ldots, a_n \in k[T]$ be $p$-polynomials such that $a_1, \ldots, a_n$ are linearly independent over $k$. 
Let $S = (s_{i, j})_{1 \leq i, j \leq n}$ be a symmetric matrix satisfying 
\[
\bm{a}(T) \cdot S  \cdot {^t}\bm{a}(T') = c(T + T') - c(T) - c(T') 
\]
for some $c(T) \in k[T]$, where $\bm{a}(T) =(a_1, \ldots, a_n) \in k[T]^n$. 
Then all diagonal entries $s_{i, i}$ of $S$ are zeroes. 
\end{lem}

\Proof 
Note that 
\begin{eqnarray*}
\lefteqn{\bm{a}(T) \cdot S  \cdot {^t}\bm{a}(T') } \\
& =& \sum_{i = 1}^r s_{i, i} \cdot a_i(T) \cdot a_i(T' )
 + \sum_{1 \leq i < j \leq n} s_{i, j} \cdot( a_i(T) \cdot a_j(T' ) + a_j(T) \cdot a_i(T') ) . 
\end{eqnarray*}
We know from Lemma 3.10 that there exists a polynomial $c^* \in k[T]$ such that 
\[
(\ast) \qquad  \sum_{i = 1}^n s_{i, i} \cdot a_i(T) \cdot a_i(T') = c^*(T + T') - c^*(T) - c^*(T') . 
\]
We can write the $p$-polynomial $a_1(T)$ appearing in the left hand side of the above equality $(\ast)$ as 
\[
 a_1(T) = \xi_d \cdot T^{p^d} + (\text{ lower order terms in $T$ }) ,
\]
where $\xi_d \in k \backslash \{ 0 \}$. 
We shall prove that $s_{i, i} = 0$ for all $1 \leq i \leq n$ by induction on $n$. 
Let $n = 1$. 
Supposing $s_{1, 1} \ne 0$, we shall find a contradiction. 
The monomial $s_{1, 1} \cdot \xi_d^2 \cdot T^{p^d} \cdot T'^{p^d}$ appears in the left hand side of the equality $(\ast)$, 
but the monomial does not appear in the right hand side of the equality $(\ast)$, because 
no monomial of total degree $2 p^d = p^{d + 1}$, where $\deg(T) = \deg(T') = 1$, appears in the right hand side of the equality $(\ast)$. 
This is a contradiction. So, we have $s_{1, 1} = 0$. 
Now, let $n \geq 2$. We may assume without loss of generality that $a_1$ has the minimum order among 
$a_1, \ldots, a_n$, i.e., $\ord(a_1) \leq \min \{ \ord(a_i) \in \mathbb{Z}_{\geq 0} \mid 1 \leq i \leq n \}$. 
There exist $\lambda_i \in k$ $(2 \leq i \leq n)$ such that 
$\ord (a_i - \lambda_i \cdot a_1) > \ord(a_1)$ for all $2 \leq i \leq n$. 
Let $A_i := a_i - \lambda_i \cdot a_1$ $(2 \leq i \leq n)$. 
Note that 
\begin{eqnarray*}
(\ast\ast) \qquad \sum_{i = 1}^n s_{i, i} \cdot a_i(T) \cdot a_i(T') 
 &=& (s_{1, 1} + \sum_{i = 2}^n s_{i, i}  \cdot \lambda_i^2 ) \cdot a_1(T) \cdot a_1(T') \\
 & & + \sum_{i = 2}^n s_{i, i} \cdot A_i(T) \cdot A_i(T') \\
 & & + \sum_{i = 2}^n s_{i, i} \cdot \lambda_i \cdot (a_1(T') \cdot A_i(T) + a_1(T) \cdot A_i(T')) . 
\end{eqnarray*}
Let 
\[
 \mu :=  s_{1, 1} + \sum_{i = 2}^n s_{i, i}  \cdot \lambda_i^2 \in k.  
\]
Supposing $\mu \ne 0$, we shall find a contradiction. 
The monomial $\mu \cdot \xi_d^2 \cdot T^{p^d} \cdot T'^{p^d}$ appears in the left hand 
side of the above equality $(\ast\ast)$. But no monomial of total degree $p^{d + 1}$ appears 
in  $c^*(T + T') - c^*(T) - c^*(T')$. This is a contradiction. So, $\mu = 0$. Thus, we have 
\begin{eqnarray*}
\begin{array}{l}
\ds \sum_{i = 2}^n s_{i, i} \cdot A_i(T) \cdot A_i(T') + \sum_{i = 2}^n s_{i, i} \cdot \lambda_i \cdot (a_1(T') \cdot A_i(T) + a_1(T) \cdot A_i(T')) \\ [5mm] 
 \qquad =  c^*(T + T') - c^*(T) - c^*(T') . 
\end{array} 
\end{eqnarray*} 
Let 
\[
 c^{**}(T) := c^*(T) - \sum_{i = 2}^n s_{i, i} \cdot \lambda_i \cdot a_1(T) \cdot A_i(T) . 
\]
Then we have 
\[
\sum_{i = 2}^n s_{i, i} \cdot A_i(T) \cdot A_i(T') = c^{**}(T + T') - c^{**}(T) - c^{**}(T') . 
\]
Clearly, all $A_i$ are $p$-polynomials and $A_2, \ldots, A_n$ are linearly independent over $k$. By the induction hypothesis, 
we have  $s_{i, i} = 0$ for all $2 \leq i \leq n$. 
Since $\mu = 0$, we also have $s_{1, 1} = 0$.  
\QED

\subsubsection{On a polynomial matrix of ${^\ell \mf{H}}^{r_1}_{r_2}$ satisfying 
the conditions {\rm (1)} and {\rm (2)} of Theorem 3.3}

Given a polynomial matrix $A(T) = \eta(a_1, \ldots, a_m, a_{m + 1}, \ldots, a_{2m},  a_{2m + 1})$  of ${^\ell \mf{H}}^{r_1}_{r_2}$ 
satisfying the conditions {\rm (1)} and {\rm (2)} of Theorem 3.3, we study, in the following lemma, 
the top right corner $a_{2m + 1}$.  

\begin{lem} 
Let $A(T) = \eta(a_1, \ldots, a_m, a_{m + 1}, \ldots, a_{2m},  a_{2m + 1}) \in {^\ell \mf{H}}^{r_1}_{r_2}$ 
be a polynomial matrix satisfying the conditions {\rm (1)} and {\rm (2)} of {\rm Theorem 3.3}. 
We define a polynomial $\alpha(T) \in k[T]$ as 
\begin{eqnarray*}
\alpha(T) := \left\{ 
\begin{array}{ll}
\ds a_{2m + 1}(T) - \sum_{1 \leq i < j \leq r_1} s_{i, j} \cdot a_i(T) \cdot a_j(T)  & \quad \text{ if \; $p = 2$},  \\ [6mm] 
\ds a_{2m + 1}(T) -  \frac{1}{2}  \cdot \bm{a}(T) \cdot S \cdot {^t} \bm{a}(T) \qquad & \quad  \text{ if  \; $p \geq 3$},
\end{array} 
\right. 
\end{eqnarray*}
where $\bm{a}(T) := (a_1(T), \ldots, a_{r_1}(T)) \in k[T]^{r_1}$.  
Then 
\[
 \alpha(T + T') - \alpha(T) - \alpha(T')
 = a_{2m + 1}(T + T') - a_{2m + 1}(T) - a_{2m + 1}(T')  - \sum_{i = 1}^{r_1} a_i(T) \cdot a_{m + i}(T'). 
\]
\end{lem}

\Proof 
In the case where $p = 2$, we have 
\begin{align*}
\lefteqn{\alpha(T + T') - \alpha(T) - \alpha(T') } \\
  &= \,a_{2m + 1}(T + T') - a_{2m + 1}(T) - a_{2m + 1}(T')  -  \sum_{1 \leq i <  j \leq r_1} s_{i, j}  \cdot ( a_i(T) \cdot a_j(T') + a_j(T) \cdot a_i(T') ) \\
 &=\, a_{2m + 1}(T + T') - a_{2m + 1}(T) - a_{2m + 1}(T')  - \bm{a}(T) \cdot S \cdot {^t} \bm{a}(T) \\
 &=\, a_{2m + 1}(T + T') - a_{2m + 1}(T) - a_{2m + 1}(T')   - \sum_{i = 1}^{r_1} a_i(T) \cdot a_{m + i}(T'). 
\intertext{In the case where $p \geq 3$, we have }
\lefteqn{\alpha(T + T') - \alpha(T) - \alpha(T') } \\
 &=\, a_{2m + 1}(T + T') - a_{2m + 1}(T) - a_{2m + 1}(T')  - \bm{a}(T) \cdot S \cdot {^t} \bm{a}(T)  \\
 &=\, a_{2m + 1}(T + T') - a_{2m + 1}(T) - a_{2m + 1}(T')  - \sum_{i = 1}^{r_1} a_i(T) \cdot a_{m + i}(T'). 
\end{align*}
\QED

\subsubsection{A proof of Theorem 3.3}

Now, we are ready to prove Theorem 3.3. 
We first prove (1) $\Longrightarrow$ (2). So, assume 
$A(T) = \eta(a_1, \ldots, a_m, a_{m + 1}, \ldots, a_{2m},  a_{2m + 1}) \in ({^\ell \mf{H}}^{r_1}_{r_2})^E$. 
We know from Lemma 3.2 that all $a_i$ $(1 \leq i \leq 2m)$ are $p$-polynomials 
and 
\[
(\ast) \qquad a_{2m + 1}(T + T') - a_{2m + 1}(T) - a_{2m + 1}(T')  = \sum_{i = 1}^{r_1} a_i(T) \cdot a_{m + i}(T'). 
\]
By Theorem 3.4, there exists a symmetric matrix $S \in \Mat(r_1, k)$ such that 
\[
 (a_{m + 1}, \ldots, a_{m + r_1}) = (a_1, \ldots, a_{r_1} ) S. 
\]
Let $\bm{a}(T) := (a_1(T), \ldots, a_{r_1}(T)) \in k[T]^{r_1}$.  So, $(\ast)$ implies 
\[
 a_{2m + 1}(T + T') - a_{2m + 1}(T) - a_{2m + 1}(T')  = \bm{a} \cdot S \cdot {^t}\bm{a}. 
\]
If $p = 2$, all diagonal entries of $S$ are zeroes (see Lemma 3.11). 
Now, the conditions (1) and (2) of Theorem 3.3 hold true. 
Let $\alpha(T) \in k[T]$ be the polynomial defined in Lemma 3.12. 
Using $(\ast)$, we then have $\alpha(T + T') - \alpha(T) - \alpha(T) = 0$, which implies 
$\alpha(T)$ is a $p$-polynomial. Thus the condition (3) of Theorem 3.3 holds true. 

We next prove (2) $\Longrightarrow$ (1). 
Since assuming the conditions (1) and (2) of Theorem 3.3 hold true, we can use Lemma 3.12. 
The condition (3) of Theorem 3.3 implies the equality $(\ast)$ holds true. 
Hence, we know from Lemma 3.2 that $A(T)$ is an exponential matrix. 
\QED

\subsection{$(\daleth_{m + 2})^E$}

Given any polynomial matrix of $\daleth_{m + 2}$, 
its equivalent form  can be found in just one of the six sets 
$\daleth_{m + 2}^{\rm I}$, $\daleth_{m + 2}^{\rm II}$, $\daleth_{m + 2}^{\rm III}$, $\daleth_{m + 2}^{\rm IV}$, 
$\daleth_{m + 2}^{\rm V}$, $\daleth_{m + 2}^{\rm VI}$ (see Theorem 2.5). 
Thus we have the following theorem: 

\begin{thm}
The following assertion {\rm (1)} and {\rm (2)} hold true: 
\begin{enumerate}
\item[\rm (1)] $(\daleth_{m + 2})^E \leadsto (\daleth_{m + 2}^{\rm I})^E \cup (\daleth_{m + 2}^{\rm II})^E 
\cup (\daleth_{m + 2}^{\rm III})^E \cup (\daleth_{m + 2}^{\rm IV})^E \cup (\daleth_{m + 2}^{\rm V})^E
\cup (\daleth_{m + 2}^{\rm VI})^E$. 
\item[\rm (2)] The six subsets $(\daleth_{m + 2}^{\rm I})^E, (\daleth_{m + 2}^{\rm II})^E, (\daleth_{m + 2}^{\rm III})^E, 
(\daleth_{m + 2}^{\rm IV})^E, (\daleth_{m + 2}^{\rm V})^E, (\daleth_{m + 2}^{\rm VI})^E$ are mutually $GL(m + 2, k)$-disjoint. 
\end{enumerate}
\end{thm}

By Lemma 1.28, we can describe the six sets $(\daleth_{m + 2}^{\rm I})^E$
$(\daleth_{m + 2}^{\rm II})^E$, $(\daleth_{m + 2}^{\rm III})^E$, $(\daleth_{m + 2}^{\rm IV})^E$, 
$(\daleth_{m + 2}^{\rm V})^E$, $(\daleth_{m + 2}^{\rm VI})^E$ as follows:

\begin{thm} 
The following assertions {\rm (1), (2), (3), (4), (5), (6)} hold true: 
\begin{enumerate}
\item[\rm (1)] $(\daleth_{m + 2}^{\rm I})^E = \daleth_{m + 2}^{\rm I} = \{ I_{m + 2} \}$.

\item[\rm (2)] Let $A(T) = \Lambda(1, j; \alpha(T)) \in \daleth_{m + 2}^{\rm II}$, where $1 \leq j \leq m + 1$. 
Then $A(T) \in (\daleth_{m + 2}^{\rm II})^E$ if and only if all entries of $\alpha(T)$ are $p$-polynomials.

\item[\rm (3)] Let $A(T) = \Lambda(i, 1; \alpha(T)) \in \daleth_{m + 2}^{\rm III}$, where $1 \leq i \leq m + 1$. 
Then $A(T) \in (\daleth_{m + 2}^{\rm III})^E$ if and only if all entries of $\alpha(T)$ are $p$-polynomials.

\item[\rm (4)] Let $A(T) = \Lambda(2, 2; \alpha(T)) \in \daleth_{m + 2}^{\rm IV} ( = ({^1}\daleth_1)^{\dim = 2})^0 )$, where   
$\alpha(T) = L(a_1 \mid a_2 \mid 0) \in {^1}\mf{L}_1^0$. 
Then $A(T) \in (\daleth_{m + 2}^{\rm IV})^E$ if and only if $a_1, a_2$ are $p$-polynomials.

\item[\rm (5)] 
Let $A(T) = \Lambda(2, 2; \alpha(T)) \in  \daleth_{m + 2}^{\rm V}  ( = ({^1}\daleth_1)^{\dim = 3} )$, where  
$\alpha(T) = L(a_1 \mid a_2 \mid a) \in ({^1}\mf{L}_1)^{\dim = 3}$.  
Then $A(T) \in (\daleth_{m + 2}^{\rm V})^E$ if and only if $a_1, a_2, a$ are $p$-polynomials.

\item[\rm (6)] 
Let $A(T) = \Lambda(j + 1, i + 1; \alpha(T) ) \in \daleth_{m + 2}^{\rm VI}$, where  
$\alpha(T) = L(a_1, \ldots, a_i \mid a_{i + 1}, \ldots, a_{i + j} \mid a) \in {^i}\mf{L}_j$ and $3 \leq i  + j \leq m$. 
Then $A(T) \in (\daleth_{m + 2}^{\rm VI})^E$ if and only if $a_1, \ldots, a_i, a_{i + 1}, \ldots, a_{i + j}, a$ are $p$-polynomials.
\end{enumerate} 
\end{thm}

\subsection{Equivalence relations of exponential matrices of Heisenberg groups}

In this subsection, we study equivalence of two exponential matrices of $({^\ell}\mf{H}^{r_1}_{r_2})^E$ (see Theorems 3.17 and 3.23 
in Subsubsection 3.4.1) and those of $(\daleth_{m + 2})^E$ (see Theorem 3.24 in Subsubsection 3.4.2).

\subsubsection{$({^\ell}\mf{H}^{r_1}_{r_2})^E$}

For considering equivalence of two exponential matrices of $({^\ell}\mf{H}^{r_1}_{r_2})^E$, 
we shall define a parametrization of  $({^\ell}\mf{H}^{r_1}_{r_2})^E$. 
For it, we prepare several notations. 

Let $R$ be a commutative ring $R$ with unity and let $n \geq 1$ be an integer. 
We denote by $GLS(n, R)$ the set of all symmetric matrices belonging to $GL(n, R)$, 
and denote by $GLS(n, R)^0$ the set of all matrices of $GLS(n, R)$ whose all diagonal entries are zeroes.  
We define a set $GL\ss(n, R)$ as 
\[
GL\ss(n, R) := 
\left\{
\begin{array}{ll}
GLS(n, R)^0 \quad & \text{ if \quad $p = 2$}, \\
GLS(n, R)              & \text{ if \quad $p \geq 3$} . 
\end{array}
\right.
\]

Let $R'$ be a non-necessarily commutative $k$-algebra. 
For integers $i \geq 1$, $i' \geq 0$, and $i'' \geq 1$, we denote by $\F(i, i', i''; R')$ denote 
the set of all elements $(f_1, \ldots, f_i, f_{i + 1}, \ldots, f_{i + i'}, g_{i + 1}, \ldots, g_{i + i''})$ 
of $R'^{i+ i' + i''}$ such that $i + i'$ elements $f_1, \ldots, f_{i + i'}$ of $R'$ are linearly independent over $k$, 
and $i + i''$ elements $f_1, \ldots, f_i, g_{i + 1}, \ldots, g_{i + i''}$ of $R'$ are also linearly independent over $k$.

For any $(\ell, r_1, r_2) \in \Omega_m$, we define sets ${^\ell}\X^{r_1}_{r_2}$ and ${^\ell}\Y^{r_1}_{r_2}$ as 
\begin{eqnarray*}
\left\{
\begin{array}{rcl}
{^\ell}\X^{r_1}_{r_2} &:=& GL\ss(r_1, k) \times \F(r_1, \ell - r_1, r_2; \mf{P}) \times \mf{P} , \\
{^\ell}\Y^{r_1}_{r_2} &:=& GL\ss(r_1, k) \times k[T]^{r_1} \times k[T]^{\ell - r_1} \times k[T]^{r_2} \times k[T] , 
\end{array}
\right. 
\end{eqnarray*}
where $\mf{P}$ is the ring consisting of all $p$-polynomials. 
Clearly, ${^\ell}\X^{r_1}_{r_2} \subset {^\ell}\Y^{r_1}_{r_2}$. 
We can define a map $\vartheta : {^\ell}\Y^{r_1}_{r_2} \to {^\ell}\W^{r_1}_{r_2}$ as 
\[
\vartheta(S, \bm{a}_1, \bm{a}_2, \bm{\alpha}_2, \alpha) := (\bm{a}_1, \bm{a}_2, \bm{\alpha}_1, \bm{\alpha}_2, \varphi) , 
\]
where 
\begin{gather*}
\bm{\alpha}_1 : = \bm{a}_1 \, S, \\
\intertext{and by writing $\bm{a}_1 = (a_1, \ldots, a_{r_1})$, we let  }
\varphi := 
\left\{
\begin{array}{ll}
 \alpha + \ds\sum_{1 \leq i < j \leq r_1} s_{i, j} \cdot  a_i \cdot a_j  \quad & \text{ if \quad $p = 2$}, \\ [6mm] 
 \alpha + \ds \frac{1}{2} \cdot \bm{a}_1 \cdot S \cdot {^t} \bm{a}_1  & \text{ if \quad $p \geq 3$} . 
\end{array}
\right. 
\end{gather*} 
We define a map $\hslash : {^\ell}\Y^{r_1}_{r_2} \to \mf{H}(m + 2, k[T])$ as $\hslash = h \circ \vartheta$.

\begin{lem}
For any $(S, \bm{a}_1, \bm{a}_2, \bm{\alpha}_2, \alpha) \in {^\ell}\Y^{r_1}_{r_2}$, 
the following conditions {\rm (1)} and {\rm (2)} are equivalent: 
\begin{enumerate} 
\item[\rm (1)] $(S, \bm{a}_1, \bm{a}_2, \bm{\alpha}_2, \alpha) \in {^\ell}\X^{r_1}_{r_2}$. 
\item[\rm (2)] $\hslash (S, \bm{a}_1, \bm{a}_2, \bm{\alpha}_2, \alpha) \in ({^\ell}\mf{H}^{r_1}_{r_2})^E$. 
\end{enumerate}
\end{lem}

\Proof The proof follows from the definition of ${^\ell}\X^{r_1}_{r_2}$ and Theorem 3.3. 
\QED 

By the above Lemma 3.15, we can regard $\hslash$ as a bijection from ${^\ell}\X^{r_1}_{r_2}$ to $({^\ell}\mf{H}^{r_1}_{r_2})^E$. 
We have the following commutative diagram, where upward arrows are inclusions: 
\begin{eqnarray*}
\xymatrix@!C=36pt{                                     &                                                 & \mf{H}(m + 2, k[T]) \\
 {^\ell}\Y^{r_1}_{r_2} \ar[r]_\vartheta  \ar[rru]^{\hslash} & {^\ell}\W^{r_1}_{r_2} \ar[r]_(.42){\cong} \ar[ru]_(0.52){h}  
 &  h( {^\ell}\W^{r_1}_{r_2}) \ar[u] \\
 {^\ell}\X^{r_1}_{r_2} \ar[r]^{\hookrightarrow} \ar[u] \ar[rrd]_{\hslash} &  {^\ell}\V^{r_1}_{r_2} \ar[r]^{\cong} \ar[u]     &  {^\ell}\mf{H}^{r_1}_{r_2} \ar[u] \\ 
 & & ( {^\ell} \mf{H}^{r_1}_{r_2} \ar[u] )^E
}
\end{eqnarray*}
Keeping the above commutative diagram in mind, 
we shall define a group ${^\ell}G^{r_1}_{r_2}$ and a map ${^\ell}\Y^{r_1}_{r_2} \times {^\ell}G^{r_1}_{r_2} \to {^\ell}\Y^{r_1}_{r_2}$, 
which make it possible to consider an equivalence relation of two exponential matrices $A(T)$ and $B(T)$ of 
$({^\ell}\mf{H}^{r_1}_{r_2})^E$ in terms of two corresponding elements $x$ and $x'$ of ${^\ell}\X^{r_1}_{r_2}$ (see Theorem 3.17).  
We mention here that the map  ${^\ell}\Y^{r_1}_{r_2} \times {^\ell}G^{r_1}_{r_2} \to {^\ell}\Y^{r_1}_{r_2}$ becomes an action 
(see Lemma 3.18).

Let ${^\ell}G^{r_1}_{r_2}$ be the set of all regular matrices $g$ of $GL(\ell + r_1 + r_2 + 1, k)$ with the form 
\[
 g = 
\left(
\begin{array}{c | c | c | c | c}
 g_{1, 1} & O          & O         & O          & O \\ 
\hline 
 O        & g_{2, 2}  & g_{2, 3}   & g_{2, 4}    & g_{2, 5} \\ 
\hline 
 O        & O         & g_{3, 3}   & O          & g_{3, 5} \\
\hline 
 O        & O         & O         & g_{4, 4}    & g_{4, 5} \\ 
\hline 
 O        & O         & O         & O           & g_{5, 5} 
\end{array}
\right), 
\]
where $g_{1, 1} \cdot {^t}g_{2, 2} = g_{5, 5}^{\oplus r_1}$ and the sizes of the square matrices of $g_{1, 1}$, $g_{2, 2}$, $g_{3, 3}$, $g_{4, 4}$, $g_{5, 5}$ are 
$r_1$, $r_1$, $\ell - r_1$, $r_2$, $1$, respectively. 

For any matrix $\Sigma = (\sigma_{i, j}) \in \Mat(n, k)$, 
we define an upper triangular matrix $\Sigma_+ = (\sigma_{i, j}^+)$ and an lower triangular matrix $\Sigma_- = (\sigma_{i, j}^{-})$ as 
\[
 \sigma_{i, j}^+ = 
\left\{
\begin{array}{ll}
 \sigma_{i, j} \quad  & \text{ if \quad  $i < j$}, \\
 0                & \text{ if \quad $i \geq j$} 
\end{array} 
\right. , 
\qquad 
 \sigma_{i, j}^- = 
\left\{
\begin{array}{ll}
 0  & \text{ if \quad  $i \leq j$}, \\
 \sigma_{i, j} \quad                 & \text{ if \quad $i > j$} 
\end{array} 
\right. . 
\]
Clearly, if $\Sigma$ is a symmetric matrix, we have ${^t}\Sigma_+ = \Sigma_-$.

For $n \geq 1$, $\bm{a} \in k[T]^n$, and $Q, \Sigma \in \Mat(n, k)$, 
we define a polynomial $\Delta_{\bm{a}, Q, \Sigma_+}$ of $k[T]$ as  
\begin{eqnarray*}
 \Delta_{\bm{a}, Q, \Sigma_+} := 
\left\{
\begin{array}{ll}
 \Tr ( \diag(\bm{a}) \cdot Q \cdot \Sigma_+ \cdot {^t}Q \cdot \diag( \bm{a} ) ) \qquad  & \text{ if \qquad $p  = 2$} , 
 \\
 0 & \text{ if \qquad $p \geq 3$}. 
\end{array}
\right. 
\end{eqnarray*}
Clearly, if $\bm{a} \in \mf{P}^n$, then $\Delta_{\bm{a}, Q, \Sigma_+} \in \mf{P}$. 
For any $c \in k$, we have $c \cdot  \Delta_{\bm{a}, Q, \Sigma_+} =  \Delta_{\bm{a}, Q, c \cdot\Sigma_+}$.  

We can define a map ${^\ell}\Y^{r_1}_{r_2} \times {^\ell}G^{r_1}_{r_2} \to {^\ell}\Y^{r_1}_{r_2}$ as 
\begin{eqnarray*}
\lefteqn{(S,  \bm{a}_1,  \bm{a}_2,   \bm{\alpha}_2, \alpha) \star g } \\
& : = & 
\left( 
\begin{array}{r}
S', 
\; \, 
(\bm{a}_1, \bm{a}_2, \bm{\alpha}_2, \alpha) 
\left(
\begin{array}{c | c | c |c}
 g_{2, 2} & g_{2, 3} & g_{2, 4} & g_{2, 5} \\ 
\hline 
 O & g_{3, 3} & O & g_{3, 5} \\ 
\hline 
 O & O & g_{4, 4} & g_{4, 5} \\ 
\hline 
 O & O & O & g_{5, 5}
\end{array}
\right)
- 
(0, 0, 0,  \Delta_{\bm{a}_1, \; g_{2, 2}, \; S'_+}  ) 
\end{array}
\right) , 
\end{eqnarray*} 
where 
\[
 S' := g_{2, 2}^{-1} \cdot S \cdot g_{1, 1} .
\]
Note that $S' \in GL\ss (r_1, k)$. In fact, $S' = g_{5, 5}^{-1} \cdot {^t} g_{1, 1} \cdot S \cdot g_{1, 1} \in GLS(r_1, k)$.  
So, assume $p = 2$. Write $g_{1, 1} = (\gamma_{i, j})$ and $S = (s_{i, j})$. 
Then the $(i, i)$-th entry of ${^t} g_{1, 1} \cdot S \cdot g_{1, 1}$ is calculated as 
\begin{eqnarray*}
 \sum_{j, j' = 1}^{r_1} \gamma_{j, i} \cdot s_{j, j'} \cdot \gamma_{j', i}
 &=& \sum_{1 \leq j < j' \leq r_1} \gamma_{j, i} \cdot s_{j, j'} \cdot \gamma_{j', i}
    + \sum_{1 \leq j \leq r_1} \gamma_{j, i} \cdot s_{j, j} \cdot \gamma_{j, i}
    +  \sum_{1 \leq j' < j \leq r_1} \gamma_{j, i} \cdot s_{j, j'} \cdot \gamma_{j', i} \\
 &=& 2 \cdot  \sum_{1 \leq j < j' \leq r_1} \gamma_{j, i} \cdot s_{j, j'} \cdot \gamma_{j', i}
    + \sum_{1 \leq j \leq r_1} \gamma_{j, i} \cdot 0 \cdot \gamma_{j, i} = 0 . 
\end{eqnarray*}

The following basic lemma shall be used in the proof of the following Theorem 3.17. 

\begin{lem}
Let $S \in \Mat(n, k)$ be a symmetric matrix. 
Let $S_1, S_2$ be matrices of $\Mat(n, k)$ satisfying $S = S_1 + S_2$ and ${^t}S_1 = S_2$. 
Write $S = (s_{i, j})$ and $S_1 = (s_{i, j}^1 )$. 
Then for all $x_1, \ldots, x_n \in k[T]$, we have 
\begin{eqnarray*}
 \sum_{1 \leq i < j \leq n} s_{i, j} x_i x_j = \sum_{1 \leq i, j \leq n} s_{i, j}^1 x_i x_j - \sum_{i = 1}^n s_{i, i}^1 x_i^2 
 = \bm{x} \cdot S_1 \cdot {^t} \bm{x} - \Tr( \diag( \bm{x} ) \cdot S_1 \cdot \diag(\bm{x})) , 
\end{eqnarray*}
where $\bm{x} := (x_1, \ldots, x_n)$. 
\end{lem}

\paragraph{Equivalence relations of exponential matrices of $({^\ell}\mf{H}^{r_1}_{r_2})^E$}

\begin{thm} 
Let $A(T)  \in ( {^\ell} \mf{H}^{r_1}_{r_2} )^E$ and 
let $B(T) \in  ({^\ell} \mf{H}^{r_1}_{r_2})^E$. 
Write $A(T) = \hslash(x)$ for some  $x = (S, \bm{a}_1, \bm{a}_2, \bm{\alpha}_2, \alpha) \in {^\ell}\X^{r_1}_{r_2}$, 
and $B(T) = \hslash(x')$ for some $x' = (S', \bm{b}_1, \bm{b}_2, \bm{\beta}_2, \beta) \in {^\ell}\X^{r_1}_{r_2}$. 
Then the following conditions {\rm (1)} and {\rm (2)} are equivalent: 
\begin{enumerate}
\item[\rm (1)] $A(T)$ and $B(T)$ are equivalent. 
\item[\rm (2)] $x' = x \star g$ 
for some $g \in {^\ell}G^{r_1}_{r_2}$.  
\end{enumerate}
\end{thm}

\Proof We first prove (1) $\Longrightarrow$ (2). 
Using the map $\vartheta : {^\ell}\Y^{r_1}_{r_2} \to {^\ell}\W^{r_1}_{r_2}$, we let
\begin{eqnarray*}
\vartheta(x) = (\bm{a}_1, \bm{a}_2, \bm{\alpha}_1, \bm{\alpha}_2, \varphi)
\qquad \text{ and } \qquad 
\vartheta(x') = (\bm{b}_1, \bm{b}_2, \bm{\beta}_1, \bm{\beta}_2, \psi) . 
\end{eqnarray*}
Write $S = (s_{i, j})$, $\bm{a}_1 = (a_1, \ldots, a_{r_1})$, $S' = (s_{i, j}')$, $\bm{b}_1 = (b_1, \ldots, b_{r_1})$. 
So, 
\begin{align}
 \bm{\alpha}_1 = \bm{a}_1 \, S, 
  \qquad     \varphi = \ds \alpha + \sum_{1 \leq i < j \leq r_1} s_{i, j} a_i a_j, 
  \qquad    \bm{b}_1  = \bm{\beta}_1 \, S',
  \qquad     \psi = \ds \beta + \sum_{1 \leq i < j \leq r_1} s_{i, j}' b_i b_j .  
  \tag{i} 
\end{align}
We know from Lemma 2.12 that 
$(\bm{a}_1, \bm{a}_2, \bm{\alpha}_1, \bm{\alpha}_2, \varphi) = (\bm{b}_1, \bm{b}_2, \bm{\beta}_1, \bm{\beta}_2, \psi) \, Q$ 
for some $Q \in {^\ell}\cQ^{r_1}_{r_2}$. 
The equality implies 
\begin{align}
 \bm{b}_1  &=  \bm{a}_1\,  Q_{1, 1},                                    \tag{ii}\\
 \bm{b}_2 &= \bm{a}_1 \, Q_{1, 2} + \bm{a}_2 \, Q_{2, 2},         \tag{iii} \\
 \bm{\beta}_1 &= \bm{\alpha}_1 \, Q_{3, 3},                         \tag{iv}\\
 \bm{\beta}_2 &= \bm{\alpha}_1 \, Q_{3, 4} + \bm{\alpha}_2 \, Q_{4, 4},                    \tag{v}\\
 \psi &= \bm{a}_1 \, Q_{1, 5} + \bm{a}_2 \, Q_{2, 5} + \bm{\alpha}_1 \, Q_{3, 5} + \bm{\alpha}_2 \, Q_{4, 5} + \varphi \, Q_{5, 5}  \tag{vi}. 
\intertext{Thus we respectively have from (ii), (iv), (i) and from (v), (i) }
  S' &= Q_{1, 1}^{-1} \, S \, Q_{3, 3},                                                           \tag{vii}\\
  \bm{\beta}_2 &= \bm{a}_1 \, S \, Q_{3, 4} + \bm{\alpha}_2 \, Q_{4, 4} .          \tag{viii}
\intertext{We shall show that $\beta$ has the following expression: }
\beta &= 
\bm{a}_1 \, (Q_{1, 5} + S \, Q_{3, 5} ) + \bm{a}_2 \, Q_{2, 5} + \bm{\alpha}_2 \, Q_{4, 5} + \alpha \, Q_{5, 5}  - \Delta_{\bm{a}_1, Q_{1, 1}, S'_+}  .       
                     \tag{ix}
\end{align}
In the case where $p = 2$, we have 
\begin{eqnarray*}
 \beta
 &\overset{\rm (i)}{ = }& \psi -  \sum_{1 \leq i < j \leq r_1} s_{i, j}' b_i b_j \\
 & \overset{\rm (vi), \; (i)}{ = } & \bm{a}_1 \,  ( Q_{1, 5} + S \, Q_{3, 5} ) + \bm{a}_2 \, Q_{2, 5}  + \bm{\alpha}_2 \, Q_{4, 5} 
  + \left( \alpha + \sum_{1 \leq i < j \leq r_1} s_{i, j} a_i a_j \right) Q_{5, 5} 
 - \sum_{1 \leq i < j \leq r_1} s_{i, j}' b_i b_j .   
\end{eqnarray*}
By the definition of ${^\ell}\cQ^{r_1}_{r_2}$ and (vii), we have $S \, Q_{5, 5}^{\oplus r_1} = S \, Q_{3, 3} \, {^t} Q_{1, 1} =  Q_{1, 1} \, S' \, {^t}Q_{1, 1}$, which implies 
$q_{5, 5} \, S = Q_{1, 1} \, S'_+ \, {^t}Q_{1, 1} + Q_{1, 1} \, S'_- \, {^t}Q_{1, 1}$, where $Q_{5, 5} = (q_{5, 5})$. 
Since ${^t} ( Q_{1, 1} \, S'_+ \, {^t}Q_{1, 1} )  = Q_{1, 1} \, S'_- \, {^t}Q_{1, 1}$, we know from Lemma 3.16 that 
\begin{eqnarray*}
q_{5, 5} \cdot \sum_{1 \leq i < j \leq r_1} s_{i, j} a_i a_j 
&=& \bm{a}_1 \cdot Q_{1, 1} \, S'_+ \, {^t} Q_{1, 1} \cdot {^t}\bm{a}_1 - \Tr (\diag(\bm{a}_1) \cdot Q_{1, 1} \, S'_+ \, {^t}Q_{1, 1} \cdot \diag(\bm{a}_1) ) \\
 &\overset{\rm (ii)}{=}& \sum_{1 \leq i < j \leq r_1} s_{i, j}' b_i b_j - \Delta_{\bm{a}_1, Q_{1, 1}, S'_+} . 
\end{eqnarray*}
Thus we have the desired expression (ix) of $\beta$. 
In the case where $p \geq 3$, (ix) follows from the following two expressions of $\psi$: 
\begin{eqnarray*}
\left\{
\begin{array}{rcl}
\psi &\overset{\rm (i), \;  (ii)}{=} & \ds \beta + \frac{1}{2} \, \bm{a}_1 \, Q_{1, 1} \, S' \, {^t}Q_{1, 1} \, {^t}\bm{a}_1 
 \overset{\rm (vii)}{=} \beta + \frac{1}{2} \, \bm{a}_1 \, S \, Q_{3, 3} \, {^t}Q_{1, 1} \, {^t}\bm{a}_1
  \overset{Q \in {^\ell}\cQ^{r_1}_{r_2}, \; {\rm (i)}}{=} \beta + (\varphi - \alpha) \, Q_{5, 5}, \\
\psi
 &\overset{\rm (vi), \; (i)}{=} & \bm{a}_1\, (Q_{1, 5} + S \, Q_{3, 5}) + \bm{a}_2 \, Q_{2, 5} + \bm{\alpha}_2\, Q_{4, 5} + \varphi \, Q_{5, 5} . 
\end{array} 
\right. 
\end{eqnarray*}

Now, letting $g$ be the matrix of $GL(\ell + r_1 + r_2 + 1, k)$ defined by 
\[
 g := 
\left(
\begin{array}{c | c | c | c | c} 
 Q_{3, 3} & O & O & O & O \\
\hline 
 O & Q_{1, 1} & Q_{1, 2} & S \, Q_{3, 4} & Q_{1, 5} + S \, Q_{3, 5} \\ 
\hline 
 O & O & Q_{2, 2} & O & Q_{2, 5} \\ 
\hline 
 O & O & O & Q_{4, 4} & Q_{4, 5} \\ 
\hline 
 O & O & O & O & Q_{5, 5}
\end{array}
\right) .
\]
we have $g \in {^\ell}G^{r_1}_{r_2}$ and 
$(S', \bm{b}_1, \bm{b}_2, \bm{\beta}_2, \beta) = (S, \bm{a}_1, \bm{a}_2, \bm{\alpha}_2, \alpha) \star g$.

We next prove (2) $\Longrightarrow$ (1). We let  
\[
 Q := 
\left(
\begin{array}{cc | cc | c}
 g_{2, 2} & g_{2, 3} & O & O& O \\
 O & g_{3, 3} & O & O & g_{3, 5} \\
\hline 
 O & O & g_{1, 1} & S^{-1} \, g_{2, 4} & S^{-1} \, g_{2, 5} \\
 O & O & O & g_{4, 4} & g_{4, 5} \\
\hline 
 O & O & O & O & g_{5, 5} 
\end{array}
\right) . 
\]
Then $Q \in {^\ell}\cQ^{r_1}_{r_2}$ and 
$(\bm{b}_1, \bm{b}_2, \bm{\beta}_1, \bm{\beta}_2, \psi) =  (\bm{a}_1, \bm{a}_2, \bm{\alpha}_1, \bm{\alpha}_2, \varphi) \, Q$, 
which implies $A(T)$ and $B(T)$ are equivalent (see Lemma 2.12). 
\QED

\paragraph{${^\ell}\Y^{r_1}_{r_2} \curvearrowleft {^\ell}G^{r_1}_{r_2}$ and ${^\ell}\X^{r_1}_{r_2} \curvearrowleft {^\ell}G^{r_1}_{r_2}$}

\begin{lem}
The following assertions {\rm (1)} and {\rm (2)} hold true: 
\begin{enumerate}
\item[\rm (1)] ${^\ell}G^{r_1}_{r_2}$ acts on ${^\ell}\Y^{r_1}_{r_2}$ from right. 
\item[\rm (2)] ${^\ell}\X^{r_1}_{r_2}$ is ${^\ell}G^{r_1}_{r_2}$-invariant, i.e., 
for all $x \in {^\ell}\X^{r_1}_{r_2}$ and $g \in {^\ell}G^{r_1}_{r_2}$, 
we have $x \star g \in {^\ell}\X^{r_1}_{r_2}$ 
\end{enumerate} 
\end{lem}

The above Lemma 3.18 and Theorem 3.17 enable us 
to describe the set of all equivalence classes of exponential matrices of 
$({^\ell}\mf{H}^{r_1}_{r_2})^E$ as a quotient ${^\ell}\X^{r_1}_{r_2}/{^\ell}G^{r_1}_{r_2}$ of the action 
${^\ell}\X^{r_1}_{r_2} \curvearrowleft {^\ell}G^{r_1}_{r_2}$. 
For proving the above Lemma 3.18, we prepare the following Lemmas 3.19, 3.20, 3.21, and 3.22.

Let $n \geq 1$. 
We denote by $\Seq$ the totally ordered set $\{ (\lambda, \mu) \in \Z^2 \mid 1 \leq \lambda < \mu \leq n \}$ whose ordering $\prec$ is defined by 
\begin{eqnarray*}
 (\lambda, \mu) \prec (\lambda', \mu') 
\Longleftrightarrow 
\left\{
\begin{array}{l}
\mu < \mu' , \\
\mu = \mu' \quad \text{ and } \quad  \lambda < \lambda' . 
\end{array} 
\right.  
\end{eqnarray*}
So, $\Seq = \{ (1, 2), (1, 3), (2, 3), (1,4),\ldots, (n - 1, n) \}$. 
Clearly, the number of elements of $\Seq$ is equal to $N := (n^2 - n)/2$.  

For any matrix $Q = (q_{i, j}) \in \Mat(n, k)$, any integer $1 \leq i \leq n$, and any $(\lambda, \mu) \in \Seq$, 
we let $i \circ (\lambda, \mu) := q_{i, \lambda} q_{i, \mu}$, and let $Q^{[1, n] \circ \Seq}$ be the matrix of $\Mat_{n, N}(k)$
whose $(i, (\lambda, \mu))$-th entries are defined by $i \circ (\lambda, \mu)$. So the $i$-th row $\bm{r}_i$ of 
$Q^{[1, n] \circ \Seq}$ has the following expression:  
\[
\bm{r}_i = ( i \circ (\lambda, \mu) )_{(\lambda, \mu ) \in \Seq}
  = (q_{i, 1} q_{i, 2}, \; q_{i, 1}q_{1, 3}, \; q_{i, 2} q_{i, 3}, \; \ldots \; ,  \; q_{i, n - 1} q_{i, n}) . 
\]

\begin{lem} 
For all $\bm{a} = (a_1, \ldots, a_n) \in k[T]^n$, and $Q, \Sigma = (\sigma_{i, j}) \in \Mat(n, k)$, 
we have 
\[
 \Delta_{\bm{a}, Q, \Sigma_+} = (a_1^2, \ldots, a_n^2 ) \cdot Q^{[1, n] \circ \Seq} \cdot {^t}(\sigma_{\lambda, \mu})_{(\lambda, \mu) \in \Seq} . 
\]
\end{lem}

\Proof  By the definition of $\Delta_{\bm{a}, Q, \Sigma_+}$, calculating 
the $(i, i)$-th entries of $Q \cdot \Sigma_+ \cdot {^t} Q$ $(1 \leq i \leq n)$, 
we have 
\begin{eqnarray*}
\Delta_{\bm{a}, Q, \Sigma_+} = \sum_{i = 1}^n a_i^2 \cdot 
 \sum_{1 \leq \lambda, \mu \leq n} q_{i, \lambda} \cdot \sigma_{\lambda, \mu} \cdot q_{i, \mu} 
= \sum_{i = 1}^n  a_i^2 \cdot \sum_{(\lambda, \mu) \in \Seq} q_{i, \lambda}  q_{i, \mu} \cdot \sigma_{\lambda, \mu} . 
\end{eqnarray*}
\QED

For any $R = (r_{i, j}) \in \Mat(n, k)$ and for $(i, j), (\lambda, \mu) \in \Seq$, 
we let $(i, j) \Join (\lambda, \mu) := r_{i, \lambda} r_{j, \mu} + r_{i, \mu} r_{j, \lambda}$ and let 
$R^{\Seq \Join \Seq}$ be the matrix of $\Mat(N, k)$ whose $((i, j), (\lambda, \mu) )$-th entries are defined by $(i, j) \Join (\lambda, \mu)$.

\begin{lem}
Let $\Sigma = (\sigma_{i, j}) \in \Mat(n, k)$ be a symmetric matrix, and 
let $\Theta = (\tau_{i, j}) \in \Mat(n, k)$ be a symmetric matrix whose all diagonal entries are zeroes. 
Assume that $R = (r_{i, j}) \in \Mat(n, k)$ satisfies $\Sigma = R \cdot \Theta \cdot {^t} R$. 
Then we have 
\[
^{t}(\sigma_{\lambda, \mu})_{(\lambda, \mu) \in \Seq} = R^{\Seq \Join \Seq} \cdot {^t} (\tau_{\lambda, \mu})_{(\lambda, \mu) \in \Seq}. 
\] 
\end{lem}

\Proof For any $(i, j) \in \Seq$, we have 
\begin{eqnarray*}
\sigma_{i, j}
 = \sum_{1 \leq \lambda, \mu \leq n} r_{i, \lambda} \cdot \tau_{\lambda, \mu} \cdot r_{j, \mu} 
 = \sum_{1 \leq \lambda < \mu \leq n} ( r_{i, \lambda} r_{j, \mu} + r_{i, \mu} r_{j, \lambda} ) \cdot \tau_{\lambda, \mu} 
 = \sum_{(\lambda, \mu) \in \Seq} (i, j) \Join (\lambda, \mu) \cdot \tau_{\lambda, \mu} .
\end{eqnarray*}
\QED

For any $Q = (q_{i, j}) \in \Mat(n, k)$, we let $Q^{\square}$ be the matrix of $\Mat(n, k)$ whose $(i, j)$-th entries are defined by 
$q_{i, j}^2$.

\begin{lem}
For all $Q, R \in \Mat(n, k)$, we have 
\[
(Q \cdot R)^{[1, n] \circ \Seq} = Q^{[1, n] \circ \Seq} \cdot R^{\Seq \Join \Seq} + Q^{\square} \cdot R^{[1, n] \circ \Seq} . 
\]
\end{lem}

\Proof Write $Q = (q_{i, j})$, $R = (r_{i, j})$, and $Q \cdot R = (c_{i, j})$. 
We can calculate the $(\xi, (\lambda, \mu))$-th entries of $(Q \cdot R)^{[1, n] \circ \Seq}$ as 
\setcounter{equation}{0}
\begin{eqnarray}
\xi \circ (\lambda, \mu)
 &=& c_{\xi, \lambda} c_{\xi, \mu}
  = \left( \sum_{i = 1}^n  q_{\xi, i} r_{i, \lambda} \right) \cdot \left( \sum_{j = 1}^n q_{\xi, j} r_{j, \mu}\right)  \nonumber \\ 
 &=& \sum_{1 \leq i, j \leq n} q_{\xi, i} q_{\xi, j} r_{i, \lambda} r_{j, \mu}   \nonumber \\
 &=& \sum_{1 \leq i < j \leq n} q_{\xi, i} q_{\xi, j} r_{i, \lambda} r_{j, \mu} 
        + \sum_{1 \leq i = j \leq n} q_{\xi, i} q_{\xi, j} r_{i, \lambda} r_{j, \mu} 
        + \sum_{1 \leq j < i \leq n} q_{\xi, i} q_{\xi, j} r_{i, \lambda} r_{j, \mu}   \nonumber \\
 &=& \sum_{(i, j) \in \Seq}  q_{\xi, i} q_{\xi, j} r_{i, \lambda} r_{j, \mu} 
        + \sum_{i = 1}^n q_{\xi, i}^2 \cdot r_{i, \lambda} r_{i, \mu} 
        + \sum_{(j, i) \in \Seq}q_{\xi, j}  q_{\xi, i} r_{j, \mu} r_{i, \lambda}  .  
\end{eqnarray}
We can calculate the $(\xi, (\lambda, \mu))$-th entries of $Q^{[1, n] \circ \Seq} \cdot R^{\Seq \Join \Seq}$ as 
\begin{eqnarray}
\sum_{(i, j) \in \Seq} ( \xi \circ (i, j) ) \cdot ( (i, j) \Join (\lambda, \mu) )
 &=& \sum_{(i, j) \in \Seq}  q_{\xi ,i} q_{\xi, j} \cdot ( r_{i, \lambda} r_{j, \mu} + r_{i, \mu} r_{j, \lambda} )   \nonumber \\
 &=&  \sum_{(i, j) \in \Seq}  q_{\xi ,i} q_{\xi, j}  r_{i, \lambda} r_{j, \mu}
           + \sum_{(i, j) \in \Seq}  q_{\xi ,i} q_{\xi, j}  r_{i, \mu} r_{j, \lambda} .  
\end{eqnarray}
We can calculate the $(\xi, (\lambda, \mu))$-th entries of $Q^\square \cdot R^{[1, n] \circ \Seq}$ as 
\begin{eqnarray}
\sum_{i = 1}^n q_{\xi, i}^2 \cdot (i \circ (\lambda, \mu)) = \sum_{i = 1}^n q_{\xi, i}^2 \cdot r_{i, \lambda} r_{i, \mu} . 
\end{eqnarray}
Hence, by (1), (2), (3), we have the desired equality. 
\QED

\begin{lem}
Let $\bm{a} = (a_1, \ldots, a_n) \in k[T]^n$, let $Q \in \Mat(n, k)$, and let $\bm{b} := \bm{a} \, Q$. 
Let $\Sigma, \Theta, R$ be as in {\rm Lemma 3.20}. 
Then we have 
\[
\Delta_{\bm{a}, Q, \Sigma_+} + \Delta_{\bm{b}, R, \Theta_+} = \Delta_{\bm{a}, Q\cdot R, \Theta_+}. 
\]
\end{lem}

\Proof Using Lemmas 3.19, 3.20, 3,21, we have 
\begin{eqnarray*}
\lefteqn{\Delta_{\bm{a}, Q, \Sigma_+} + \Delta_{\bm{b}, R, \Theta_+} } \\ 
 &= & (a_1^2, \ldots, a_n^2) \cdot Q^{[1, n] \circ \Seq} \cdot {^t}(\sigma_{\lambda, \mu})_{(\lambda, \mu) \in \Seq} 
 + (b_1^2, \ldots, b_n^2) \cdot R^{[1, n] \circ \Seq} \cdot  {^t}(\tau_{\lambda, \mu})_{(\lambda, \mu) \in \Seq}  
\\
 & = &  (a_1^2, \ldots, a_n^2) \cdot Q^{[1, n] \circ \Seq} \cdot R^{\Seq \Join \Seq} \cdot {^t}(\tau_{\lambda, \mu})_{(\lambda, \mu) \in \Seq} + (a_1^2, \ldots, a_n^2) \cdot Q^\square \cdot R^{[1, n] \circ \Seq} \cdot  {^t}(\tau_{\lambda, \mu})_{(\lambda, \mu) \in \Seq}  \\ 
 &=& (a_1^2, \ldots, a_n^2) \cdot (Q^{[1, n] \circ \Seq} \cdot R^{\Seq \Join \Seq} + Q^\square \cdot R^{[1, n] \circ \Seq}) 
\cdot {^t}(\tau_{\lambda, \mu})_{(\lambda, \mu) \in \Seq} \\
 &=& (a_1^2, \ldots, a_n^2) \cdot (Q \cdot R)^{[1, n] \circ \Seq}  \cdot {^t}(\tau_{\lambda, \mu})_{(\lambda, \mu) \in \Seq} \\
 &=& \Delta_{\bm{a}, Q\cdot R, \Theta}. 
\end{eqnarray*}
\QED

For any $g \in {^\ell}G^{r_1}_{r_2}$, let $g_{(1, 1)}$ be the submatrix of $g$ obtained by deleting from $g$ the rows and columns through $g_{1, 1}$, i.e., 
\[
g_{(1, 1)} 
:=
\left(
\begin{array}{ c | c | c | c}
 g_{2, 2} & g_{2, 3} & g_{2, 4} & g_{2, 5} \\
\hline 
  O & g_{3, 3} & O & g_{3, 5} \\
\hline 
  O & O &  g_{4, 4} & g_{4, 5} \\
\hline 
  O & O & O & g_{5, 5}
\end{array}
\right) \in GL(\ell + r_2 + 1, k) . 
\]

Now, we prove Lemma 3.18. 
In order to prove that ${^\ell}G^{r_1}_{r_2}$ acts from right on ${^\ell}\Y^{r_1}_{r_2}$, 
we choose an arbitrary element $y$ of ${^\ell}\Y^{r_1}_{r_2}$ 
and arbitrary two elements $g, h$ of ${^\ell}G^{r_1}_{r_2}$. 
Write $y, y \star g, (y \star g) \star h \in {^\ell}\Y^{r_1}_{r_2}$ as 
\begin{eqnarray*}
\left\{ 
\begin{array}{ccl}
y &=& (S, \bm{a}_1, \bm{a}_2, \bm{\alpha}_2, \alpha) , \\
y \star g &=& (S', \bm{b}_1, \bm{b}_2, \bm{\beta}_2, \beta) , \\
(y \star g) \star h  &=& (S'', \bm{c}_1, \bm{c}_2, \bm{\gamma}_2, \gamma) . 
\end{array}
\right. 
\end{eqnarray*} 
So, we have 
\begin{eqnarray*}
\lefteqn{(y \star g) \star h } \\
 &=& (h_{2, 2}^{-1} \, g_{2, 2}^{-1} \, S \, g_{1, 1} \, h_{1, 1}, \;  (\bm{a}_1, \bm{a}_2, \bm{\alpha}_2, \alpha) g_{(1, 1)} h_{(1, 1)}
        - (0, 0, 0, \Delta_{\bm{a}_1, \; g_{2, 2}, \; S'_+} \cdot h_{5, 5} + \Delta_{\bm{b}_1, \; h_{2, 2}, \; S''_+}))  
\end{eqnarray*}
and 
\begin{eqnarray*}
y \star (g \cdot h) 
 = ( (g_{2, 2} \, h_{2, 2} )^{-1} \, S \, g_{1, 1} \, h_{1,1}, \; (\bm{a}_1, \bm{a}_2, \bm{\alpha}_2, \alpha) (gh)_{(1, 1)} 
 - (0, 0, 0, \Delta_{\bm{a}_1,\;  g_{2, 2} \cdot h_{2, 2}, \; S''_+})) . 
\end{eqnarray*}
For proving $(y \star g) \star h = y \star (g \cdot h)$, we have only to show that  
\[
 \Delta_{\bm{a}_1, \;  g_{2, 2}, \; S'_+} \cdot h_{5, 5} + \Delta_{\bm{b}_1, \; h_{2, 2}, \; S''_+}
 =  \Delta_{\bm{a}_1, \; g_{2, 2} \cdot h_{2, 2}, \; S''_+ }, 
\] 
where $S' = g_{2, 2}^{-1} \, S \, g_{1, 1}$ and $S'' = h_{2, 2}^{-1} \, g_{2, 2}^{-1} \, S \, g_{1, 1} \, h_{1, 1}$. 
Note that $h_{5, 5} \, S' = h_{2, 2} \cdot S'' \cdot {^t}h_{2, 2}$. 
We know from Lemma 3.22 that 
\begin{eqnarray*}
 \Delta_{\bm{a}_1, \; g_{2, 2}, \; S'_+} \cdot h_{5, 5} + \Delta_{\bm{b}_1, \; h_{2, 2}, \; S''_+}  
 = \Delta_{\bm{a}_1, \; g_{2, 2}, \; h_{5, 5} \, S'_+} + \Delta_{\bm{b}_1, \; h_{2, 2}, \; S''_+}  
 = \Delta_{\bm{a}_1, \; g_{2, 2} \cdot h_{2, 2}, \; S''_+} . 
\end{eqnarray*}
Clearly, $y \star I_{\ell + r_1 + r_2 + 1} = y$. 
Thus, the assertion (1) is proved.

We next prove the assertion (2). 
Choose an arbitrary element $x$ of ${^\ell}\X^{r_1}_{r_2}$ and an arbitrary element $g$ of ${^\ell}G^{r_1}_{r_2}$. 
Write $x = (S, \bm{a}_1, \bm{a}_2, \bm{\alpha}_2, \alpha)$. 
Therefore, $(\bm{a}_1, \bm{a}_2, \bm{\alpha}_2) \in \F(r_1, \ell - r _1, r_2, \mf{P})$. 
Since 
\[
 (\bm{b}_1, \bm{b}_2)
 = (\bm{a}_1, \bm{a}_2) 
 \left( 
 \begin{array}{cc}
 g_{2, 2} & g_{2, 3} \\
 O & g_{3, 3}
 \end{array}
 \right)
\qquad \text{ and } \qquad  
(\bm{b}_1, \bm{\beta}_2)
 = 
(\bm{a}_1, \bm{\alpha}_2)
\left(
\begin{array}{cc}
 g_{2, 2} & g_{2, 4} \\
  O &  g_{4, 4}  
\end{array}
\right), 
\]
we have $(\bm{b}_1, \bm{b}_2, \bm{\beta}_2) \in \F(r_1, \ell - r_1, r_2, \mf{P})$. 
\QED

Now, we have the following theorem: 

\begin{thm}
For all $(\ell, r_1, r_2) \in \Omega_m$, we have 
\begin{eqnarray*}
 ({^\ell} \mf{H}^{r_1}_{r_2} )^E / \sim & \cong & {^\ell}\X^{r_1}_{r_2} \, /\, {^\ell} G^{r_1}_{r_2} . 
\end{eqnarray*}
\end{thm}

\subsubsection{$(\daleth_{m + 2})^E$}

We are interested in an equivalence relation of two exponential matrices belonging to any one of 
the five sets $(\daleth_{m + 2}^{\rm II})^E$, $(\daleth_{m + 2}^{\rm III})^E$, $(\daleth_{m + 2}^{\rm IV})^E$, 
$(\daleth_{m + 2}^{\rm V})^E$, $(\daleth_{m + 2}^{\rm VI})^E$. 

We obtain the following theorem from Theorems 2.14, 2.15, 2.17, 2.19, 2.22.

\begin{thm}
The following assertions {\rm (1), (2), (3), (4), (5)} hold true. 
\begin{enumerate} 
\item[\rm (1)] 
Let $A(T) = \Lambda(1, j; \alpha(T)) \in (\daleth_{m + 2}^{\rm II})^E$ 
and $B(T) = \Lambda(1, j'; \beta(T)) \in (\daleth_{m + 2}^{\rm II})^E$, where $1 \leq j, j' \leq m + 1$. 
Then the following conditions {\rm (1.1)} and {\rm (1.2)} are equivalent: 
\begin{enumerate}
\item[\rm (1.1)] $A(T)$ and $B(T)$ are equivalent. 
\item[\rm (1.2)] $j = j'$, and $\beta(T) = \alpha(T) Q$ for some $Q \in GL(j, k)$. 
\end{enumerate}

\item[\rm (2)] 
Let $A(T) = \Lambda(i, 1; \alpha(T)) \in (\daleth_{m + 2}^{\rm III})^E$ and 
$B(T) = \Lambda(i', 1; \beta(T)) \in (\daleth_{m + 2}^{\rm III})^E$, where $2 \leq i, i' \leq m + 1$. 
Then the following conditions {\rm (2.1)} and {\rm (2.2)} are equivalent: 
\begin{enumerate}
\item[\rm (2.1)] $A(T)$ and $B(T)$ are equivalent. 
\item[\rm (2.2)] $i = i'$, and $\beta(T) = Q \cdot \alpha(T)$ for some $Q \in GL(i, k)$. 
\end{enumerate}

\item[\rm (3)] 
Let $A(T) = \Lambda(2, 2; \alpha(T)) \in (\daleth_{m + 2}^{\rm IV})^E$ 
and let $B(T) = \Lambda(2, 2; \beta(T)) \in (\daleth_{m + 2}^{\rm IV})^E$, where 
$\alpha(T) = L(a_1 \mid a_2 \mid 0) \in {^1}\mf{L}_1^0$ 
and $\beta(T) = L(b_1 \mid b_2 \mid 0) \in {^1}\mf{L}_1^0$, where $a_1, a_2, b_1, b_2 \in \mf{P}$.   
Then the following conditions {\rm (3.1)} and {\rm (3.2)} are equivalent: 
\begin{enumerate}
\item[\rm (3.1)] $A(T)$ and $B(T)$ are equivalent. 
\item[\rm (3.2)] $(b_1, b_2) = (a_1, a_2) \, g $ for some $g \in GLDD(2, k)$. 
\end{enumerate}

\item[\rm (4)] 
Let $A(T) = \Lambda(2, 2; \alpha(T)) \in (\daleth_{m + 2}^{\rm V})^E$ 
and let $B(T) = \Lambda(2, 2; \beta(T)) \in (\daleth_{m + 2}^{\rm V})^E$, where 
$\alpha(T) = L(a_1 \mid a_2 \mid a) \in ({^1}\mf{L}_1)^{\dim = 3}$, 
$\beta(T) = L(b_1 \mid b_2 \mid b) \in ({^1}\mf{L}_1)^{\dim = 3}$, and $a_1, a_2, a, b_1, b_2, b \in \mf{P}$.   
Then the following conditions {\rm (4.1)} and {\rm (4.2)} are equivalent: 
\begin{enumerate}
\item[\rm (4.1)] $A(T)$ and $B(T)$ are equivalent. 
\item[\rm (4.2)] $(b_1, b_2, b) = (a_1, a_2, a) \, g$ for some $g \in V(1, 1, 1; k)$. 
\end{enumerate}

\item[\rm (5)] 
Let $A(T) = \Lambda(j + 1, i + 1; \alpha(T) ) \in ( \daleth_{m + 2}^{\rm VI} )^E$ 
and let $B(T) = \Lambda( j' + 1, i' + 1; \beta(T) ) \in ( \daleth_{m + 2}^{\rm VI})^E$, where 
$\alpha(T) = L(a_1, \ldots, a_i \mid a_{i + 1}, \ldots, a_{i + j} \mid a) \in {^i}\mf{L}_j$, $a_1, \ldots, a_{i + j}, a \in \mf{P}$,  
and $\beta(T) = L(b_1, \ldots, b_{i'} \mid b_{i' + 1}, \ldots, b_{i' + j'} \mid b) \in {^{i'}}\mf{L}_{j'}$, $b_1, \ldots, b_{i' + j'}, b \in \mf{P}$.   
Then the following conditions {\rm (5.1)} and {\rm (5.2)} are equivalent: 
\begin{enumerate}
\item[\rm (5.1)] $A(T)$ and $B(T)$ are equivalent. 
\item[\rm (5.2)] $(i, j) = (i', j')$, and 
$(b_1, \ldots, b_i, b_{i + 1}, \ldots, b_{i + j}, b) = (a_1, \ldots, a_i, a_{i + 1}, \ldots, a_{i + j}, a) \, g$ for some $g \in V(i, j, 1)$. 
\end{enumerate} 
\end{enumerate}
\end{thm}

\section{Exponential matrices of size four-by-four}

\subsection{$\E(4, k[T])^E$}

\begin{thm} 
We have 
\[
\E(4, k[T])^E
 \leadsto 
 \mf{U}_{[4]}^E \; \cup \;  
 \mf{U}_{[3, 1]}^E \; \cup \;  
 \mf{U}_{[1, 3]}^E \; \cup \; 
  \mf{A}(1, 3)^E \; \cup \; \mf{A}(2, 2)^E \; \cup \; \mf{A}(3, 1)^E \; \cup \;  
 ({^2}\mf{H}^2_0)^E. 
\]
\end{thm}

\Proof 
We know from Lemma 1.12 that 
\[
 \E(4, k[T]) \leadsto \mf{U}_4^E
 = \bigcup_{\text{ $[d_1, \ldots, d_t]$ is an ordered partition of $4$} } \mf{U}_{[d_1, \ldots, d_t]}^E. 
\]
We can make a list of all ordered partitions of $4$ as follows: 
\begin{eqnarray*}
[4], \qquad [3, 1], \qquad [2, 2], \qquad [1, 3], \qquad [2, 1, 1], \qquad [1, 2, 1] \qquad [1, 1, 2] ,  \qquad [1, 1, 1, 1] . 
\end{eqnarray*}
Thus 
\[
\E(4, k[T])
 \leadsto
 \mf{U}_{[4]}^E \; \cup \;  
 \mf{U}_{[3, 1]}^E \; \cup \;  
 \mf{U}_{[2, 2]}^E \; \cup \;  
 \mf{U}_{[1, 3]}^E \; \cup \; 
 \mf{U}_{[2, 1, 1]}^E \; \cup \;   
 \mf{U}_{[1, 2, 1]}^E \; \cup \;  
 \mf{U}_{[1, 1, 2]}^E \; \cup \;  
 \mf{U}_{[1, 1, 1, 1]}^E . 
\]
Note that $\mf{U}_{[2, 2]}^E \cup \mf{U}_{[2, 1, 1]}^E \cup \mf{U}_{[1, 1, 2]}^E \cup \mf{U}_{[1, 1, 1, 1]}^E  \subset \mf{H}(4, k[T])^E$ and 
$\mf{U}_{[1, 2, 1]}^E \subset \mf{A}(2, 2)^E$. 
Recalling Theorem 3.1, we have 
\begin{eqnarray*}
\lefteqn{\E(4, k[T]) } \\
& \leadsto & 
 \mf{U}_{[4]}^E \; \cup \;  
 \mf{U}_{[3, 1]}^E \; \cup \;  
 \mf{U}_{[1, 3]}^E \; \cup \; 
 \mf{H}(4, k[T])^E \; \cup \;   
 \mf{A}(2, 2)^E  \\
& \leadsto & 
 \mf{U}_{[4]}^E \; \cup \;  
 \mf{U}_{[3, 1]}^E \; \cup \;  
 \mf{U}_{[1, 3]}^E \; \cup \; 
 \left(  (\daleth_4)^E \; \cup \; \bigcup_{(\ell, r_1, r_2) \in \Omega_2} ({^\ell}\mf{H}^{r_1}_{r_2})^E \right) \; \cup \;   
 \mf{A}(2, 2)^E .  
\end{eqnarray*} 
Since $\daleth_4 \subset  \mf{A}(1, 3) \, \cup \, \mf{A}(2, 2) \, \cup \, \mf{A}(3, 1)$, 
we have 
\[
(\daleth_4)^E \subset  \mf{A}(1, 3)^E \, \cup \, \mf{A}(2, 2)^E \, \cup \, \mf{A}(3, 1)^E. 
\] 
Note that $\Omega_2 =  \{ (1, 1, 0), \; (1, 1, 1), \;  (2, 1, 0), \;  (2, 2, 0)\}$ and that 
\begin{eqnarray*}
({^1}\mf{H}^1_0)^E \leadsto (\mf{U}_{[3, 1]})^E, \qquad 
({^1}\mf{H}^1_1)^E \leadsto (\mf{U}_{[1, 3]})^E, \qquad 
({^2}\mf{H}^1_0)^E \leadsto (\mf{U}_{[3, 1]})^E . 
\end{eqnarray*}
Hence we have 
\begin{eqnarray*}
\E(4, k[T])^E 
\leadsto 
 \mf{U}_{[4]}^E \; \cup \;  
 \mf{U}_{[3, 1]}^E \; \cup \;  
 \mf{U}_{[1, 3]}^E \; \cup \; 
  \mf{A}(1, 3)^E \; \cup \; \mf{A}(2, 2)^E \; \cup \; \mf{A}(3, 1)^E \; \cup \;  
 ({^2}\mf{H}^2_0)^E . 
\end{eqnarray*} 
\QED

The following is a corollary of the above Theorem 4.1: 

\begin{cor} 
We have 
\begin{eqnarray*}
\lefteqn{\E(4, k[T])^E} \\
 &\leadsto&
\left\{ 
\begin{array}{rl}

 \mf{A}(1, 3)^E \; \cup \; \mf{A}(2, 2)^E \; \cup \; \mf{A}(3, 1)^E \; \cup \;  
 ({^2}\mf{H}^2_0)^E   & \text{ if \quad $p = 2$} , \\

 \mf{J}_{[3, 1]}^E \; \cup \;  
 \mf{J}_{[1, 3]}^E \; \cup \; 
 \mf{A}(1, 3)^E \; \cup \; \mf{A}(2, 2)^E \; \cup \; \mf{A}(3, 1)^E \; \cup \;  
 ({^2}\mf{H}^2_0)^E   & \text{ if \quad $p = 3$}, \\

 \mf{J}_{[4]}^E \; \cup \;  
 \mf{J}_{[3, 1]}^E \; \cup \;  
 \mf{J}_{[1, 3]}^E \; \cup \; 
 \mf{A}(1, 3)^E \; \cup \; \mf{A}(2, 2)^E \; \cup \; \mf{A}(3, 1)^E \; \cup \;  
 ({^2}\mf{H}^2_0)^E  & \text{ if \quad $p \geq 5$} . 

\end{array}
\right. 
\end{eqnarray*}
\end{cor}

\Proof Recall Corollaries 1.15, 1.22 and 1.26. 
\QED

We know forms of exponential matrices belonging to any one of the sets $\mf{J}_{[4]}^E$ (see Lemma 1.16), 
$\mf{J}_{[3, 1]}^E$ (see Lemma 1.23), 
$\mf{J}_{[1, 3]}^E$ (see Lemma 1.27),  
$\mf{A}(1, 3)^E$, $\mf{A}(2, 2)^E$, $\mf{A}(3, 1)^E$ (see Lemma 1.28), and 
$({^2}\mf{H}^2_0)^E$ (see Theorem 3.3). 

Now we can describe the following classification of exponential matrices of size four-by-four:

\subsection{A classification of exponential matrices of size four-by-four, up to equivalence}

\subsubsection{$p = 2$}

\begin{thm}
Assume $p = 2$. Any exponential matrix of $\Mat(4, k[T])$ is equivalent to one of the following exponential matrices 
{\rm (1), (2), (3), (4)}: 
\begin{enumerate}

\item[\rm (1)] 
$
\left(
\begin{array}{cccc}
 1 & a & b & c \\
 0 & 1 & 0 & 0 \\
 0 & 0 & 1 & 0 \\ 
 0 & 0 & 0 & 1 
\end{array}
\right)
\qquad 
(\text{ $a, b, c$ are $p$-polynomials }) . 
$

\item[\rm (2)] 
$
\left(
\begin{array}{cccc}
 1 & 0 & a & b \\
 0 & 1 & c & d \\
 0 & 0 & 1 & 0 \\ 
 0 & 0 & 0 & 1 
\end{array}
\right)
\qquad 
(\text{ $a, b, c, d$ are $p$-polynomials }) . 
$

\item[\rm (3)] 
$
\left(
\begin{array}{cccc}
 1 & 0 & 0 & a \\
 0 & 1 & 0 & b \\
 0 & 0 & 1 & c \\ 
 0 & 0 & 0 & 1 
\end{array}
\right)
\qquad 
(\text{ $a, b, c$ are $p$-polynomials }) . 
$

\item[\rm (4)] 
$
\left( 
\begin{array}{cccc}
 1 & a & b & \mu ab+ c \\
 0 & 1 & 0 & \mu b \\
 0 & 0 & 1 & \mu a \\
 0 & 0 & 0 & 1
\end{array}
\right)  \qquad 
\left(  
\begin{array}{l}
\text{$a, b, c$ are $p$-polynomials, } \\
\text{$a, b$ are linearly independent over $k$}, \\
\text{$\mu \in k$}
\end{array}
\right) . 
$
\end{enumerate}
\end{thm}

\subsubsection{$p = 3$}

\begin{thm}
Assume $p  =3$. Any exponential matrix $A(T)$ of $\Mat(4, k[T])$ is equivalent to one of the following exponential matrices 
{\rm (1), (2), (3), (4), (5), (6)}: 
\begin{enumerate}

\item[\rm (1)] 
$
\left(
\begin{array}{cccc}
 1 & a & \frac{1}{2} a^2 + b & c \\
 0 & 1 & a & 0 \\
 0 & 0 & 1 & 0 \\ 
 0 & 0 & 0 & 1 
\end{array}
\right)
\qquad 
\left(  
\begin{array}{l}
\text{$a, b, c$ are $p$-polynomials}, \\
\text{$a \ne 0$}
\end{array}
\right) . 
$

\item[\rm (2)] 
$
\left(
\begin{array}{cccc}
 1 & 0 & 0 & c \\
 0 & 1 & a & \frac{1}{2} a^2 + b \\
 0 & 0 & 1 & a \\ 
 0 & 0 & 0 & 1 
\end{array}
\right)
\qquad 
\left(  
\begin{array}{l}
\text{$a, b, c$ are $p$-polynomials}, \\
\text{$a \ne 0$}
\end{array}
\right) . 
$

\item[\rm (3)] 
$
\left(
\begin{array}{cccc}
 1 & a & b & c \\
 0 & 1 & 0 & 0 \\
 0 & 0 & 1 & 0 \\ 
 0 & 0 & 0 & 1 
\end{array}
\right)
\qquad 
(\text{ $a, b, c$ are $p$-polynomials }) . 
$

\item[\rm (4)] 
$
\left(
\begin{array}{cccc}
 1 & 0 & a & b \\
 0 & 1 & c & d \\
 0 & 0 & 1 & 0 \\ 
 0 & 0 & 0 & 1 
\end{array}
\right)
\qquad 
(\text{ $a, b, c, d$ are $p$-polynomials }) . 
$

\item[\rm (5)] 
$
\left(
\begin{array}{cccc}
 1 & 0 & 0 & a \\
 0 & 1 & 0 & b \\
 0 & 0 & 1 & c \\ 
 0 & 0 & 0 & 1 
\end{array}
\right)
\qquad 
(\text{ $a, b, c$ are $p$-polynomials }) . 
$

\item[\rm (6)] 
$
\left( 
\begin{array}{cccc}
 1 & a & b & \frac{\lambda}{2}  a^2 + \mu a b  + \frac{\nu}{2} b^2 + c \\
 0 & 1 & 0 & \lambda a + \mu  b \\
 0 & 0 & 1 & \mu a + \nu b \\
 0 & 0 & 0 & 1
\end{array}
\right)  \qquad 
\left(  
\begin{array}{l}
\text{$a, b, c$ are $p$-polynomials, } \\
\text{$a, b$ are linearly independent over $k$}, \\
\text{$\lambda, \mu ,\nu \in k$, and $\lambda \nu - \mu^2\ne 0$}
\end{array}
\right) . 
$
\end{enumerate}
\end{thm}

\subsubsection{$p \geq 5$}

\begin{thm}
Assume $p  \geq 5$. Any exponential matrix $A(T) \in \Mat(4, k[T])$ is equivalent to one of the following exponential matrices 
{\rm (1), (2), (3), (4), (5), (6), (7)}: 
\begin{enumerate}

\item[\rm (1)] 
$
\left( 
\begin{array}{cccc}
 1 & a & \frac{1}{2} a^2 + b & \frac{1}{6} a^3 + a b + c  \\
 0 & 1 & a & \frac{1}{2} a^2 + b \\
 0 & 0 & 1 & a \\
 0 & 0 & 0 & 1
\end{array}
\right) 
\qquad 
\left(  
\begin{array}{l}
\text{$a, b, c$ are $p$-polynomials}, \\
\text{$a \ne 0$}
\end{array}
\right) . 
$

\item[\rm (2)] 
$
\left(
\begin{array}{cccc}
 1 & a & \frac{1}{2} a^2 + b & c \\
 0 & 1 & a & 0 \\
 0 & 0 & 1 & 0 \\ 
 0 & 0 & 0 & 1 
\end{array}
\right)
\qquad 
\left(  
\begin{array}{l}
\text{$a, b, c$ are $p$-polynomials}, \\
\text{$a \ne 0$}
\end{array}
\right) . 
$

\item[\rm (3)] 
$
\left(
\begin{array}{cccc}
 1 & 0 & 0 & c \\
 0 & 1 & a & \frac{1}{2} a^2 + b \\
 0 & 0 & 1 & a \\ 
 0 & 0 & 0 & 1 
\end{array}
\right)
\qquad 
\left(  
\begin{array}{l}
\text{$a, b, c$ are $p$-polynomials}, \\
\text{$a \ne 0$}
\end{array}
\right) . 
$

\item[\rm (4)] 
$
\left(
\begin{array}{cccc}
 1 & a & b & c \\
 0 & 1 & 0 & 0 \\
 0 & 0 & 1 & 0 \\ 
 0 & 0 & 0 & 1 
\end{array}
\right)
\qquad 
(\text{ $a, b, c$ are $p$-polynomials }) . 
$

\item[\rm (5)] 
$
\left(
\begin{array}{cccc}
 1 & 0 & a & b \\
 0 & 1 & c & d \\
 0 & 0 & 1 & 0 \\ 
 0 & 0 & 0 & 1 
\end{array}
\right)
\qquad 
(\text{ $a, b, c, d$ are $p$-polynomials }) . 
$

\item[\rm (6)] 
$
\left(
\begin{array}{cccc}
 1 & 0 & 0 & a \\
 0 & 1 & 0 & b \\
 0 & 0 & 1 & c \\ 
 0 & 0 & 0 & 1 
\end{array}
\right)
\qquad 
(\text{ $a, b, c$ are $p$-polynomials }) . 
$

\item[\rm (7)] 
$
\left( 
\begin{array}{cccc}
 1 & a & b & \frac{\lambda}{2}  a^2 + \mu a b  + \frac{\nu}{2} b^2 + c \\
 0 & 1 & 0 & \lambda a + \mu  b \\
 0 & 0 & 1 & \mu a + \nu b \\
 0 & 0 & 0 & 1
\end{array}
\right)  \qquad 
\left(  
\begin{array}{l}
\text{$a, b, c$ are $p$-polynomials, } \\
\text{$a, b$ are linearly independent over $k$}, \\
\text{$\lambda, \mu ,\nu \in k$, and $\lambda \nu - \mu^2\ne 0$}
\end{array}
\right) . 
$
\end{enumerate}
\end{thm}

\section{On modular representations of elementary abelian $p$-groups}

\subsection{$\Rep(E, n)$, $\Rep(E_r, n)$}

For any group $G$ and any integer $n \geq 1$, we denote by $\Rep(G, n)$ the set of all 
$n$-dimensional representations $\rho : G \to GL(n, k)$ of $G$. 

Let 
\[
E := \bigoplus_{i \geq 0} \Z/p\Z
\]
be the direct sum of a family of $p$-cyclic groups indexed by the set of all non-negative integers. 
For any integer $r \geq 1$, let $E_r$ be the group
\[
E_r := \bigoplus_{i = 0}^{r - 1} \Z/p\Z.  
\]
We say that $E_r$ is an {\it elementary abelian $p$-group of rank $r$}. 
Defining an inclusion $i_{r, \infty} : E_r \hookrightarrow E$ as $([i_0], \ldots, [i_{r-1}]) \mapsto ([i_0], \ldots, [i_{r - 1}], [0], [0], \ldots)$, 
we can define a map $\pi_r : \Rep(E, n) \to \Rep(E_r, n)$ as $\pi_r(\rho) :=  \rho \circ i_{r, \infty}$. 
The map $\pi_r$ has a right split map $\iota_r : \Rep(E_r, n) \to \Rep(E, n)$, i.e., $\pi_r \circ \iota_r = \id_{\Rep(E_r, n)}$. 
In fact, define $\iota_r : \Rep(E_r, n) \to \Rep(E, n)$ as 
\[
 \iota_r(\varrho) := \rho, 
\qquad 
\rho(e_i) := 
\begin{cases}
\varrho(e_i) \quad & \text{ if  \quad $0 \leq i \leq r - 1$}, \\
I_n   & \text{ if \quad $i \geq r$},  \\
\end{cases}
\] 
where each $e_i$ is the element of $E$ such that the $i$-th entry is $[1]$ and the other entires are $[0]$. 
Clearly, $\iota_r$ is injective. 
Thus we have $\bigcup_{r \geq 1} \Rep(E_r, n) \subset \Rep(E, n) $, and 
we can define a map $f : \bigcup_{r \geq 1} \Rep(E_r, n)  \to \E(n, k[T])$ as 
\[
 f(\rho) := \prod_{i \geq 1} \Exp( T^{p^{i - 1}} \cdot(  \rho(e_{i - 1}) - I_n) ) . 
\] 
Clearly, $f$ is bijective. 

For all integers $n, r \geq 1$, we define the sets ${\cal U}_{n, r}$, $\cU_n$, ${\cal N}_{n, r}$ and $\cN_n$ as 
\begin{eqnarray*}
\left\{
\begin{array}{rcl}
 {\cal U}_{n, r}  &:=& \ds \left\{ (U_i)_{i \geq 1} \in \prod_{i \geq 1} \Mat(n, k) \; \left| \; 
\begin{array}{ll}
 U_i^p = I_n                   & (\, 1 \leq \forall i \leq r \,), \\
 U_j U_\ell = U_\ell U_j \; & (\, 1 \leq \forall j, \ell \leq r \,), \\
 U_i = I_n & (\, \forall i \geq r + 1) 
\end{array}
\right. 
\right\} , \\ [8mm]

\cU_n &:=& \ds \bigcup_{r \geq 1} \cU_{n, r}, \\  [8mm]

 {\cal N}_{n, r} &:=& \ds \left\{  (N_i)_{i \geq 1} \in \prod_{i \geq 1} \Mat(n, k) \; 
\left| \; 
\begin{array}{ll}
 N_i^p = O_n                 & (\, 1 \leq \forall i \leq r \,), \\ 
 N_j N_\ell = N_\ell N_j \; & (\, 1 \leq \forall j, \ell \leq r \,) , \\
 N_i = O_n & (\, \forall i \geq r + 1 \, )  
\end{array} 
\right. 
\right\} , \\ [8mm]

\cN_n &:=& \ds \bigcup_{r \geq 1} \cN_{n, r}. \\ 
\end{array}
\right. 
\end{eqnarray*}

We can define a map $g : \bigcup_{r \geq 1} \Rep(E_r, n) \to {\cal U}_n$
 as $g(\rho) := (\rho(e_{ i - 1}) )_{ i \geq 1}$, 
and can also define a map $h : {\cal U}_n \to {\cal N}_n$ as $h( (U_i)_{i \geq 1} ) := ( U_i - I_n )_{i \geq 1}$.   
Clearly, both $g$ and $h$ are bijective. 
Let $\exp : {\cal N}_n \to \E(n, k[T])$ be the bijection defined by $\exp := f \circ g^{-1} \circ h^{-1}$, i.e., 
\[
 \exp( N_1, N_2, \ldots ) := \prod_{i \geq 1} \Exp(T^{p^{i  - 1}} N_i) . 
\]
We have the following commutative diagram: 
\begin{eqnarray*}
\xymatrix{
  {\cal N}_{n, r} \ar[r]^{\subset} & {\cal N}_n   \ar[rdd]^{\exp}_{\cong}  \\
  {\cal U}_{n, r} \ar[r]^{\subset} \ar[u]^h_{\cong} & {\cal U}_n  \ar[u]^h_{\cong} \\
  \Rep(E_r, n) \ar[r] \ar[u]^g_{\cong} & \ds\bigcup_{r \geq 1} \Rep(E_r, n)  \ar[r]_(0.55){f}^(0.55){\cong}  \ar[u]^{g}_{\cong}   &  \E(n, k[T]) 
}
\end{eqnarray*}

Define an equivalence relation $\sim$ on $\prod_{i \geq 1} \Mat(n, k)$ as follows: 
Two elements $(A_i)_{i \geq 1}$ and $(B_i)_{i \geq 1}$ of $\prod_{i \geq 1} \Mat(n, k)$ are {\it equivalent} 
if there exists $P \in GL(n, k)$ such that $(P^{-1} A_i P)_{i \geq 1} = (B_i)_{i \geq 1}$. 
We similarly define equivalence relations on $\cU_n$, $\cU_{n, r}$, $\cN_n$, $\cN_{n, r}$. 
Then the following commutative diagram is induced: 
\begin{eqnarray*}
\xymatrix{
 {\cal N}_{n, r}/ \sim \ar[r]^{\subset}           &  {\cal N}_n / \sim \ar[rdd]^{\exp} \\
 {\cal U}_{n, r} / \sim \ar[r]^{\subset}  \ar[u]^{h}_{\cong}      &   {\cal U}_n / \sim  \ar[u]_h^{\cong}   \\
 \Rep(E_r, n) / \sim \ar[r]^(0.45){\iota_r} \ar[u]^{g}_\cong     & \ds \bigcup_{r \geq 1} \Rep(E_r, n) / \sim \ar[r]^(0.55)f_(0.55){\cong}  \ar[u]_g^{\cong} &  \E(n, k[T]) / \sim
}
\end{eqnarray*}

\subsection{$\mf{j}_{[n]}$, $\mf{j}_{[n, 1]}$, $\mf{j}_{[1, n]}$, $\mf{a}(i_1, i_2, i_3)$, ${^\ell}\mf{h}^{r_1}_{r_2}(S)$, 
$\daleth_{m + 2}^\flat$}

For any $A(T) \in \E(n, k[T])$, there exists $(N_1, \ldots, N_r, O_n, O_n, \ldots) \in \cN_{n, r}$ such that 
$A(T) = \exp(N_1, \ldots, N_r, O_n, O_n, \ldots )$ (see Lemma 1.5). 
For any exponential matrix $A(T)$ belonging to one of the five sets 
$\mf{J}_{[n]}^E$, $\mf{J}_{[n, 1]}^E$, $\mf{J}_{[1, n]}^E$, $\mf{A}(i_1, i_2, i_3)^E$, $( {^\ell}\mf{H}^{r_1}_{r_2} )^E$, $(\daleth_{m + 2})^E$, 
we consider the exponential expression of $A(T)$ (see the following Lemmas 5.1, 5.2, 5.3, 5.4, 5.5, 5.6). 
For the consideration, we introduce finite dimensional $k$-algebras $\mf{j}_{[n]}$, $\mf{j}_{[n, 1]}$, $\mf{j}_{[1, n]}$, $\mf{a}(i_1, i_2, i_3)$, 
${^\ell}\mf{h}^{r_1}_{r_2}(S)$, $\daleth_{m + 2}^\flat$.

In this subsection, if we say that a subset $\mf{g}$ of $\Mat(n, k)$ is a {\it subalgebra} of $\Mat(n, k)$, 
we do not assume $\mf{g}$ is unital. For any subalgebra $\mf{g}$ of $\Mat(n, k)$, 
we let $\mf{g}^{\oplus r}$ denote the direct sum of the $r$-copies of $\mf{g}$. 
We can regard $\mf{g}^{\oplus r}$ as a subalgebra of $\prod_{i \geq1} \Mat(n, k)$ through 
the inclusion $\mf{g}^{\oplus r} \to \prod_{i \geq 1} \Mat(n, k)$ defined by $(g_1, \ldots, g_r) \mapsto (g_1, \ldots, g_r, O_n, O_n, \ldots)$.  

\subsubsection{$\mf{j}_{[n]}$}

Let $n \geq 2$ be an integer. 
We define a subset $\mf{j}_{[n]}$ of $\Mat( n + 1, k)$ as 
\begin{eqnarray*}
\mf{j}_{[n]} := 
\left\{ 
\left.
\sum_{i = 1}^{n - 1} s_i \nu_n^i \in \Mat(n, k) \; \right| \;  s_1, \ldots, s_{n - 1} \in k 
\right\} . 
\end{eqnarray*}
Clearly, $\mf{j}_{[n]}$ becomes a commutative subalgebra of $\Mat(n, k)$.  
For all $N \in \mf{j}_{[n]}$, we have $N^n = O$.

\begin{lem} 
Assume $1 \leq n \leq p$. 
Then the following assertions {\rm (1)} and {\rm (2)} hold true: 
\begin{enumerate} 
\item[\rm (1)]${\mf{j}_{[n]}}^{\oplus r} \subset {\cal N}_{n, r}$ $(\forall r \geq 1)$. 
\item[\rm (2)] $\ds \mf{J}_{[n]}^E \subset \bigcup_{r \geq 1} \exp( {\mf{j}_{[n]}}^{\oplus r} )$. 
\end{enumerate} 
\end{lem}

\Proof The proof of assertion (1) is straightforward. We shall prove assertion (2).  
Choose any element $A(T)$ of $\mf{J}_{[n]}^E$. We can express $A(T)$ as 
$A(T) = \Exp(\sum_{i = 1}^{n - 1} f_i \nu_n^i)$ for some $p$-polynomials $f_i$ $(1 \leq i \leq n - 1)$ 
(see Lemma 1.16). 
We can express the polynomial matrix $\sum_{i = 1}^{n - 1} f_i \nu_n^i$ as 
$\sum_{i = 1}^{n - 1} f_i \nu_n^i = \sum_{i = 1}^{r} T^{p^{i - 1}} N_i$, where $N_1, \ldots, N_r \in \mf{j}_{[n]}$. 
Thus $A(T) = \exp(N_1, \ldots, N_r, O_n, O_n, \ldots ) \in \exp( {\mf{j}_{[n]} }^{\oplus r})$. 
\QED

\subsubsection{$\mf{j}_{[n, 1]}$}

Let $n \geq 2$ be an integer. 
We define a subset $\mf{j}_{[n, 1]}$ of $\Mat(n + 1, k)$ as 
\begin{eqnarray*}
\mf{j}_{[n, 1]} := 
\left\{
\left. 
\left(
\begin{array}{c | c}
 \ds \sum_{i = 1}^{n - 1} s_i \nu_n^i & s_n \bm{e}_1\\
\hline 
 \bm{0} & 0
\end{array}
\right)
\in \Mat( n + 1, k) 
\; 
\right|
\;
s_1, \ldots, s_n \in k  
\right\} . 
\end{eqnarray*}
Clearly, $\mf{j}_{[n, 1]}$ becomes a commutative subalgebra of $\Mat(n + 1, k)$. 
For all $N \in \mf{j}_{[n, 1]}$, we have $N^n = O$.

\begin{lem}
Assume $1 \leq n \leq p$. 
Then the following assertions {\rm (1)} and {\rm (2)} hold true: 
\begin{enumerate} 
\item[\rm (1)] ${ \mf{j}_{[n, 1]} }^{\oplus r} \subset {\cal N}_{n + 1, r}$ $( \forall r \geq 1 )$.  
\item[\rm (2)] $\ds \mf{J}_{[n, 1]}^E \subset \bigcup_{r \geq 1} \exp({ \mf{j}_{[n, 1]} }^{\oplus r})$. 
\end{enumerate} 
\end{lem} 

\Proof The proof of assertion (1) is straightforward. 
We shall prove assertion (2). By Lemma 5.1, we can write $A(T)$ as 
\[
 A(T)  
 = 
\left(
 \begin{array}{c | c}
\exp(N) & \bm{0} \\
\hline 
 \bm{0} & 1
 \end{array} 
\right)
\left(
 \begin{array}{c | c}
 I_n & f_n \bm{e}_1 \\
\hline 
 \bm{0} & 1
 \end{array} 
\right), 
\]
where $N = (N_i)_{i \geq 1} \in \bigcup_{r \geq 1} {\mf{j}_{[n]}}^{\oplus r} (\subset {\cal N}_n)$ and $f_n$ is a $p$-polynomial. 
Write $f_n = \sum_{i \geq 1} c_i T^{p^{i-1}}$, where $c_i \in k$ for all $i \geq 1$. 
Letting 
\[
 \widehat{N}_i := 
\left(
\begin{array}{c | c}
 N_i  & c_i \bm{e}_1 \\ 
\hline 
 \bm{0} & 0
\end{array} 
\right)  
\in \mf{j}_{[n, 1]} 
\qquad ( \forall i \geq 1) , 
\]
we have $A(T)  = \Exp\left(\sum_{i \geq 1} T^{p^{i-1}} \widehat{N}_i \right)$.   
Clearly, there exists an integer $r \geq 1$ such that $\widehat{N}_i = O$ for all $i \geq r + 1$. 
Thus $A(T) = \exp(\widehat{N}_1, \ldots, \widehat{N}_r, O_{n + 1}, O_{n + 1}, \ldots ) \in \exp( {\mf{j}_{[n, 1]}}^{\oplus r} )$. 
\QED

\subsubsection{$\mf{j}_{[1, n]}$}

Let $n \geq 2$ be an integer. 
We define a subset $\mf{j}_{[1, n]}$ of $\Mat(n + 1, k)$ as 
\begin{eqnarray*}
\mf{j}_{[1, n]} := 
\left\{
\left. 
\left(
\begin{array}{c | c}
0 & s_n \transpose\bm{e_n} \\ 
\hline 
\bm{0} & \ds \sum_{i = 1}^{n - 1} s_i \nu_n^i  
\end{array}
\right)
\in \Mat( n + 1, k) 
\; 
\right|
\;
s_1, \ldots, s_n \in k  
\right\} . 
\end{eqnarray*}
Clearly, $\mf{j}_{[1, n]}$ becomes a commutative subalgebra of $\Mat(n + 1, k)$. 
For all $N \in \mf{j}_{[1, n]}$, we have $N^n = O$.

\begin{lem}
Assume $1 \leq p \leq n$. 
Then the following assertions {\rm (1)} and {\rm (2)} hold true: 
\begin{enumerate} 
\item[\rm (1)] ${\mf{j}_{[1, n]}}^{\oplus r} \subset {\cal N}_{n + 1, r}$ $( \forall r \geq 1 )$.  
\item[\rm (2)] $\ds \mf{J}_{[1, n]}^E \subset \bigcup_{r \geq 1} \exp({ \mf{j}_{[1, n]} }^{\oplus r})$. 
\end{enumerate} 
\end{lem}

\Proof Arguing as in the proof of Lemma 5.2, we can prove.  
\QED

\subsubsection{$\mf{a}(i_1, i_2, i_3)$}

For integers $i_1, i_2, i_3$ with $i_1 \geq 1$, $i_2 \geq 0$ and $i_3 \geq 1$, 
we define a subset $\mf{a}(i_1, i_2, i_3)$ of $\Mat(i_1 + i_2 + i_3, k)$ as 
\[
 \mf{a}(i_1, i_2, i_3) := 
\left\{ 
\left(
\left. 
\begin{array}{c | c | c}
 O_{i_1} & O_{i_1, i_2} &  \alpha \\
\hline 
 O_{i_2, i_1} & O_{i_2} & O_{i_2, i_3} \\
\hline 
 O_{i_3, i_1} & O_{i_3, i_2} & O_{i_3}  
\end{array}
\right)
\in \Mat(i_1 + i_2 + i_3, k)
\;\right|\;
\alpha \in \Mat_{i_1, i_3}(k)
\right\} . 
\]
We frequently use the notation $\mf{a}(i_1, i_3)$ in place of $\mf{a}(i_1, i_2, i_3)$ if 
we can understand the value of $i_1 + i_2 + i_3$ from the context. 
Clearly, $\mf{a}(i_1, i_2, i_3)$ becomes a commutative subalgebra of $\Mat(i_1 + i_2 + i_3, k)$. 
For all $N \in \mf{a}(i_1, i_2, i_3)$, we have $N^2 = O$.

\begin{lem}
Then the following assertions {\rm (1)} and {\rm (2)} hold true: 
\begin{enumerate} 
\item[\rm (1)] $\mf{a}(i_1, i_2, i_3)^{\oplus r} \subset {\cal N}_{i_1 + i_2 + i_3, r}$ $( \forall r \geq 1 )$.  
\item[\rm (2)] $\ds \mf{A}(i_1, i_2, i_3)^E \subset \bigcup_{r \geq 1} \exp({ \mf{a}(i_1, i_2, i_3) }^{\oplus r})$. 
\end{enumerate} 
\end{lem}

\Proof The proof of assertion (1) is straightforward. We shall prove assertion (2). 
Choose any element $A(T)$ of $\mf{a}(i_1, i_2, i_3)^E$. 
We can write $A(T)$ as $A(T) = \Lambda(i_1, i_2, i_3; \alpha(T) )$ for some $\alpha(T) \in \Mat_{i_1, i_3}( k[T] )$. 
Clearly, $A(T) = \Exp(A(T) - I_{i_1 + i_2 + i_3})$.
We know from Lemma 1.28 that all $A(T) - I_{i_1 + i_2 + i_3} = \sum_{i = 1}^r T^{p^{i-1}} N_i$ 
for some matrices $N_1, \ldots, N_r$ of $\mf{a}(i_1, i_2, i_3)$. 
Thus we have $A(T) = \exp(N_1, \ldots, N_r, O_n, O_n, \ldots ) \in \exp( \mf{a}(i_1, i_2, i_3)^{\oplus r})$. 
\QED

\subsubsection{${^\ell}\mf{h}^{r_1}_{r_2}(S)$}

Let $m \geq 1$ and let $(\ell, r_1, r_2) \in \Omega_m$. 
For any $S \in GL\ss(r_1, k)$, 
we let ${^\ell}\mf{h}^{r_1}_{r_2}(S)$ be the set of all matrices $N$ of $\Mat(m + 2, k)$ 
satisfying the following conditions (1), (2), (3): 
\begin{enumerate} 
\item[\rm (1)] $N = \eta(a_1, \ldots, a_{2m + 1}) - I_{m + 2}$ for some $a_1, \ldots, a_{2m + 1} \in k$. 
\item[\rm (2)] $(a_{m + 1}, \ldots, a_{m + r_1}) = (a_1, \ldots, a_{r_1}) S$. 
\item[\rm (3)] $a_i = 0$ for all integers $i$ within one of the following ranges 
\[
 \ell + 1 \leq i \leq m, \qquad m + r_1 + 1 \leq i \leq m + \ell, \qquad m + \ell + r_2 + 1 \leq i \leq 2m. 
\] 
\end{enumerate} 
Clearly, ${^\ell}\mf{h}^{r_1}_{r_2}(S)$ becomes a commutative subalgebra of $\Mat(m + 2, k)$, 
and for any $N \in {^\ell}\mf{h}^{r_1}_{r_2}(S)$ we have $I_{m + 2} + N \in \mf{H}(m + 2, k)$. 
If $p \geq 3$, for all $N, N' \in {^\ell}\mf{h}^{r_1}_{r_2}(S)$, we have $N^2 \cdot N' = O$. 
If $p = 2$, we have $N^2 = O$ for all $N \in {^\ell}\mf{h}^{r_1}_{r_2}(S)$.

\begin{lem}
The following assertions {\rm (1)} and {\rm (2)} hold true: 
\begin{enumerate} 
\item[\rm (1)] ${{^\ell}\mf{h}^{r_1}_{r_2}(S)}^{\oplus r} \subset \cN_{m + 2, r}$ $(\forall r \geq 1)$. 
\item[\rm (2)] $\ds ({^\ell}\mf{H}^{r_1}_{r_2})^E \subset \bigcup_{S \in GL\ss(r_1, k)} \exp({ {^\ell}\mf{h}^{r_1}_{r_2}(S) }^{\oplus r} ) 
\subset \mf{H}(m + 2, k[T])^E$. 
\end{enumerate} 
\end{lem}

\Proof The proof of assertion (1) is straightforward. We shall prove assertion (2). 
Choose any $A(T) = \eta(a_1, \ldots, a_{2m + 1}) \in  {^\ell}\mf{H}^{r_1}_{r_2} $. 
Then $a_1, \ldots, a_{2m + 1}$ satisfies the conditions (1), (2) and (3) in Theorem 3.3. 
So, we can express $A(T)$ as 
\begin{eqnarray*}
A(T) 
&=& 
\left( 
\prod_{1 \leq i \leq r_1}  \Exp 
\left(
\begin{array}{ c | c | c }
 0 & a_i \cdot {^t}\bm{e_i} & 0 \\ 
\hline 
 \bm{0}  & O_m     &  S \cdot (a_i \bm{e}_i) \\ 
\hline 
    0 & \bm{0} & 0 
\end{array}
\right) 
\right) 
\cdot 
\left( 
\prod_{r_1 + 1 \leq i \leq \ell} \Exp
\left(
\begin{array}{ c | c | c }
 0 & a_i \cdot {^t}\bm{e_i} & 0 \\ 
\hline 
 \bm{0}  & O_m     &  \bm{0} \\ 
\hline 
    0 & \bm{0} & 0 
\end{array}
\right) 
\right) \\
& & 
\cdot 
\left( 
\prod_{m + \ell + 1 \leq i \leq m + \ell + r_2} 
\Exp 
\left(
\begin{array}{ c | c | c }
 0 & \bm{0} & 0 \\ 
\hline 
 \bm{0}  & O_m     &  a_i  \cdot \bm{e}_{i - m } \\ 
\hline 
    0 & \bm{0} & 0 
\end{array}
\right) 
\right)
\cdot 
\left( 
\Exp 
\left(
\begin{array}{ c | c | c }
 0 & \bm{0} & \alpha(T) \\ 
\hline 
 \bm{0}  & O_m     &  \bm{0} \\ 
\hline 
    0 & \bm{0} & 0 
\end{array}
\right) 
\right) . 
\end{eqnarray*}
Write the $p$-polynomials $a_i$ ($1 \leq i \leq 2m$) and $\alpha(T)$ as 
$a_i = \sum_{j \geq 0} c_{i, j} T^{p^j}$ and $\alpha(T) = \sum_{j \geq 0} d_j T^{p^j}$. 
Let 
\begin{eqnarray*}
N_j &=& 
\left( 
\sum_{1 \leq i \leq r_1} 
\left(
\begin{array}{ c | c | c }
 0 & c_{i, j} \cdot {^t}\bm{e_i} & 0 \\ 
\hline 
 \bm{0}  & O_m     &  S \cdot (c_{i, j} \bm{e}_i) \\ 
\hline 
    0 & \bm{0} & 0 
\end{array}
\right) 
\right)
 + 
\left(
\sum_{r_1 + 1 \leq i \leq \ell} 
\left(
\begin{array}{ c | c | c }
 0 & c_{i, j} \cdot {^t}\bm{e_i} & 0 \\ 
\hline 
 \bm{0}  & O_m     &  \bm{0} \\ 
\hline 
    0 & \bm{0} & 0 
\end{array}
\right) 
\right) \\
& &
 + 
\left(
\sum_{m + \ell + 1 \leq i \leq m + \ell + r_2} 
 \left(
\begin{array}{ c | c | c }
 0 & \bm{0} & 0 \\ 
\hline 
 \bm{0}  & O_m     &  c_{i, j}  \cdot \bm{e}_{i - m } \\ 
\hline 
    0 & \bm{0} & 0 
\end{array}
\right) 
\right) 
+ 
\left(
\begin{array}{ c | c | c }
 0 & \bm{0} & d_j \\ 
\hline 
 \bm{0}  & O_m     &  \bm{0} \\ 
\hline 
    0 & \bm{0} & 0 
\end{array}
\right) 
\qquad ( \forall j \geq 0 ) .
\end{eqnarray*} 
Each $N_j \in {^\ell}\mf{h}^{r_1}_{r_2}(S)$, and there exists an integer $r \geq 1$ such that $N_j = O$ for all $j \geq r$. 
Thus we have $A(T) = \exp(N_0, \ldots, N_{r - 1}, O_n, O_n, \ldots) \in \exp({ {^\ell}\mf{h}^{r_1}_{r_2}(S) }^{\oplus r}  ) $. 
It is straightforward to prove $ \bigcup_{S \in GL\ss(r_1, k)} \exp({ {^\ell}\mf{h}^{r_1}_{r_2}(S) }^{\oplus r} ) 
\subset \mf{H}(m + 2, k[T])^E$. 
\QED

\subsubsection{$\daleth_{m + 2}^\flat$}

For any $m \geq 1$, we denote by $\daleth_{m + 2}^\flat$ the set of all matrices $N$ of $\Mat(m + 2, k)$ 
satisfying the following conditions (1) and (2): 
\begin{enumerate}
\item[\rm (1)] $N  = \eta(a_1, \ldots, a_{2m + 1}) - I_{m + 2}$ for some $a_1, \ldots, a_{2m + 1} \in k$. 
\item[\rm (2)] $\eta(a_1, \ldots, a_{2m + 1}) \in \daleth_{m + 2}$.  
\end{enumerate} 
Clearly, $\daleth_{m + 2}^\flat$ is a commutative subalgebra of $\Mat(m + 2, k)$, and 
for any $N \in \daleth_{m + 2}^\flat$ we have $I_{m + 2} + N \in \mf{H}(m + 2, k)$ and $N^2 = O_{m + 2}$.

\begin{lem}
The following assertions {\rm (1)} and {\rm (2)} hold true: 
\begin{enumerate} 
\item[\rm (1)] $(\daleth_{m + 2}^\flat)^{\oplus r} \subset \cN_{m + 2, r}$ $(\forall r \geq 1)$. 
\item[\rm (2)] $\ds (\daleth_{m + 2})^E \subset \bigcup_{r \geq 1} \exp((\daleth_{m + 2}^\flat)^{\oplus r})
\subset \mf{H}(m + 2, k[T])^E$. 
\end{enumerate} 
\end{lem}

\Proof The proof follows from Lemma 1.28 (see the proof of Lemma 5.4). 
\QED

\subsection{A classification of modular representations of elementary abelian $p$-groups into Heisenberg groups, up to equivalence}

The following theorem gives a classification of modular representations of elementary abelian $p$-groups into Heisenberg groups, 
up to equivalence.

\begin{thm}
Let $r, m \geq 1$ be integers. 
Then the following assertions {\rm (1)} and {\rm (2)} hold true: 
\begin{enumerate}
\item[\rm (1)] For any $N = (N_1, \ldots, N_r) \in \cN_{m + 2, r}$ belonging to one of the following finite dimensional algebras:  
\[
 (\daleth_{m + 2}^\flat)^{\oplus r}, \qquad \qquad 
 {{^\ell}\mf{h}^{r_1}_{r_2}(S)}^{\oplus r} \quad (\, S \in GL\ss(r_1, k) \,) , 
\]
the map $\varrho : (\Z/p\Z)^r \to GL(m + 2, k)$ defined by 
\[
\varrho([i_1], \ldots, [i_r]) := (I_{m + 2} + N_1)^{i_1} \cdots (I_{m + 2} + N_r)^{i_r}
\]
becomes a group homomorphism, and $\varrho( (\Z/p\Z)^r ) \subset \mf{H}(m + 2, k)$.

\item[\rm (2)] Let $\varrho : (\Z/p\Z)^r \to GL(m + 2, k)$ be a modular representation satisfying 
$\varrho( (\Z/p\Z)^r ) \subset \mf{H}(m + 2, k)$.  
Let $N := ( \varrho(e_1) - I_{m + 2}, \ldots, \varrho(e_r) - I_{m + 2}) \in \cN_{m + 2, r}$. 
Then there exists $N' \in \cN_{n, r}$ such that $N'$ is equivalent to $N$, and $N'$ belongs to one of the following 
finite dimensional algebras:  
\[
 (\daleth_{m + 2}^\flat)^{\oplus r}, \qquad \qquad 
 {{^\ell}\mf{h}^{r_1}_{r_2}(S)}^{\oplus r} \quad (\, S \in GL\ss(r_1, k) \,) , 
\]

\end{enumerate}
\end{thm}

\Proof Assertion (1) is clear from the definitions of $\cN_{n, r}$, $\daleth_{m +2}^\flat$ and ${^\ell}\mf{h}^{r_1}_{r_2}$. 
Assertion (2) follows from Theorem 3.1 and Lemmas 5.6, 5.5. 
\QED

\subsection{A classification of modular representations of elementary abelian $p$-groups in dimension four, up to equivalence}

In this subsection, we write down a classification of four-dimensional modular representations of 
elementary abelian $p$-groups, up to equivalence. 
Let $r \geq 1$ be an integer. 
We know in Section 5.1 that 
\[
 \Rep(E_r, 4)/\sim \; \hookrightarrow \; \E(4, k[T])/\sim
\]
So, for any $\varrho : E_r \to GL(4, k)$, we have an exponential matrix $A(T) \in \E(4, k[T])$ such that 
$A(T) = \exp(N_1, \ldots, N_r, O_4, O_4, \ldots )$, where $(N_i)_{1 \leq i \leq r} := (\varrho(e_i) - I_4 )_{1 \leq i \leq r} \in \cN_{4, r}$. 
By Corollary 4.2, there exists a regular matrix $P$ of $GL(4, k)$ such that 
\begin{eqnarray*}
P^{-1} A(T) P \in 
\left\{ 
\begin{array}{rl} 
 \mf{A}(1, 3)^E \; \cup \; \mf{A}(2, 2)^E \; \cup \; \mf{A}(3, 1)^E \; \cup \;  
 ({^2}\mf{H}^2_0)^E   & \text{if \, $p = 2$} , \\

 \mf{J}_{[3, 1]}^E \; \cup \;  
 \mf{J}_{[1, 3]}^E \; \cup \; 
 \mf{A}(1, 3)^E \; \cup \; \mf{A}(2, 2)^E \; \cup \; \mf{A}(3, 1)^E \; \cup \;  
 ({^2}\mf{H}^2_0)^E   & \text{if \, $p = 3$}, \\

 \mf{J}_{[4]}^E \; \cup \;  
 \mf{J}_{[3, 1]}^E \; \cup \;  
 \mf{J}_{[1, 3]}^E \; \cup \; 
 \mf{A}(1, 3)^E \; \cup \; \mf{A}(2, 2)^E \; \cup \; \mf{A}(3, 1)^E \; \cup \;  
 ({^2}\mf{H}^2_0)^E  & \text{if \, $p \geq 5$} . 

\end{array}
\right. 
\end{eqnarray*}
Note that $P^{-1} A(T) P = \exp(P^{-1} N_1 P, \ldots, P^{-1} N_r P, O_4, O_4, \ldots)$. 
We know from Subsection 5.2 that 
$N' := (P^{-1} N_1 P, \ldots, P^{-1} N_r P)$ satisfies 
\begin{eqnarray*}
N' \in 
\left\{ 
\begin{array}{rl} 

 \ds {\mf{a}(3, 1)}^{\oplus r}  \cup {\mf{a}(2, 2)}^{\oplus r}  \cup  {\mf{a}(1, 3)}^{\oplus r} \cup   
 \bigcup_{S \in GL\ss(2, k) } { {^2}\mf{h}^{0}_{2}(S)}^{\oplus r}  
 & \text{if \, $p = 2$} , \\

 {\mf{j}_{[3, 1]}}^{\oplus r} \cup   
 {\mf{j}_{[1, 3]}}^{\oplus r} \cup  
 \ds {\mf{a}(3, 1)}^{\oplus r}  \cup {\mf{a}(2, 2)}^{\oplus r}  \cup  {\mf{a}(1, 3)}^{\oplus r} \cup   
 \bigcup_{S \in GL\ss(2, k) } { {^2}\mf{h}^{0}_{2}(S)}^{\oplus r}  
 & \text{if \, $p = 3$}, \\
 
 {\mf{j}_{[4]}}^{\oplus r} \cup 
 {\mf{j}_{[3, 1]}}^{\oplus r} \cup   
 {\mf{j}_{[1, 3]}}^{\oplus r} \cup  
 \ds {\mf{a}(3, 1)}^{\oplus r}  \cup {\mf{a}(2, 2)}^{\oplus r}  \cup  {\mf{a}(1, 3)}^{\oplus r} \cup   
 \bigcup_{S \in GL\ss(2, k) } { {^2}\mf{h}^{0}_{2}(S)}^{\oplus r}  
  & \text{if \, $p \geq 5$} . 
\end{array}
\right. 
\end{eqnarray*} 
Now, we are ready to write down the classification of modular representations $\varrho : (\Z/p\Z)^r \to GL(4, k)$, 
up to equivalence.

\subsubsection{$p = 2$}

\begin{thm}
Assume $p = 2$. Let $r \geq 1$ be an integer and let $\varrho : (\Z/p\Z)^r \to GL(4, k)$ be a modular representation. 
Then $\varrho$ is equivalent to a representation $\varrho' :  (\Z/p\Z)^r \to GL(4, k)$ such that $\varrho'(e_i)$ $(1 \leq i \leq r)$ 
have one of the following forms {\rm (1), (2), (3), (4)}: 
\begin{enumerate}
\item[\rm (1)] 
$
\left(
\begin{array}{cccc}
 1 & \alpha_i & \beta_i & \gamma_i \\
 0 & 1 & 0 & 0 \\
 0 & 0 & 1 & 0 \\
 0 & 0 & 0 & 1
\end{array} 
\right) 
\quad (1 \leq i \leq r).
$
\item[\rm (2)] 
$
\left(
\begin{array}{cccc}
 1 & 0 & \alpha_i & \beta_i \\
 0 & 1 & \gamma_i & \delta_i \\
 0 & 0 & 1 & 0 \\
 0 & 0 & 0 & 1
\end{array} 
\right) 
\quad (1 \leq i \leq r).
$
\item[\rm (3)] 
$
\left(
\begin{array}{cccc}
 1 & 0 & 0 & \alpha_i \\
 0 & 1 & 0 & \beta_i \\
 0 & 0 & 1 & \gamma_i \\
 0 & 0 & 0 & 1
\end{array} 
\right) 
\quad (1 \leq i \leq r).
$
\item[\rm (4)] 
$
\left(
\begin{array}{cccc}
 1 & \alpha_i & \beta_i & \gamma_i \\
 0 & 1 & 0 & \mu \beta_i \\
 0 & 0 & 1 & \mu \alpha_i  \\
 0 & 0 & 0 & 1
\end{array} 
\right) 
\quad (1 \leq i \leq r) 
$ 
\quad 
having a constant $\mu \in k$ with $\mu \ne 0$. 

\end{enumerate}  
\end{thm}

\subsubsection{$p = 3$}

\begin{thm}
Assume $p = 3$. Let $r \geq 1$ be an integer and let $\varrho : (\Z/p\Z)^r \to GL(4, k)$ be a modular representation. 
Then $\varrho$ is equivalent to a representation $\varrho' :  (\Z/p\Z)^r \to GL(4, k)$ such that $\varrho'(e_i)$ $(1 \leq i \leq r)$ 
have one of the following forms {\rm (1), (2), (3), (4), (5), (6)}: 
\begin{enumerate}

\item[\rm (1)] 
$
\left(
\begin{array}{cccc}
 1 & \alpha_i & \beta_i & \gamma_i \\
 0 & 1 & \alpha_i & 0 \\
 0 & 0 & 1 & 0 \\
 0 & 0 & 0 & 1
\end{array} 
\right) 
\quad (1 \leq i \leq r). 
$
\item[\rm (2)] 
$
\left(
\begin{array}{cccc}
 1 & 0 & 0 & \gamma_i \\
 0 & 1 & \alpha_i & \beta_i \\
 0 & 0 & 1 & \alpha_i \\
 0 & 0 & 0 & 1
\end{array} 
\right) 
\quad (1 \leq i \leq r). 
$
\item[\rm (3)] 
$
\left(
\begin{array}{cccc}
 1 & \alpha_i & \beta_i & \gamma_i \\
 0 & 1 & 0 & 0 \\
 0 & 0 & 1 & 0 \\
 0 & 0 & 0 & 1
\end{array} 
\right) 
\quad (1 \leq i \leq r).
$
\item[\rm (4)] 
$
\left(
\begin{array}{cccc}
 1 & 0 & \alpha_i & \beta_i \\
 0 & 1 & \gamma_i & \delta_i \\
 0 & 0 & 1 & 0 \\
 0 & 0 & 0 & 1
\end{array} 
\right) 
\quad (1 \leq i \leq r).
$
\item[\rm (5)] 
$
\left(
\begin{array}{cccc}
 1 & 0 & 0 & \alpha_i \\
 0 & 1 & 0 & \beta_i \\
 0 & 0 & 1 & \gamma_i \\
 0 & 0 & 0 & 1
\end{array} 
\right) 
\quad (1 \leq i \leq r).
$
\item[\rm (6)] 
$
\left(
\begin{array}{cccc}
 1 & \alpha_i & \beta_i & \gamma_i \\
 0 & 1 & 0 & \lambda \alpha_i + \mu \beta_i \\
 0 & 0 & 1 & \mu \alpha_i + \nu \beta_i \\
 0 & 0 & 0 & 1
\end{array} 
\right) 
\quad (1 \leq i \leq r) 
$ 
\quad 
having constants $\lambda, \mu, \nu \in k$ with $\lambda \nu - \mu^2 \ne 0$. 

\end{enumerate}  
\end{thm}

\subsubsection{$p \geq 5$}

\begin{thm}
Assume $p = 5$. Let $r \geq 1$ be an integer and let $\varrho : (\Z/p\Z)^r \to GL(4, k)$ be a modular representation. 
Then $\varrho$ is equivalent to a representation $\varrho' :  (\Z/p\Z)^r \to GL(4, k)$ such that $\varrho'(e_i)$ $(1 \leq i \leq r)$ 
have one of the following forms {\rm (1), (2), (3), (4), (5), (6), (7)}: 
\begin{enumerate}

\item[\rm (1)] 
$
\left(
\begin{array}{cccc}
 1 & \alpha_i & \beta_i & \gamma_i \\
 0 & 1 & \alpha_i & \beta_i \\
 0 & 0 & 1 & \alpha_i \\
 0 & 0 & 0 & 1
\end{array} 
\right) 
\quad (1 \leq i \leq r). 
$
\item[\rm (2)] 
$
\left(
\begin{array}{cccc}
 1 & \alpha_i & \beta_i & \gamma_i \\
 0 & 1 & \alpha_i & 0 \\
 0 & 0 & 1 & 0 \\
 0 & 0 & 0 & 1
\end{array} 
\right) 
\quad (1 \leq i \leq r). 
$
\item[\rm (3)] 
$
\left(
\begin{array}{cccc}
 1 & 0 & 0 & \gamma_i \\
 0 & 1 & \alpha_i & \beta_i \\
 0 & 0 & 1 & \alpha_i \\
 0 & 0 & 0 & 1
\end{array} 
\right) 
\quad (1 \leq i \leq r). 
$
\item[\rm (4)] 
$
\left(
\begin{array}{cccc}
 1 & \alpha_i & \beta_i & \gamma_i \\
 0 & 1 & 0 & 0 \\
 0 & 0 & 1 & 0 \\
 0 & 0 & 0 & 1
\end{array} 
\right) 
\quad (1 \leq i \leq r).
$
\item[\rm (5)] 
$
\left(
\begin{array}{cccc}
 1 & 0 & \alpha_i & \beta_i \\
 0 & 1 & \gamma_i & \delta_i \\
 0 & 0 & 1 & 0 \\
 0 & 0 & 0 & 1
\end{array} 
\right) 
\quad (1 \leq i \leq r).
$
\item[\rm (6)] 
$
\left(
\begin{array}{cccc}
 1 & 0 & 0 & \alpha_i \\
 0 & 1 & 0 & \beta_i \\
 0 & 0 & 1 & \gamma_i \\
 0 & 0 & 0 & 1
\end{array} 
\right) 
\quad (1 \leq i \leq r).
$
\item[\rm (7)] 
$
\left(
\begin{array}{cccc}
 1 & \alpha_i & \beta_i & \gamma_i \\
 0 & 1 & 0 & \lambda \alpha_i + \mu \beta_i \\
 0 & 0 & 1 & \mu \alpha_i + \nu \beta_i \\
 0 & 0 & 0 & 1
\end{array} 
\right) 
\quad (1 \leq i \leq r) 
$ 
\quad 
having constants $\lambda, \mu, \nu \in k$ with $\lambda \nu - \mu^2 \ne 0$. 

\end{enumerate}  
\end{thm}

\vspace{1cm}

\begin{flushright}
\begin{tabular}{l}
 Faculty of Education,\\
 Shizuoka University,\\
 836 Ohya, Suruga-ku,\\
 Shizuoka 422-8529, Japan\\
e-mail: tanimoto.ryuji@shizuoka.ac.jp
\end{tabular}
\end{flushright}


\begin{thebibliography}{20}


\bibitem{Basev} 
V. A. Ba\v{s}ev, 
{\it Representations of the group $Z_2 \times Z_2$ in a field of characteristic $2$}, 
Dokl. Akad. Nauk SSSR {\bf 141}, 1961, 1015--1018. 




\bibitem{Benson}
D. J. Benson, 
{\it Representations of Elementary Abelian $p$-Groups and Vector Bundles}, 
Cambridge Tracts in Mathematics, 208. Cambridge University Press, Cambridge, 2017. 




\bibitem{Bondarenko-Drozd} 
V. M. Bondarenko, J. A. Drozd, 
{\it The representation type of finite groups} (Russian), Zap. Nau\v{c}n. Semin. Leningrad. Otdel. Mat. Inst. Steklov. (LOMI) 71 (1977) 24--41. English translation: J. Sov. Math. {\bf 20} (1982) 2515--2528.




\bibitem{Campbell-Wehlau}
 H. E. A. Campbell,  D. L. Wehlau,
{\it Modular Invariant Theory}, 
Encyclopaedia of Mathematical Sciences 139,  
Invariant Theory and Algebraic Transformation Groups VIII, Springer-Verlag, Berlin, 2011, xiv+233 pp.




\bibitem{Campbell-Shank-Wehlau}
H. E. A. Campbell, R. J. Shank, D. L. Wehlau, 
{\it Rings of invariants for modular representations of elementary abelian $p$-groups}, 
Transform. Groups {\bf 18} (2013), no. 1, 1--22. 



\bibitem{Fauntleroy}
A. Fauntleroy, 
{\it On Weitzenb\"{o}ck's theorem in positive characteristic}, 
Proc. Amer. Math. Soc. {\bf 64} (1977), no. 2, 209--213. 






\bibitem{Lucas}
E. Lucas, 
{\it Sur les congruences des nombres eul\'{e}riens et les coefficients diff\'{e}rentiels des functions trigonom\'{e}triques suivant un module premier}, 
Bull. Soc. Math. France {\bf 6} (1878), 49--54.



\bibitem{Miyanishi 1978}
M. Miyanishi, 
Curves on rational and unirational surfaces, 
Tata Institute of Fundamental Research Lectures on Mathematics and Physics, vol. 60, 
Tata Institute of Fundamental Research, Bombay, 1978. 




\bibitem{Ringel} 
C. M. Ringel, 
{\it The representation type of local algebras}, 
Proceedings of the International Conference on Representations of Algebras (Carleton Univ., Ottawa, Ont., 1974), Paper No. 22, 24 pp. Carleton Math. Lecture Notes, No. 9, Carleton Univ., Ottawa, Ont., 1974. 






\bibitem{Tanimoto 2008}
R. Tanimoto, 
{\it An algorithm for computing the kernel of a locally finite iterative higher derivation}, 
J. Pure Appl. Algebra {\bf 212} (2008), no. 10, 2284--2297.




\bibitem{Tanimoto 2013}
R. Tanimoto, 
{\it Representations of $\G_a$ of codimension two}, 
Affine algebraic geometry, 279--284, World Sci. Publ., Hackensack, NJ, 2013. 




\bibitem{Tanimoto 2018}
R. Tanimoto, 
{\it A note on the Weitzenb\"ock problem in dimension four}, 
Communications in Algebra {\bf 46} (2018), no 2, 588--596. 




\end{thebibliography}
\end{document}